%% file: final.tex
\documentclass[graybox]{svmult}

\usepackage{mathptmx}       
\usepackage{helvet}         
\usepackage{courier}        

\usepackage{makeidx}         
\usepackage{graphicx}        
\usepackage{multicol}        
\usepackage[bottom]{footmisc}

\usepackage{amsmath}
\usepackage{amssymb}

\usepackage{graphicx} 
\usepackage{array}
\usepackage{float}
\usepackage{subfig}
\usepackage{mathtools}
\usepackage{hyperref}
\usepackage{booktabs}
\usepackage{longtable}

\renewcommand{\chaptername}{}
\renewcommand{\thechapter}{}

\newcommand{\BGMId}{\mathrm{Id}}
\newcommand{\BGMlbest}{\ensuremath{\lambda_{\text{\,best}}}}
\newcommand{\BGMldefault}{\ensuremath{\lambda_{\text{\,default}}}}
\newcommand{\BGMvspace}[1]{\vspace{- #1 cm}}

\makeindex             

\begin{document}

\title{Numerical explorations of feasibility algorithms for
finding points in the intersection of finite sets}
\titlerunning{Numerical explorations of feasibility algorithms}
\author{Heinz H.\ Bauschke, Sylvain Gretchko, and  Walaa M.\
Moursi}
\institute{Heinz H. Bauschke \at Mathematics, UBCO, ASC 352, 3187 University Way, 
Kelowna, B.C.~V1V 1V7, Canada, \email{heinz.bauschke@ubc.ca}
\and Sylvain Gretchko \at Mathematics, UBCO, ASC 352, 3187 University Way, 
Kelowna, B.C.~V1V 1V7, Canada, \email{sylvain.gretchko@gmail.com},
\and Walaa M.\ Moursi \at Electrical Engineering, Stanford University, 
350 Serra Mall, Stanford, CA 94305, USA,\\
and
Mansoura University, Faculty of Science, 
Mathematics Department, Mansoura 35516, Egypt\\
 \email{wmoursi@stanford.edu}
}
%
%
\maketitle

\abstract*{
Projection methods are popular algorithms for iteratively solving
feasibility problems in Euclidean or even Hilbert spaces.
They employ (selections of) nearest point mappings to generate sequences
that are designed to approximate a point in the intersection of a collection of 
constraint sets. Theoretical properties of projection methods are fairly
well understood when the underlying constraint sets are convex; however,
convergence results for the nonconvex case are more complicated and typically
only local.
In this paper, we explore the perhaps simplest instance of
 a feasibility algorithm, namely when each constraint set consists of only finitely
 many points. We numerically investigate four constellations:
either few or many constraint sets, with either few or many points. 
Each constellation is tackled by four popular projection methods each of which
features a tuning parameter. 
We examine the behaviour 
for a single and for a multitude of orbits, and we also consider local 
and global behaviour. Our findings demonstrate the importance of the choice of the algorithm 
and that of the tuning parameter.
}
\abstract{
Projection methods are popular algorithms for iteratively solving
feasibility problems in Euclidean or even Hilbert spaces.
They employ (selections of) nearest point mappings to generate sequences
that are designed to approximate a point in the intersection of a collection of 
constraint sets. Theoretical properties of projection methods are fairly
well understood when the underlying constraint sets are convex; however,
convergence results for the nonconvex case are more complicated and typically
only local.
In this paper, we explore the perhaps simplest instance of
 a feasibility algorithm, namely when each constraint set consists of only finitely
 many points. We numerically investigate four constellations:
either few or many constraint sets, with either few or many points. 
Each constellation is tackled by four popular projection methods each of which
features a tuning parameter. 
We examine the behaviour 
for a single and for a multitude of orbits, and we also consider local 
and global behaviour. Our findings demonstrate the importance of the choice of the algorithm 
and that of the tuning parameter.
}

\begin{keywords} 
cyclic Douglas--Rachford algorithm, 
Douglas--Rachford algorithm,
extrapolated parallel projection method,
method of cyclic projections, 
nonconvex feasibility problem, 
optimization algorithm, 
projection
\end{keywords}

\noindent{\bf AMS 2010 Subject Classification:} 
49M20, 49M27, 49M37, 65K05, 65K10, 90C25, 90C26, 90C30

\section{Introduction}

\runinhead{Background}

Let $X$ be a Euclidean space (i.e., a finite-dimensional Hilbert
space), with inner product $\langle\cdot,\cdot\rangle$ and norm $\|\cdot\|$. 
The \emph{feasibility problem} is a common 
problem in science and
engineering: given finitely many 
closed subsets $C_1,\ldots,C_m$ of $X$, 
it asks to 
\begin{equation}
\label{fp}
\text{Find $x\in C := C_1\cap \cdots \cap C_m$.} 
\tag{FP}
\end{equation}
We henceforth assume that the intersection $C$ is nonempty. 
Algorithms for solving \eqref{fp} exist when the constraint sets $C_i$ allow for
simple projectors $P_{C_i}$ (i.e., nearest-point mappings). When $C_i$ is convex,
then the projector $P_{C_i}$
is a nice (firmly nonexpansive and single-valued) operator defined on the entire
space $X$; when $C_i$ is not convex, then $P_{C_i}$ is nonempty and set-valued. 
For notational simplicity, we will
use $P_{C_i}$ to denote an arbitrary but fixed selection of the set-valued projector. 
(If $S$ is a subset of $X$, then $P_S(x)$ is a minimizer of the function $s\mapsto \|x-s\|$,
where $s\in S$. For other notions not explicitly defined in this paper, we refer the reader to 
\cite{BGMbib-BBor}.)

Assuming that the operators $P_{C_1},\ldots,P_{C_m}$ are readily available and
implementable, one may try to solve \eqref{fp} iteratively by generating a sequence
$(x_k)_{k\in \mathbb{N}}$ 
of vectors in $X$ that employs the projection operators $P_{C_i}$ in some fashion to produce
the next update. 
There are hundreds of papers dealing with algorithms for solving
convex or nonconvex feasibility problems. Thus, we refrain from providing
a comprehensive list of references and rather point to 
the following recent books and ``meta'' papers as starting points:
\cite{BGMbib-BBor,BGMbib-BC,BGMbib-BK,BGMbib-Ceg,BGMbib-CZak,BGMbib-CZ,BGMbib-Com97a,BGMbib-Com97b}.
(We note that the recent manuscript \cite{BGMbib-BDL} deals with a feasibility problem
where one set is a doubleton.)
The convergence theory in the nonconvex case is much more challenging
and usually of local character.

\runinhead{Goal of this paper}

\emph{The goal of this paper is to showcase the surprising numerical complexity 
of the most simple instance of \eqref{fp}}; namely,
when each constraint set 
\begin{svgraybox}
\begin{equation*}
\text{$C_i$ contains a 
\emph{finite} number of points.}
\end{equation*}
\end{svgraybox}
In this case, the projection operator is very easy to implement --- 
this is achieved by 
simply measuring the distance of the point to each point in $C_i$ and returning
the closest one. 
Furthermore, we will restrict ourselves to the simple case when 
the underlying space
\begin{svgraybox}
\begin{equation*}
    \text{$X=\mathbb{R}^2$}
\end{equation*}
\end{svgraybox}
is simply the \emph{Euclidean plane}. Even in this setting,
the difficulty and richness of the dynamic behaviour is
impressively illustrated. 

\emph{It is our hope that the complexity revealed will spark further
analytical research in feasibility algorithms with the goal
to explain the observed complexity and ultimately to aid 
in the design of new algorithms for solving difficult feasibility problems.}

\runinhead{Organization of the paper}

The remainder of the paper is organized as follows.
In Section~\ref{BGMsec:4con}, we present the four constellations
we will use for our numerical exploration throughout the remainder of
the paper. These constellations correspond to feasibility problems
that we will attempt to solve using the algorithms listed in
Section~\ref{BGMsec:4alg}. Section~\ref{BGMsec:setup}
provides details on the implementation and execution 
of the numerical experiments.
The ``best'' tuning parameter $\BGMlbest$ is determined
in Section~\ref{BGMsec:det}. We then track
typical orbits of the algorithms in Section~\ref{BGMsec:track}.
Local and global behaviour is investigated in Section~\ref{BGMsec:locglob}.
Some interesting (and beautiful) behaviour outside of 
the main numerical experiments are collected in 
Section~\ref{BGMsec:divert}.
The final Section~\ref{BGMsec:last} contains some concluding remarks. 

\section{The four constellations}

\label{BGMsec:4con}

Even though we restrict ourselves already to 
finitely many constraint sets with finitely many points in the Euclidean plane,
the infinitely many possibilities to experiment make it a daunting task to
explore this space. 
We opted to probe this universe as follows. 

\begin{figure}
	\begin{tabular}{ccc}
		\subfloat[Few sets with few points]{\includegraphics[scale=0.33]{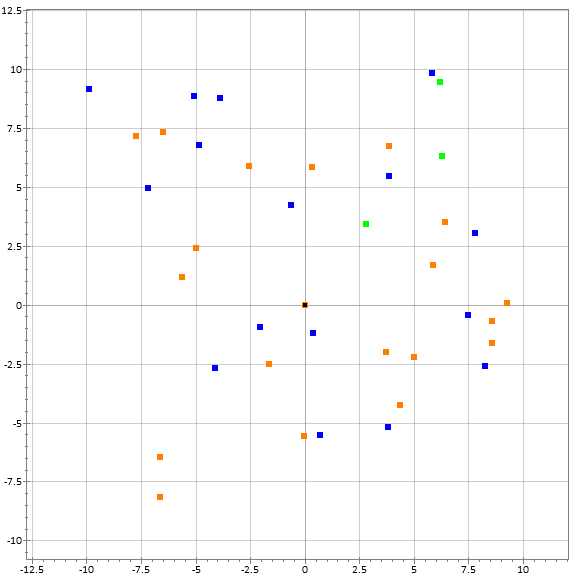}} &~~~~~~~~~&
		\subfloat[Few sets with many points]{\includegraphics[scale=0.33]{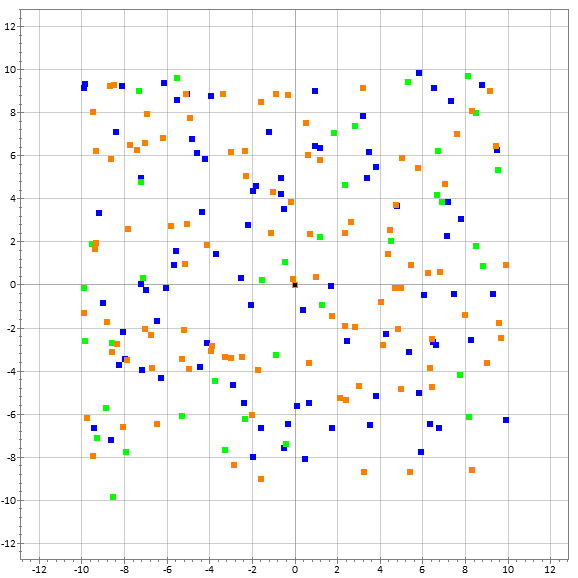}}\\
		\subfloat[Many sets with few points]{\includegraphics[scale=0.33]{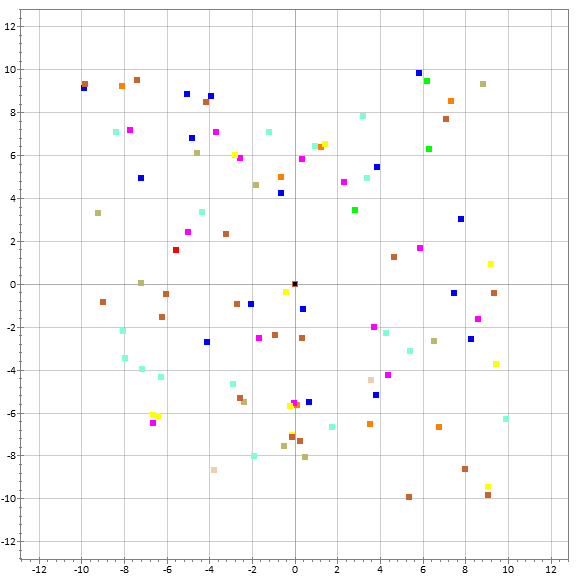}} & &
		\subfloat[Many sets with many points]{\includegraphics[scale=0.33]{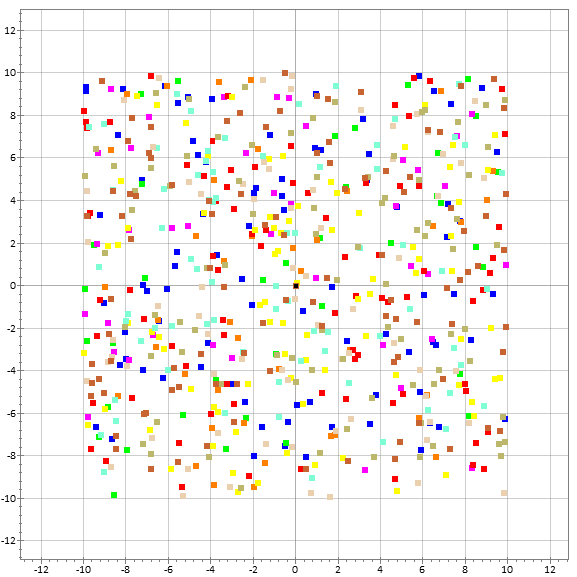}}
	\end{tabular}
	\caption{The four constellations explored in this paper. See Section~\ref{BGMsec:4con} for further information.}
	\label{fig:BGM_All_Four_Constellations}
\end{figure}

The points in each constraint set $C_i$ are chosen randomly.
We will ensure that the origin belongs to each set $C_i$ 
\begin{svgraybox}
\begin{equation*}
0 \in C_i \subset [-10,10]\times[-10,10]
\end{equation*}
\end{svgraybox}
to have a consistent feasibly problem with
\begin{svgraybox}
\begin{equation*}
C = C_1\cap \cdots \cap C_m = \{0\}.
\end{equation*}
\end{svgraybox}
We will focus on two alternatives for the 
\emph{number of constraint sets}, either ``few'' or ``many''.
We will also consider constraint sets
with a \emph{maximum number of points} in the constraint sets, 
either ``few'' or ``many''.
From now on, we will use the following language:
\begin{svgraybox}
\begin{itemize}
    \item The number of \textbf{few sets} is $3$.
    \item The number of \textbf{many sets} is $10$.
    \item The number of \textbf{few points} is $20$.
    \item The number of \textbf{many points} is $100$. 
\end{itemize}
\end{svgraybox}
This will give rise to \emph{four constellations}:
few sets with few points, few sets with many points,
many sets with few points, and many sets with many points.
The four constellations used in our numerical experiments are shown in
Figure~\ref{fig:BGM_All_Four_Constellations}.


\section{The four feasibility algorithms}

\label{BGMsec:4alg}

We will numerically solve instances of \eqref{fp} using four 
algorithms which
we briefly review in this section. 
While there is a myriad of competing algorithms available, our selection consists of trustworthy
``work horses'' that have been employed elsewhere and for which
the convergence theory in the \emph{convex} case is fairly well understood. 
Each of these algorithms has a ``tuning'' parameter $\lambda$ in the range
$\left]0,2\right[$.
The \emph{default value} 
$\BGMldefault$ is $1$. Guided by experiments, we will also (numerically) look for the ``best'' 
value $\BGMlbest$. We now turn to these four algorithms. Each algorithm
will have a \emph{governing sequence} driving the iteration, and a (possibly different)
\emph{monitored sequence} which is meant to find a solution of \eqref{fp}.

\runinhead{Cyclic Projections (CycP)}

Given $x_0\in X$, the governing sequence is defined by 
\begin{svgraybox}
\begin{equation}
x_{k+1} \coloneqq \big((1-\lambda)\BGMId+\lambda P_{C_m}\big)\circ \cdots \circ
\big((1-\lambda)\BGMId+\lambda P_{C_1}\big)x_{k}. 
\end{equation}
\end{svgraybox}
The default parameter is $\BGMldefault=1$, from the range $\left]0,2\right[$. 
The sequence monitored is $\big(\tfrac{1}{m}\sum_{i=1}^mP_{C_i}x_{k}\big)_{k\in\mathbb{N}}$. 
Selected references: \cite{BGMbib-BBor, BGMbib-BC, BGMbib-BK,BGMbib-Ceg,BGMbib-CCCDH, BGMbib-CZ,
BGMbib-Com97b}. 

\runinhead{Extrapolated Parallel Projections (ExParP)}

Given $x_0\in X$, the governing and monitored sequence is defined by 
\begin{svgraybox}
\begin{equation}
x_{k+1} \coloneqq x_k + \lambda\cdot \frac{\sum_{i=1}^m\|x_k-P_{C_i}x_k\|^2}{\|\sum_{i=1}^m(x_k-P_{C_i}x_k)\|^2}
\sum_{i=1}^m(P_{C_i}x_k-x_k)
\end{equation}
\end{svgraybox}
if $x_k\notin C$; $x_{k+1} = x_k$ otherwise. 
The default parameter is $\BGMldefault=1$, from the range $\left]0,2\right[$. 
Selected references: \cite{BGMbib-BC,BGMbib-BCK,BGMbib-Com97a}.

\runinhead{Douglas--Rachford (DR)}

Given $x_0\in X$, $\mathbf{x}_0 \coloneqq (x_{0,1},\ldots,x_{0,m}) = (x_0,\ldots,x_0) 
\in \mathbf{X} \coloneqq X^m$, 
$\mathbf{x}_k = (x_{k,1},\ldots,x_{k,m})\in \mathbf{X}$, 
and
$\bar{x}_k \coloneqq \tfrac{1}{m}\sum_{i=1}^mx_{k,i}$, 
the next iterate is 
$\mathbf{x}_{k+1} = (x_{k+1,1},\ldots,x_{k+1,m})$, where
\begin{svgraybox}
\begin{equation}
(\forall i\in\{1,\ldots,m\})\quad
x_{k+1,i} \coloneqq x_{k,i} +\lambda\big(P_{C_i}(2\bar{x}_k-x_{k,i})-\bar{x}_k\big).
\end{equation}
\end{svgraybox}
The default parameter is $\BGMldefault=1$, from the range $\left]0,2\right[$. 
The sequence monitored is $(\bar{x}_k)_{k\in\mathbb{N}}$. 
Selected references: \cite{BGMbib-BC,BGMbib-BM,BGMbib-EckBer,BGMbib-ERT,BGMbib-LM}. 

\runinhead{Cyclic Douglas--Rachford (CycDR)}

Given $x_0\in X$, the governing sequence is defined by 
\begin{svgraybox}
\begin{multline}
x_{k+1} \coloneqq \big((1-\tfrac{\lambda}{2})P_{C_m}+
\tfrac{\lambda}{4}\big(\BGMId + R_{C_1}R_{C_{m}}\big)\big)\circ \cdots \circ \\
\big((1-\tfrac{\lambda}{2})P_{C_2}+\tfrac{\lambda}{4}\big(\BGMId + R_{C_3}R_{C_{2}}\big)\big)\circ
\big((1-\tfrac{\lambda}{2})P_{C_1}+\tfrac{\lambda}{4}\big(\BGMId + R_{C_2}R_{C_{1}}\big)\big)x_{k}.
\end{multline}
\end{svgraybox}
The default parameter is $\BGMldefault=1$, from the range $\left]0,2\right[$. 
The sequence monitored is $\big(\tfrac{1}{m}\sum_{i=1}^mP_{C_i}x_{k}\big)_{k\in\mathbb{N}}$. 
Selected references: \cite{BGMbib-Luke,BGMbib-DP18a,BGMbib-LST}.
(For other cyclic version of DR, see \cite{BGMbib-BNP,BGMbib-BT}.
Also, if $m=2$ and $C_1=C_2=\{0\}$, then $(x_k)_{k\in\mathbb{N}} =
((\lambda/2)^kx_0)_{k\in\mathbb{N}}$ is actually unbounded when 
$x_0\neq 0$ and $\lambda>2$.) 

\section{Setting up the numerical explorations}

\label{BGMsec:setup}

\runinhead{Stopping criteria}
The feasibility measure 
\begin{equation*}
d: X \to \mathbb{R}_+: x \mapsto \sqrt{\frac{\sum_{i=1}^m \lVert x - P_{C_i}x\rVert^2 }{\sum_{i=1}^m \lVert x_0 - P_{C_i}x_0\rVert^2}},
\end{equation*}
where $x_0\in X\smallsetminus C$, vanishes exactly when 
$x\in C$. 
We stop each algorithm with monitored sequence 
$(y_k)_{k\in \mathbb{N}}$ 
either when
\begin{equation*}
d(y_k) < \epsilon := 10^{-6}
\end{equation*}
or when the maximum number of iterations, which we set to 1000, is reached.
These values were chosen to allow a reasonable exploration of the feasibility problem 
while maintaining computational efficiency. 

\begin{figure}
	\centering
		\begin{tabular}{c}
			\subfloat[]{\includegraphics[scale=0.3]{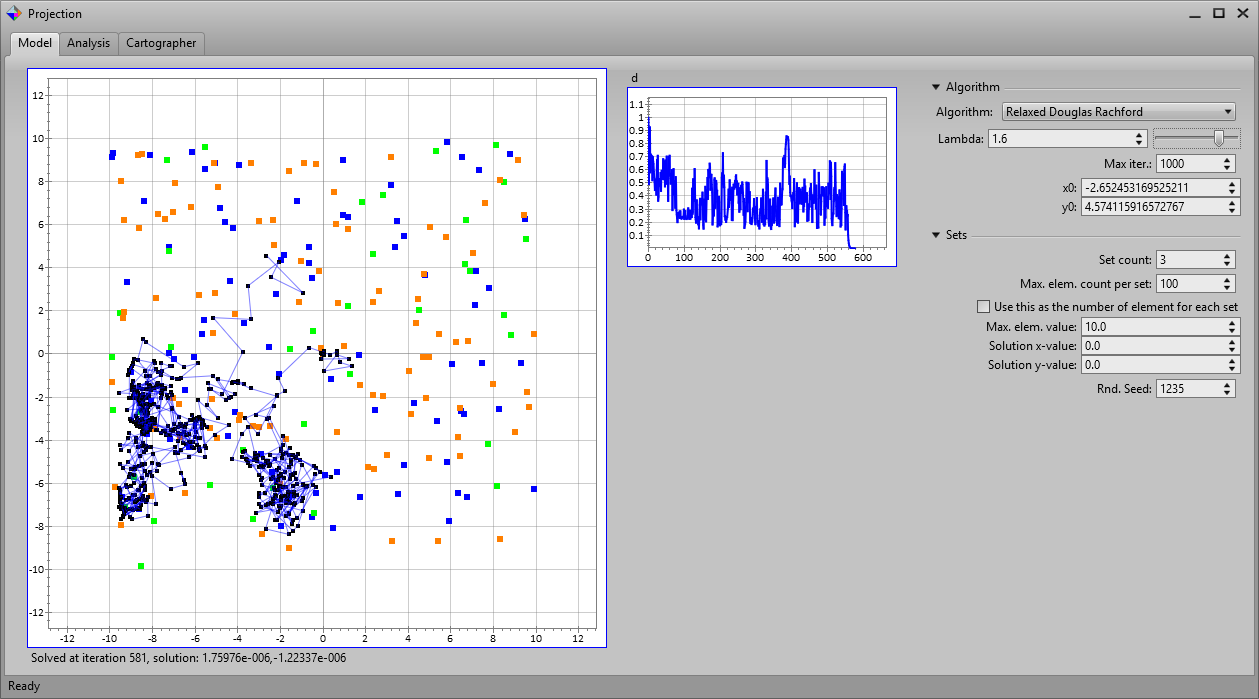}} \\
			\subfloat[]{\includegraphics[scale=0.3]{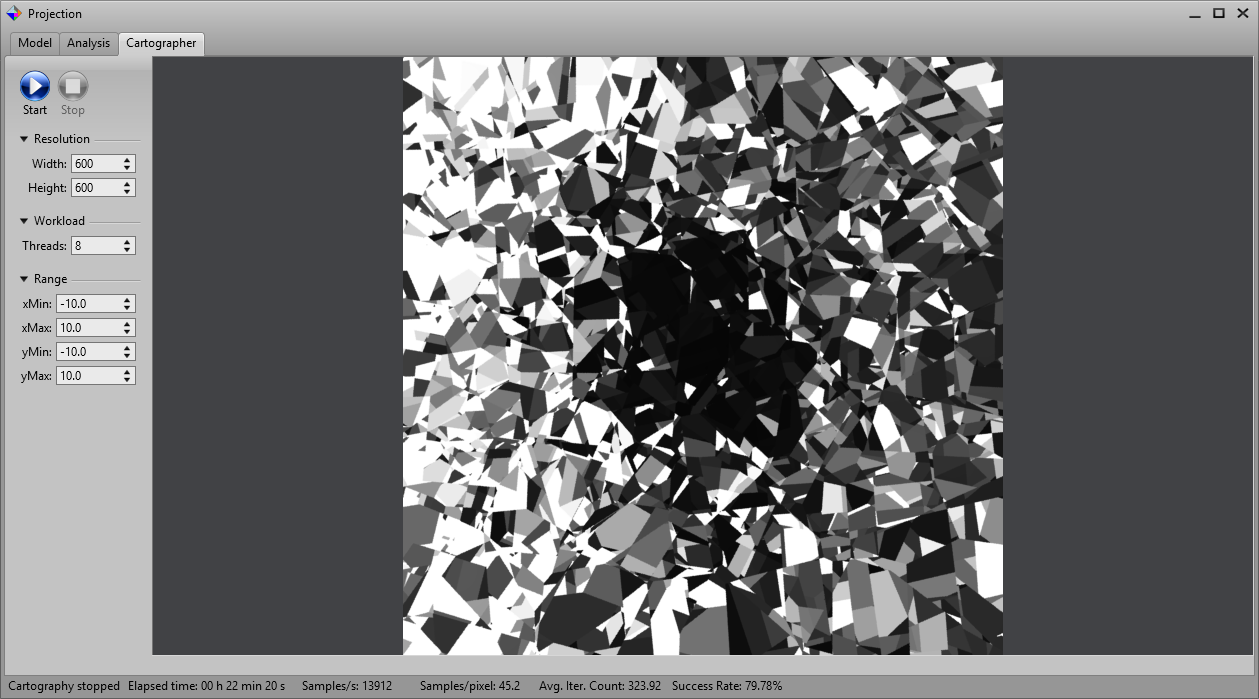}}
		\end{tabular}
		\caption{The software developed for this work. Setting the constellation of points and the algorithm to be used is done in the main tab, shown in (a). The generation of a performance plot is done in the cartographer tab, shown in (b).}
		\label{fig:BGM_screenshots}
\end{figure}

\runinhead{Details on program}
A program was developed in C++ to run the different experiments, 
see Figure~\ref{fig:BGM_screenshots} for two screenshots which 
we describe next. In the main
tab of the user interface one can select the algorithm to be used and set up
the problem to be solved by choosing the number of sets and the maximum
number of elements per set. By clicking on the diagram showing the current
constellation of points, the user can select a starting point and immediately
observe the resulting orbit being rendered over the constellation. The graph
of the feasibility measure $d$, corresponding to the current orbit, is also
displayed.\\
The Cartographer tab allows the exploration of a very large number of
starting points to construct a picture of the performance of a given
algorithm. This two-dimensional plot shows for each starting point the number
of iterations required to solve the problem, ranging from zero (black) to the
maximum number of iterations allowed (white). The plot is generated
progressively and uses Quasi-Monte Carlo sampling for the selection of the
starting points. This is the most computationally intensive part of the
software, and it is fully multi-threaded to take advantage of modern
processor architectures.

\vspace{-1.0cm}

\section{Determining the ``best'' parameter $\BGMlbest$}
\label{BGMsec:det}

\begin{figure}[H]
	\begin{tabular}{cccc}
		\subfloat[CycP]{\includegraphics[scale=0.13]{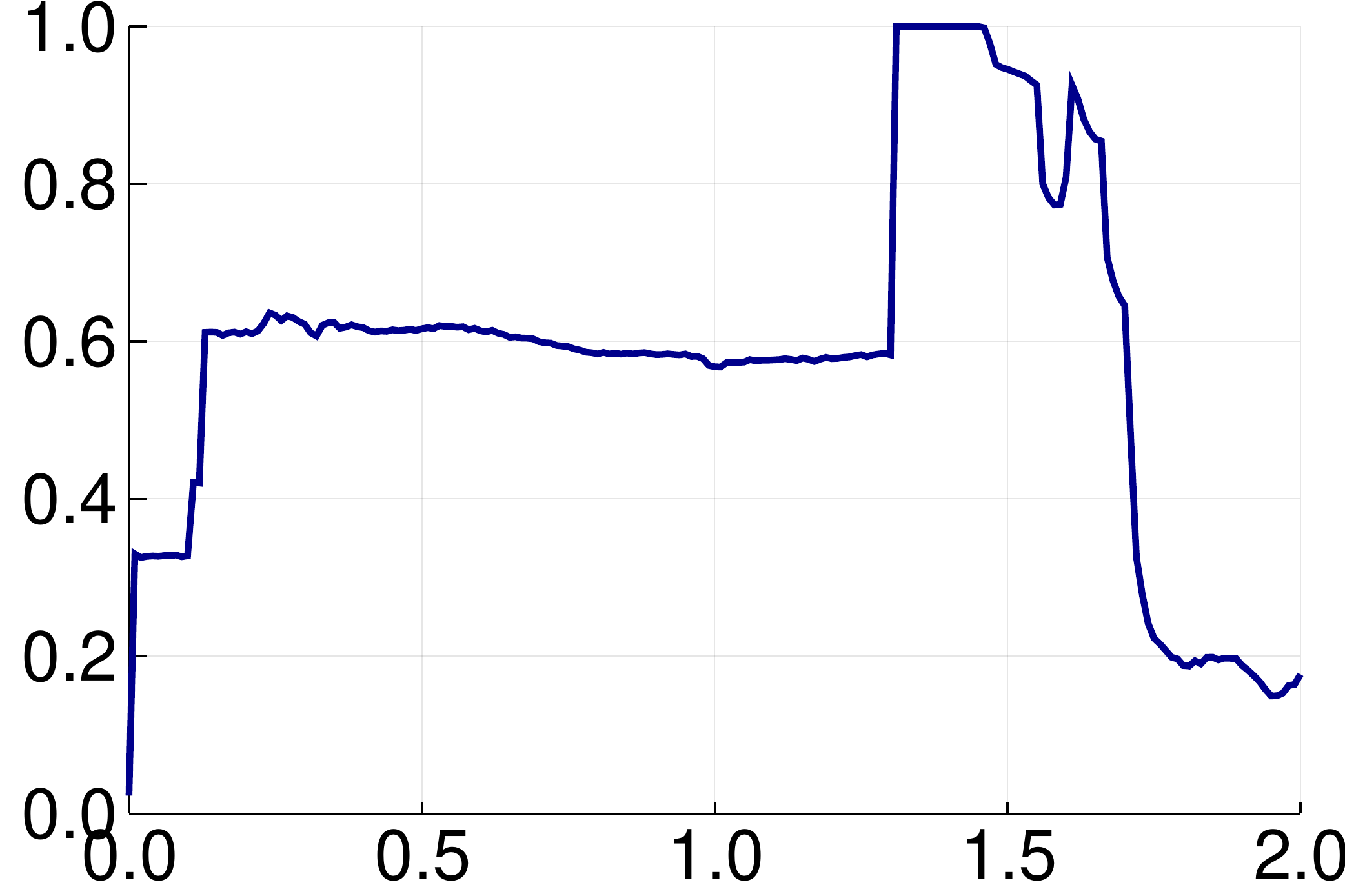}} & 
		\subfloat[ExParP]{\includegraphics[scale=0.13]{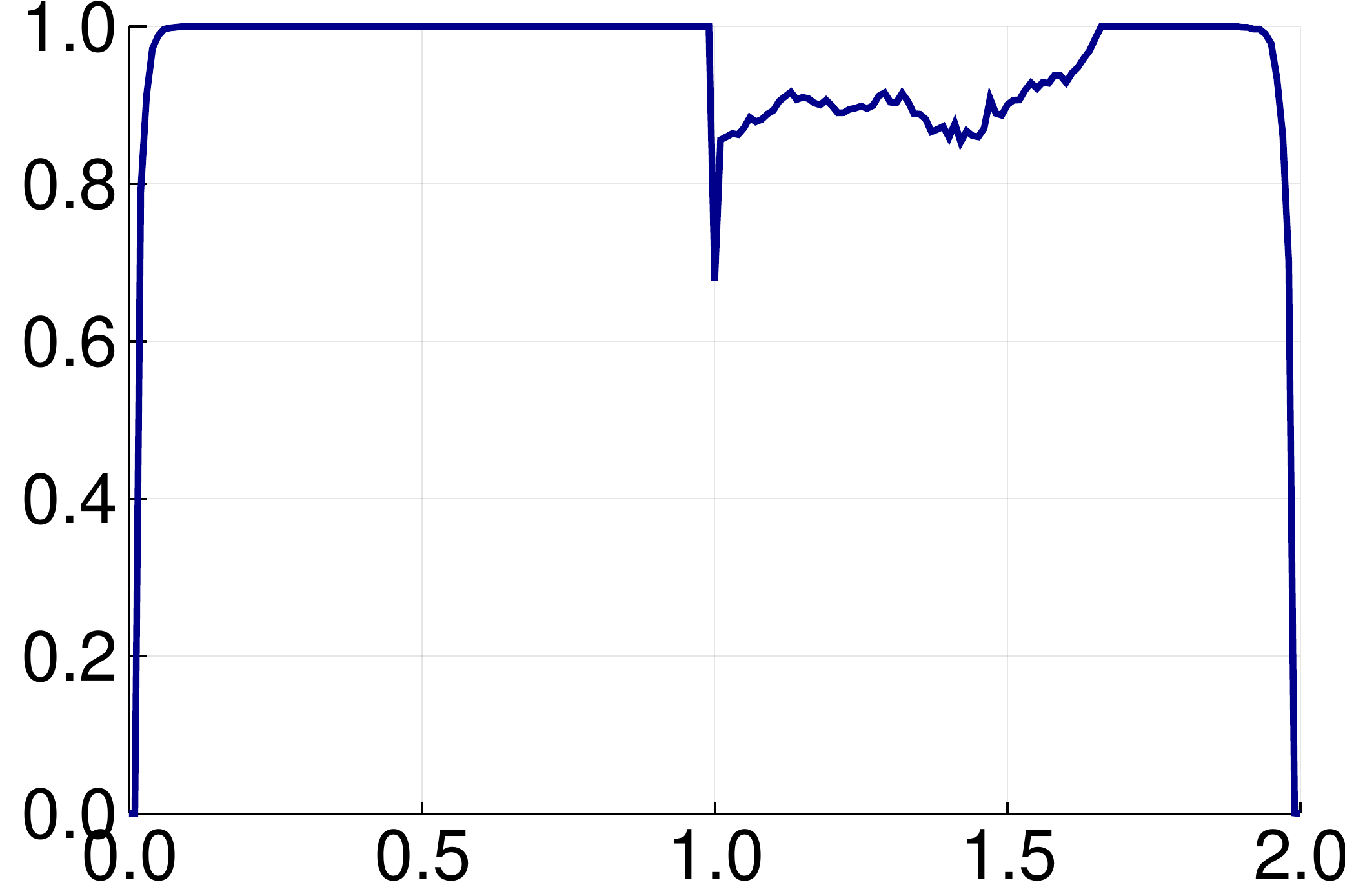}}
		& \subfloat[DR]{\includegraphics[scale=0.13]{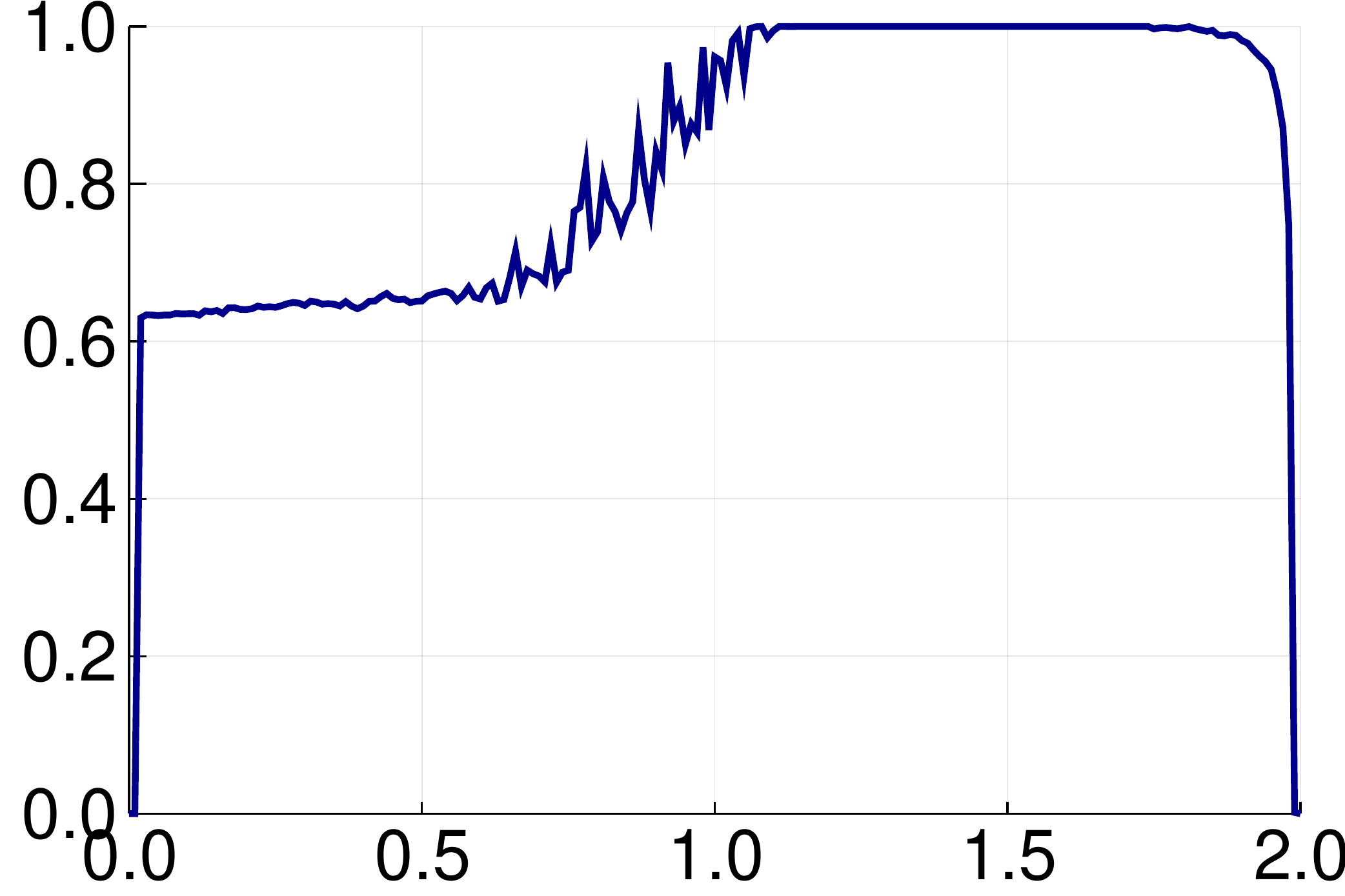}} & 
		\subfloat[CycDR]{\includegraphics[scale=0.13]{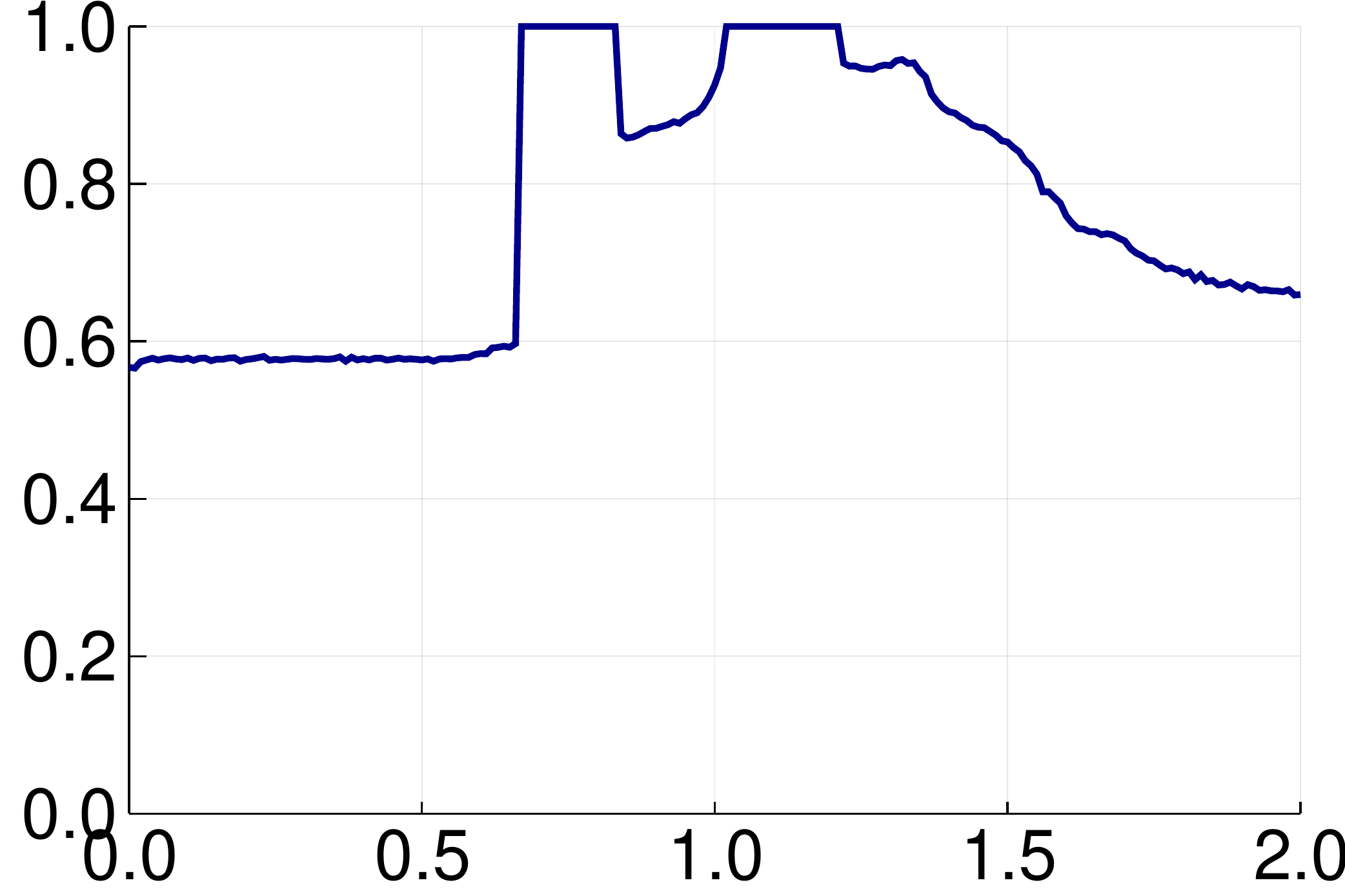}}
	\end{tabular}
	\caption{Success rates in terms of $\lambda$
		for the few sets with few points constellation.}
	\label{fig:BGM_FewSets_FewPts_LambdaCurves}
\end{figure}

\begin{figure}[H]
	\begin{tabular}{cccc}
		\subfloat[CycP]{\includegraphics[scale=0.13]{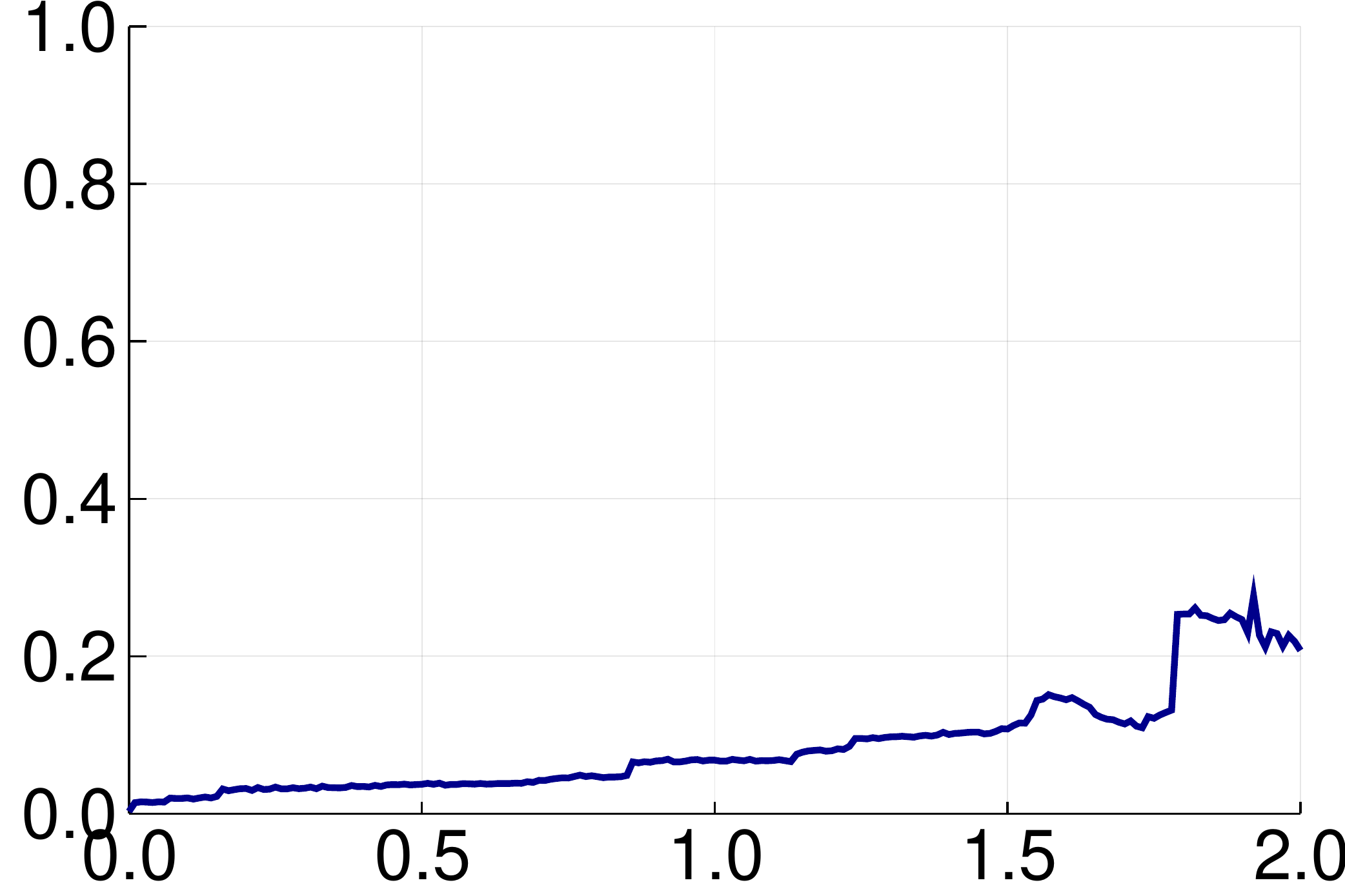}} & 
		\subfloat[ExParP]{\includegraphics[scale=0.13]{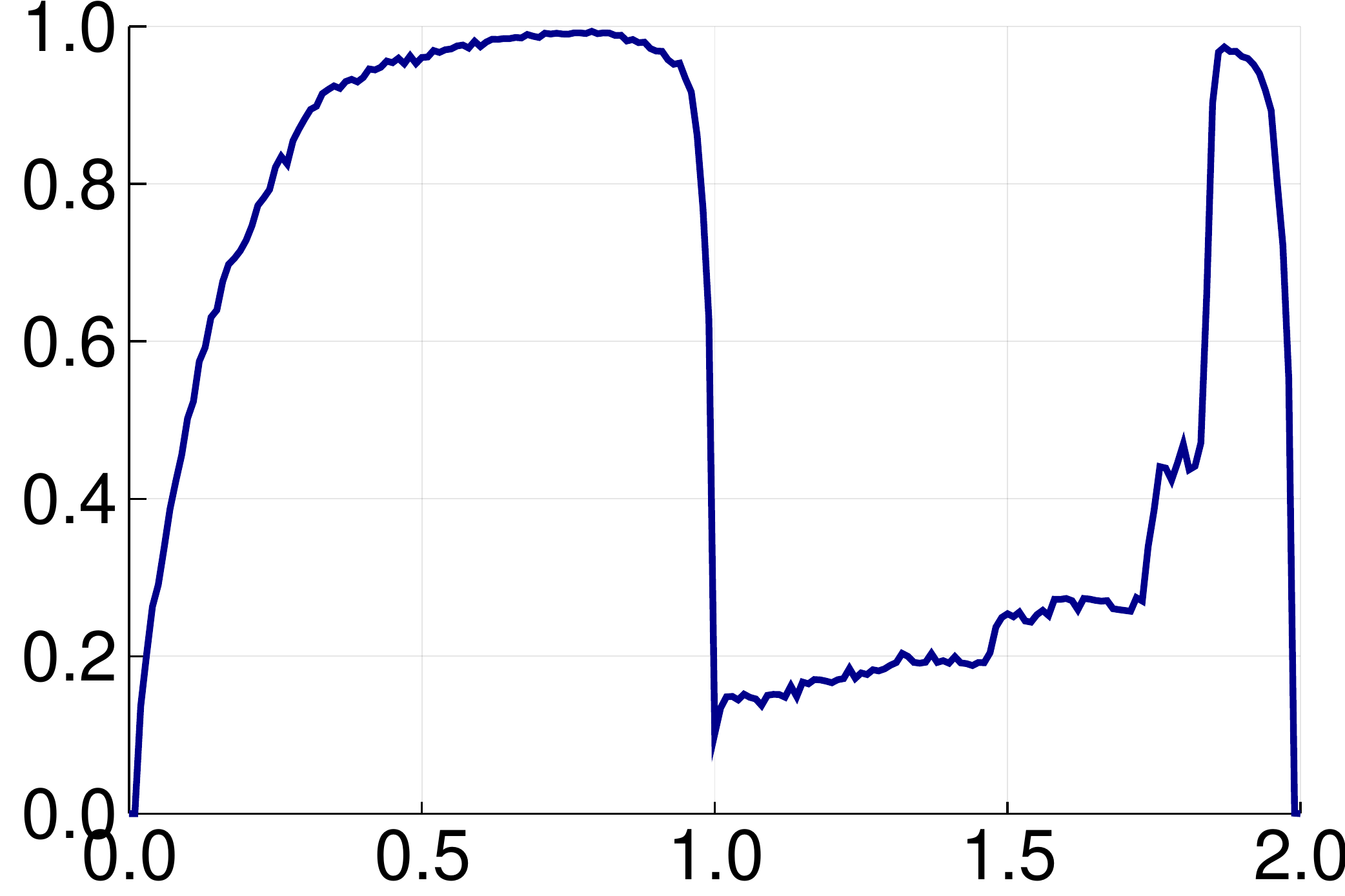}}
		& \subfloat[DR]{\includegraphics[scale=0.13]{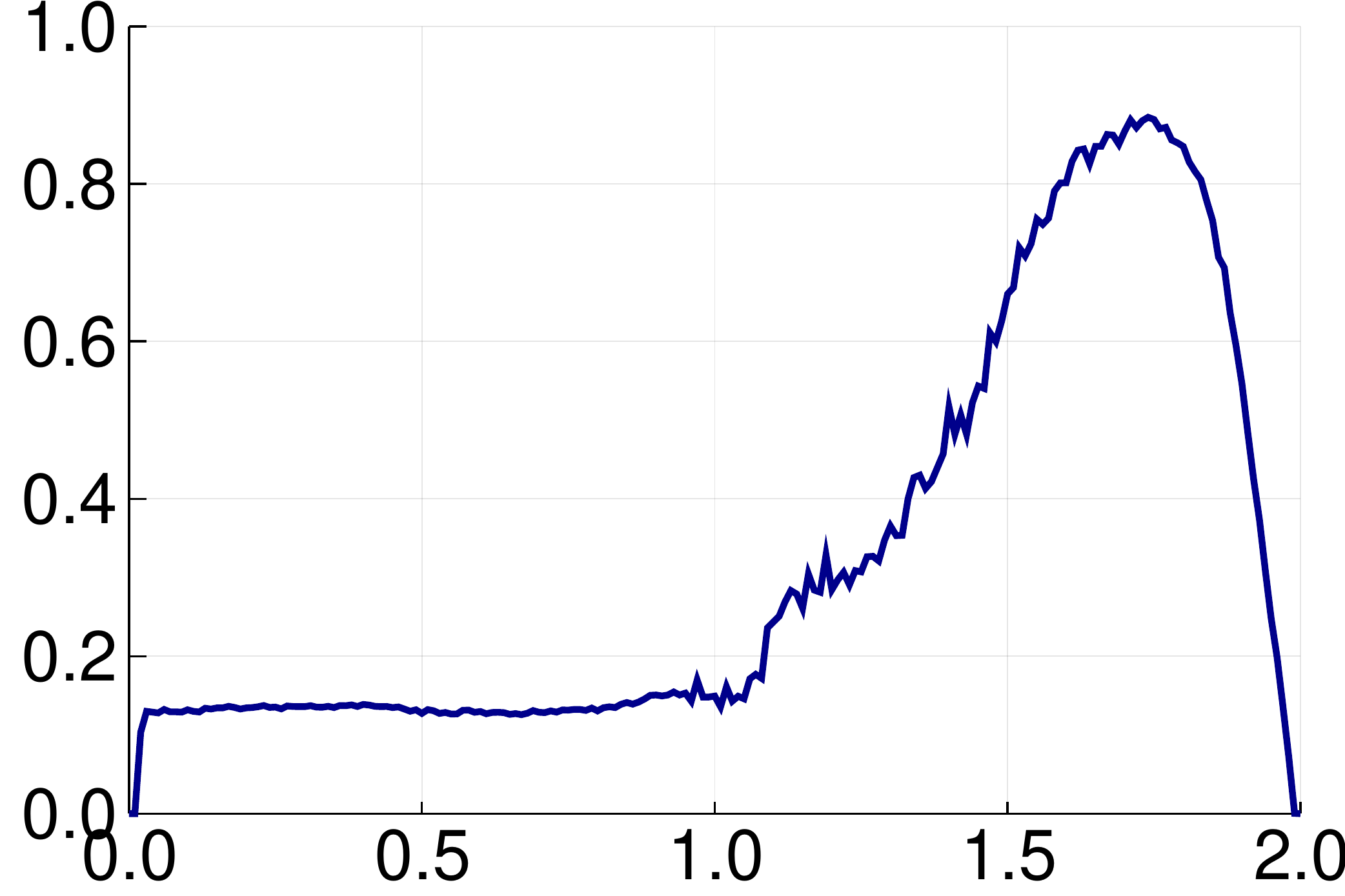}} & 
		\subfloat[CycDR]{\includegraphics[scale=0.13]{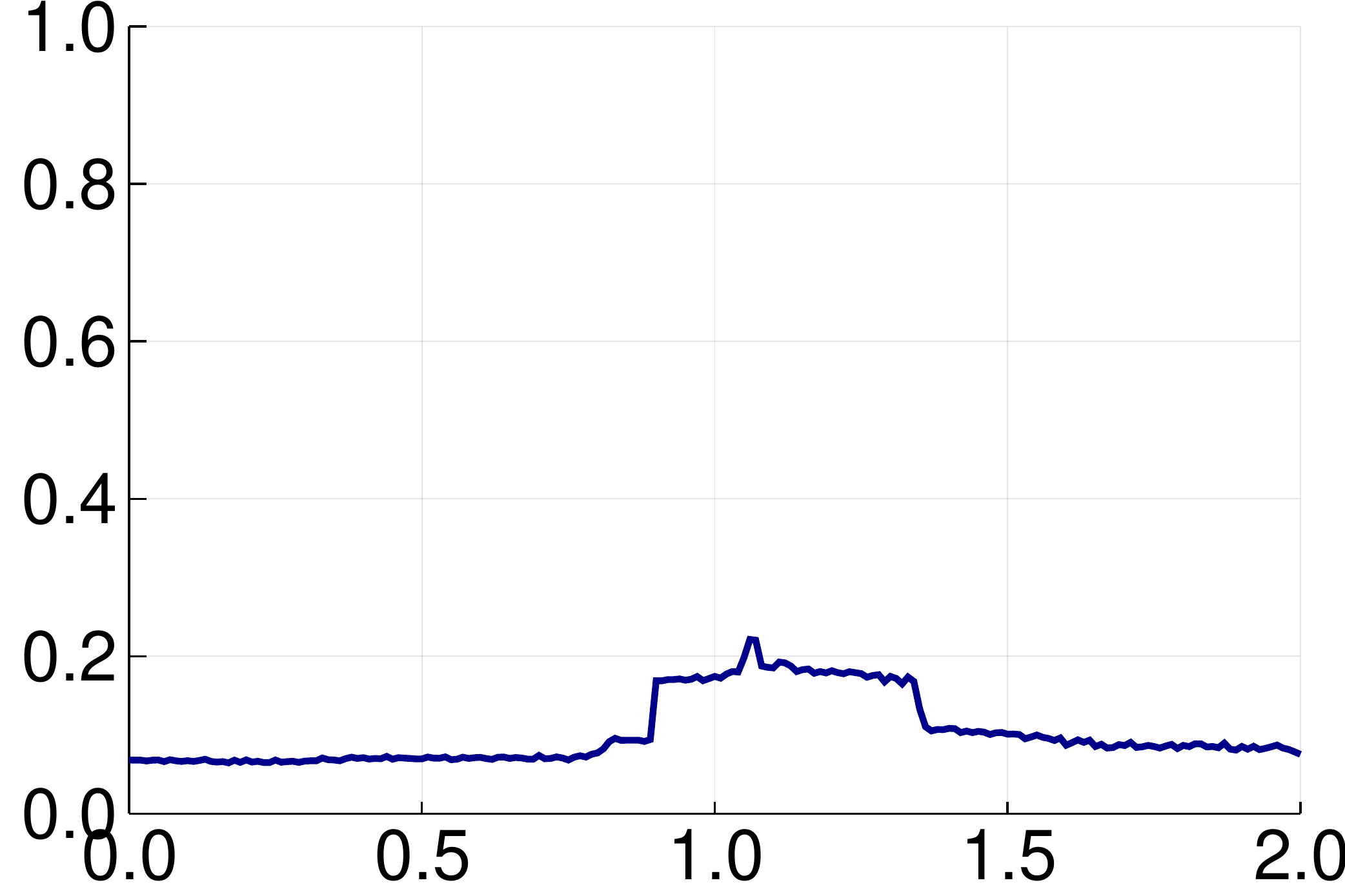}}
	\end{tabular}
	\caption{Success rates in terms of $\lambda$
		for the few sets with many points constellation.}
	\label{fig:BGM_FewSets_ManyPts_LambdaCurves}
\end{figure}

\begin{figure}[H]
	\begin{tabular}{cccc}
		\subfloat[CycP]{\includegraphics[scale=0.13]{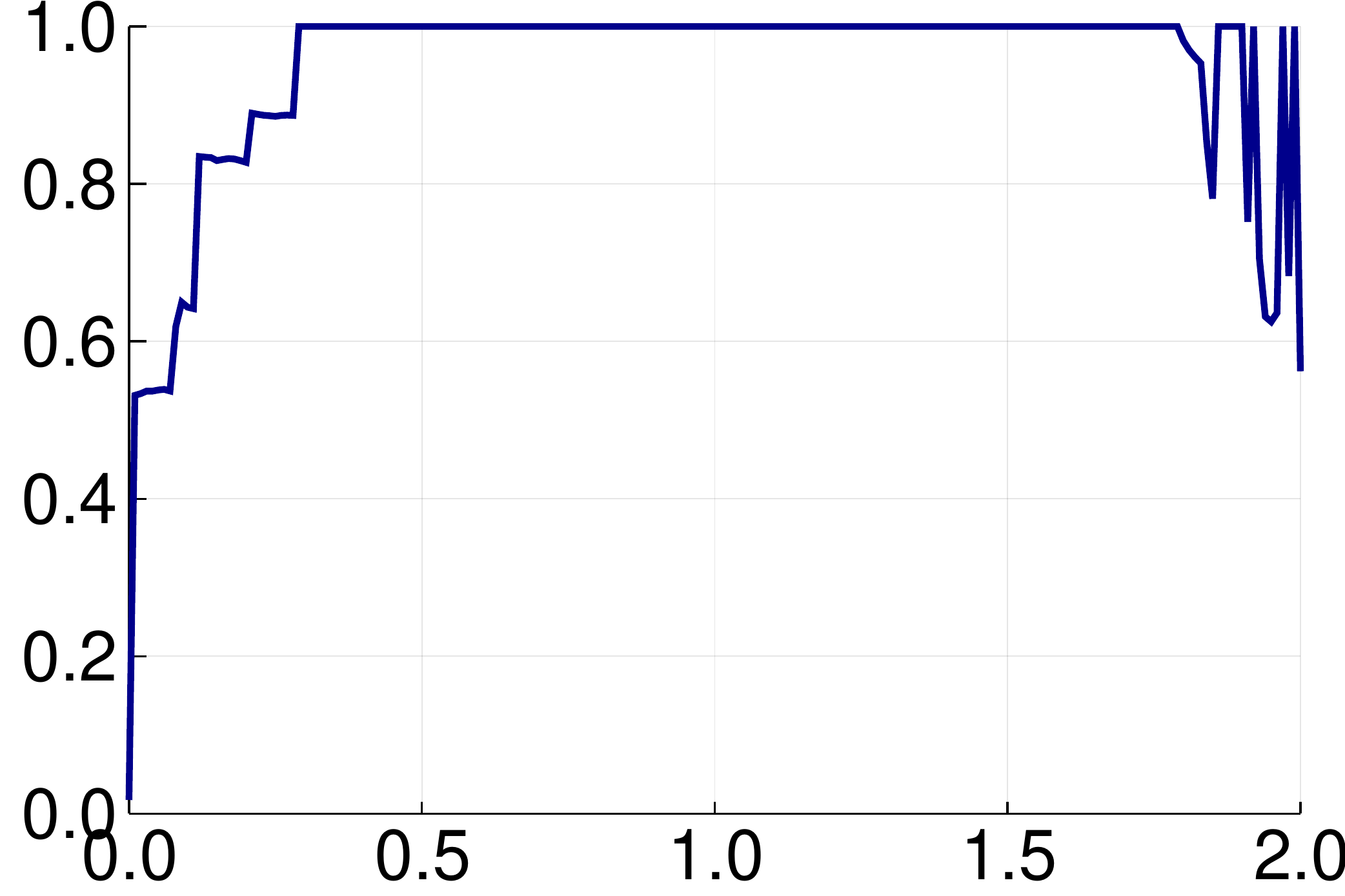}} & 
		\subfloat[ExParP]{\includegraphics[scale=0.13]{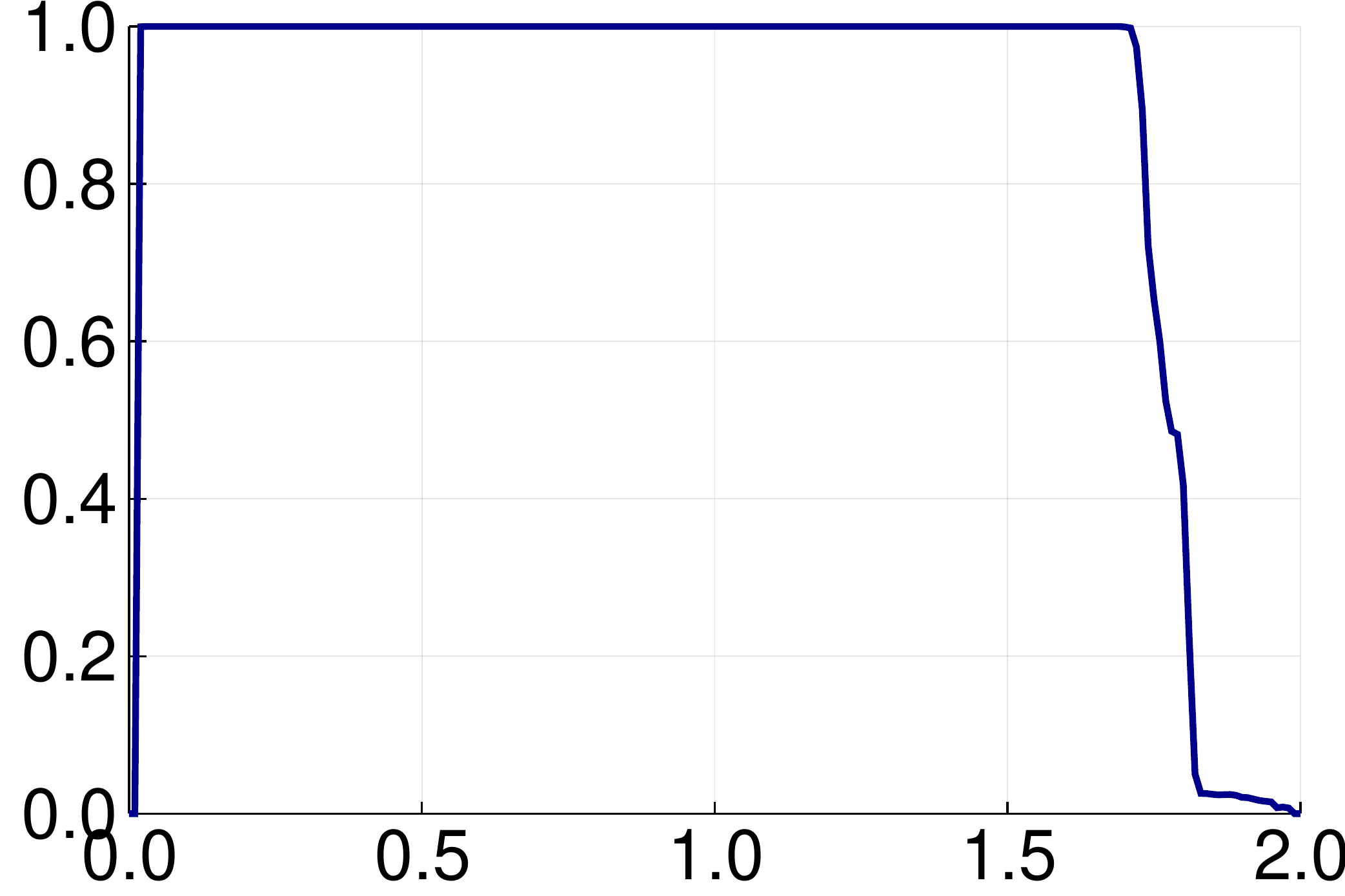}}
		& \subfloat[DR]{\includegraphics[scale=0.13]{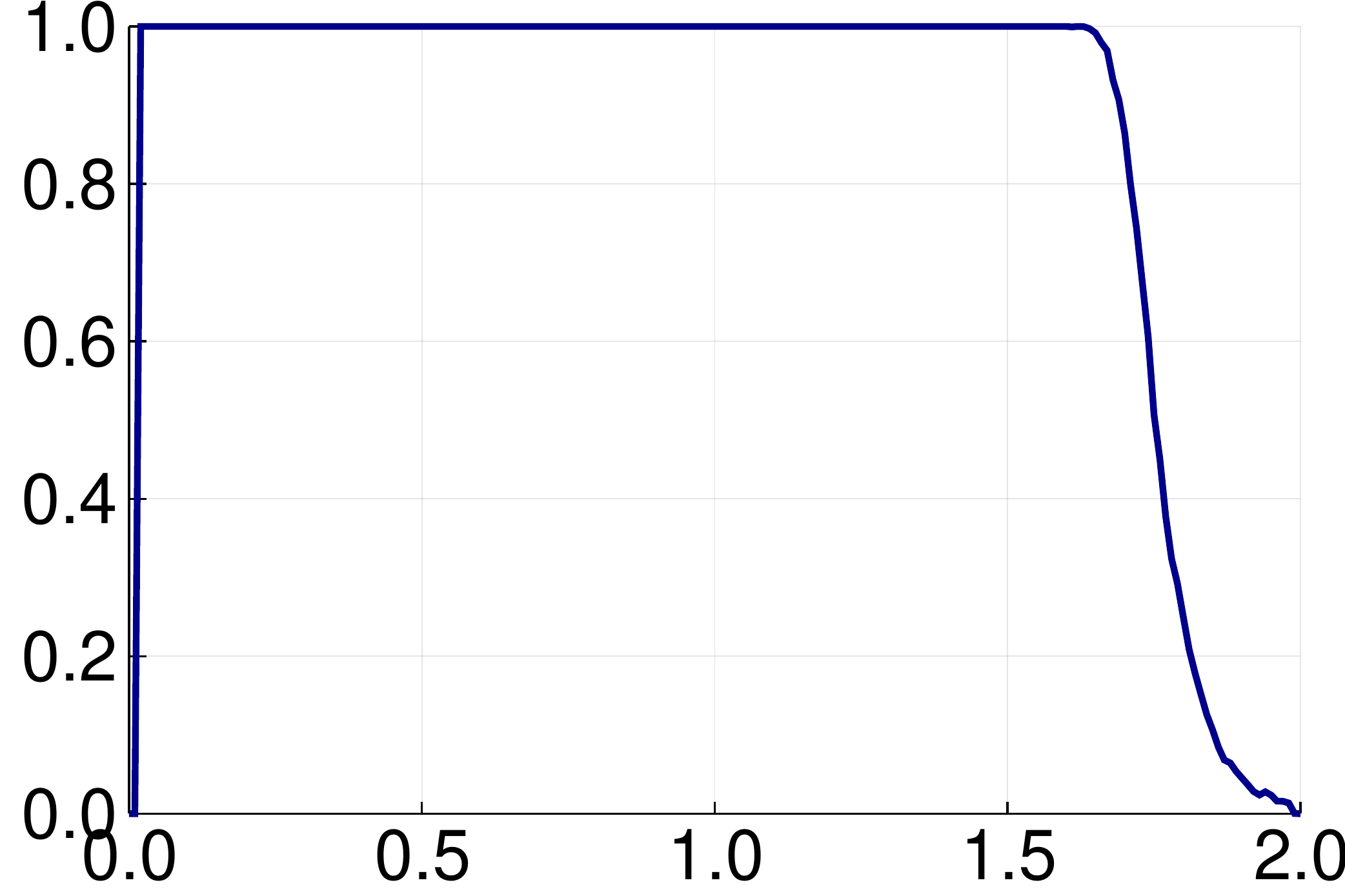}} & 
		\subfloat[CycDR]{\includegraphics[scale=0.13]{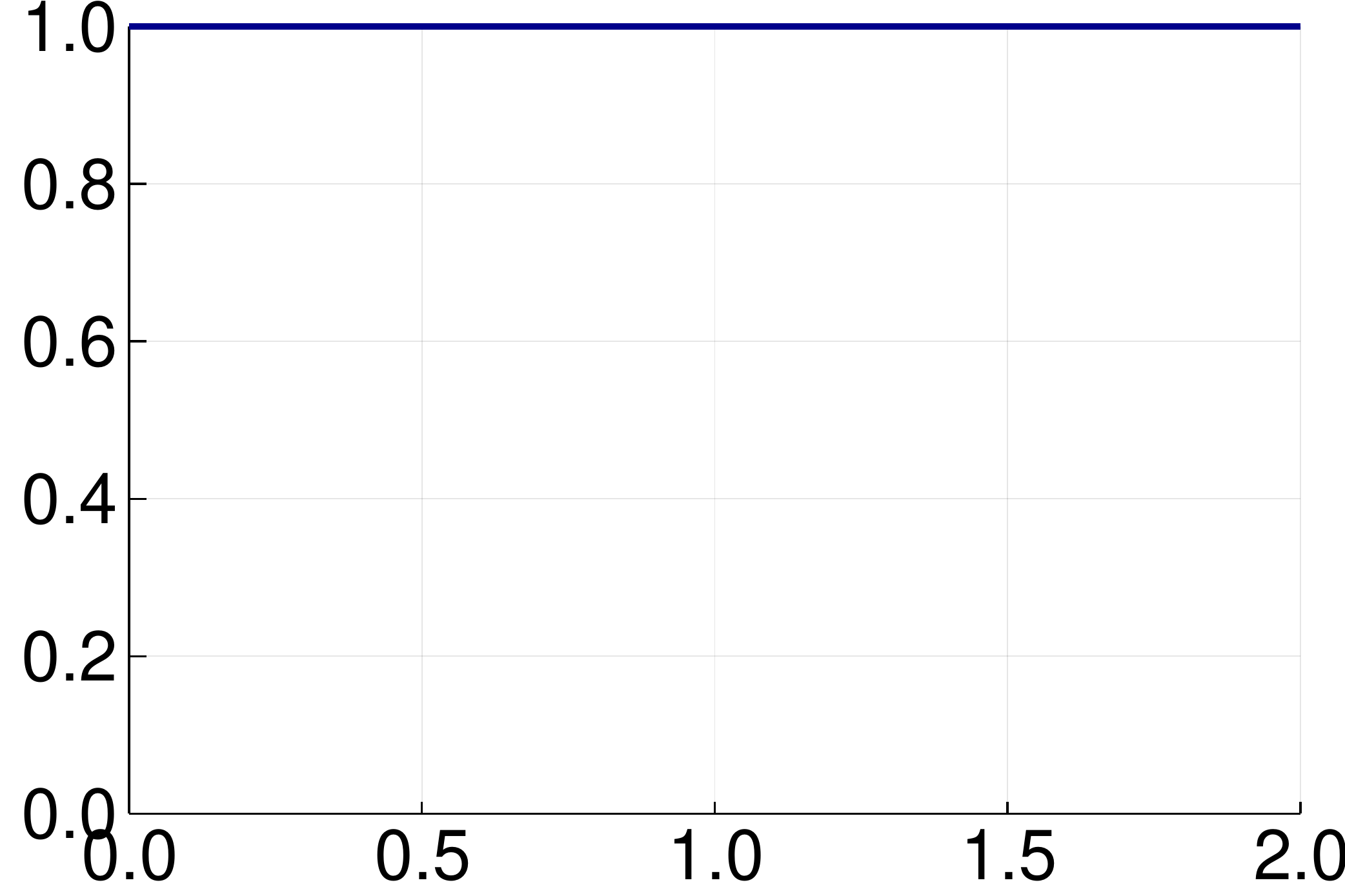}}
	\end{tabular}
	\caption{Success rates in terms of $\lambda$
		for the many sets with few points constellation.}
	\label{fig:BGM_ManySets_FewPts_LambdaCurves}
\end{figure}

\begin{figure}[H]
	\begin{tabular}{cccc}
		\subfloat[CycP]{\includegraphics[scale=0.13]{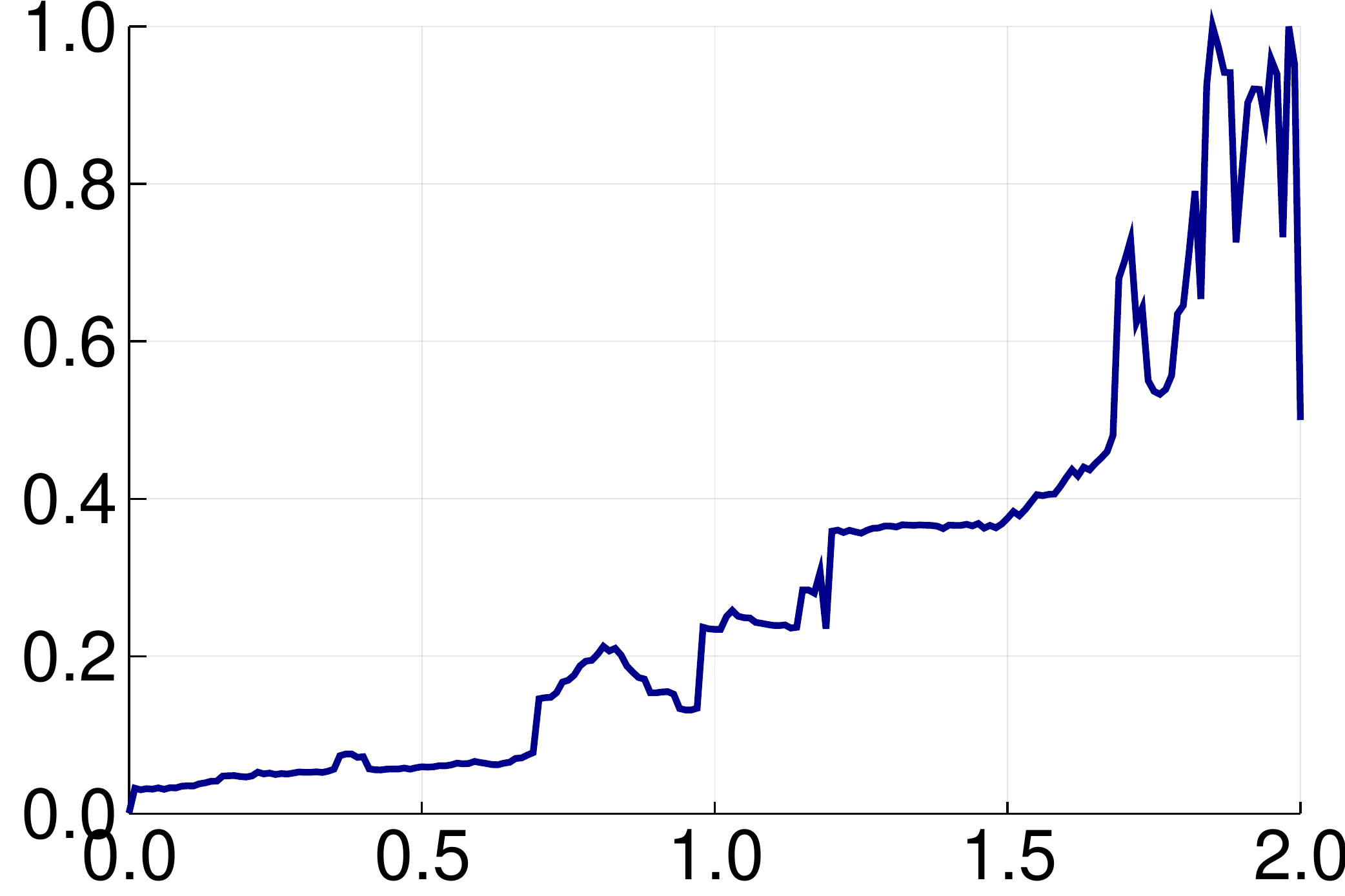}} & 
		\subfloat[ExParP]{\includegraphics[scale=0.13]{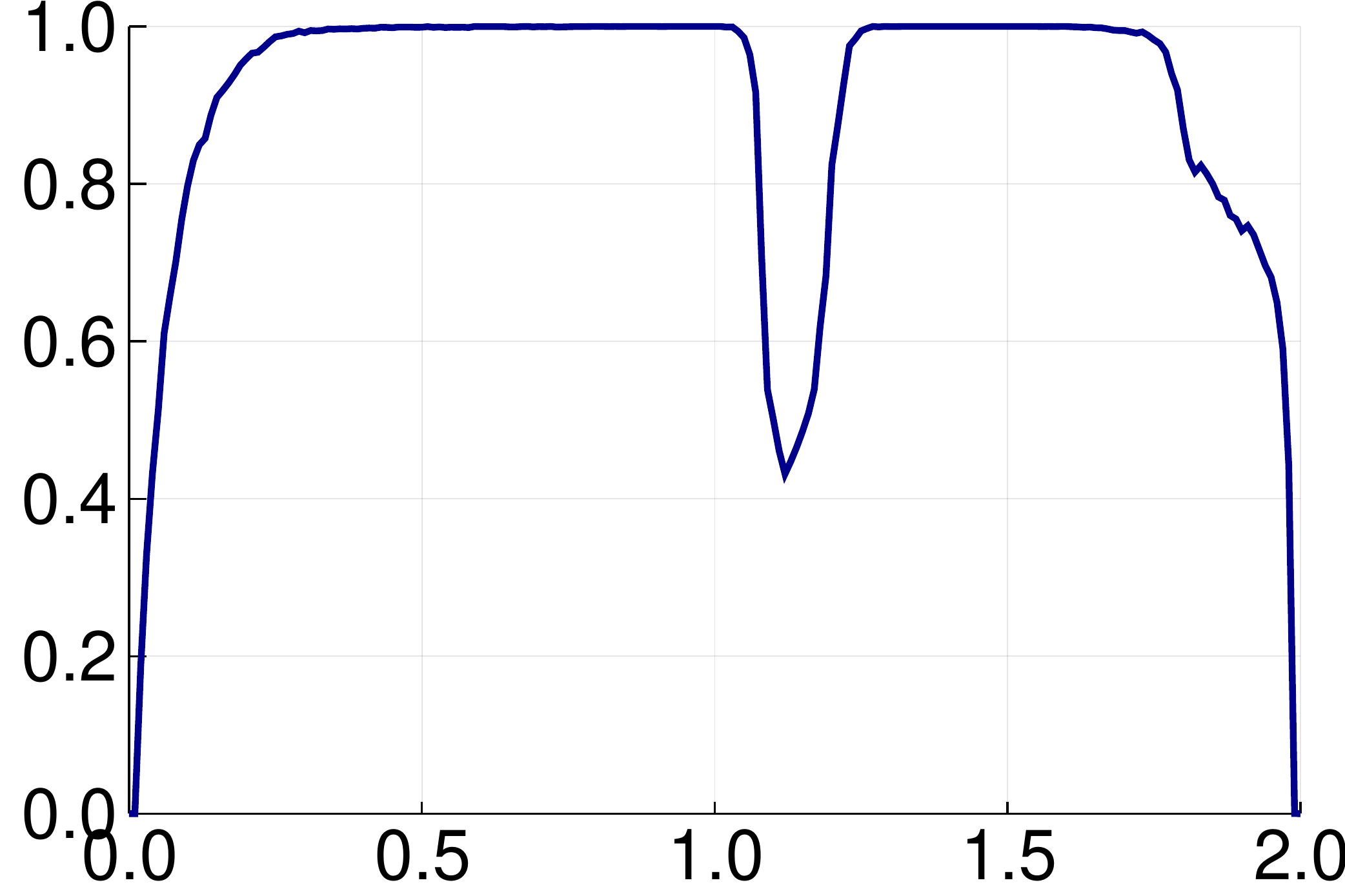}}
		& \subfloat[DR]{\includegraphics[scale=0.13]{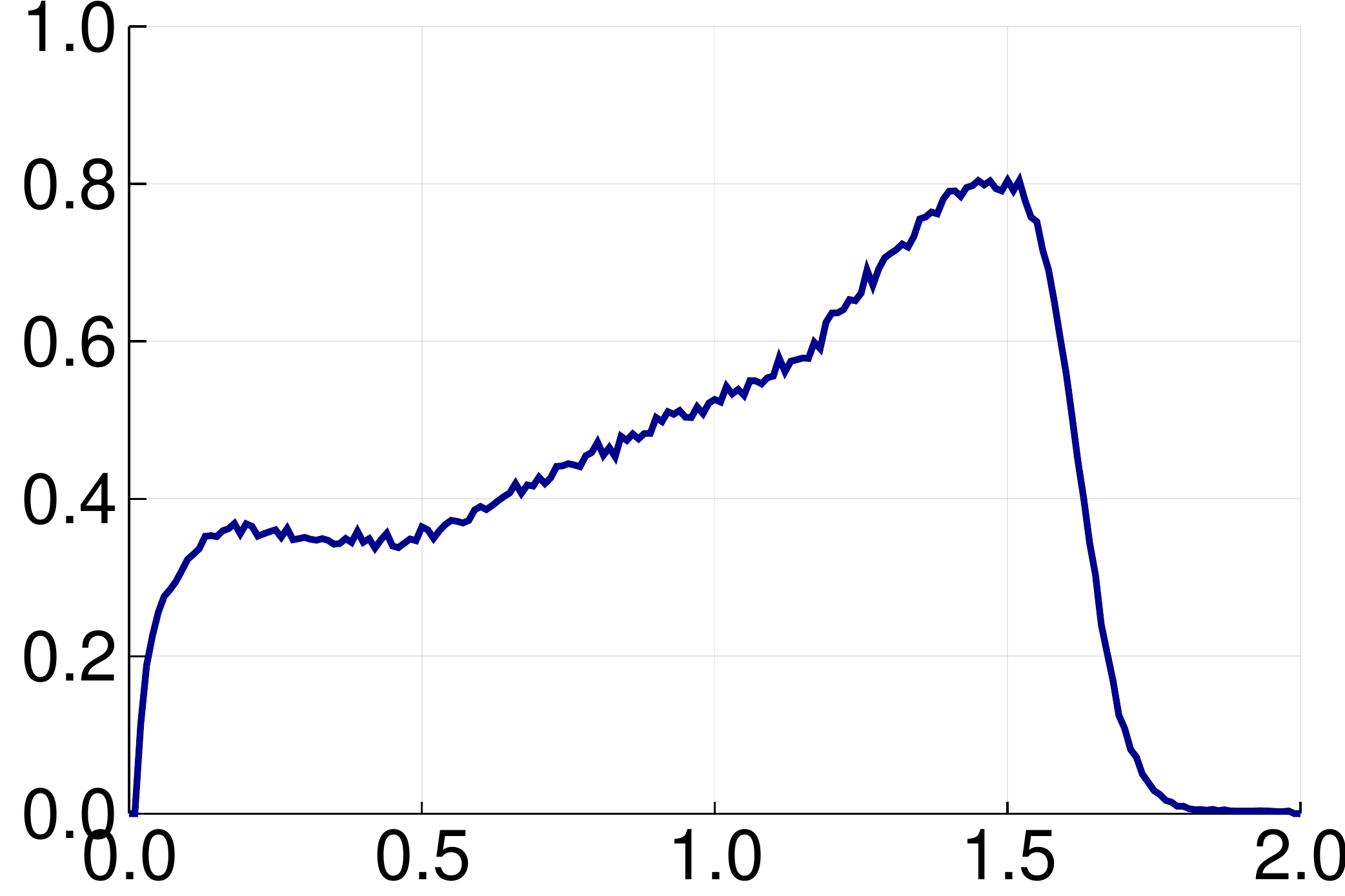}} & 
		\subfloat[CycDR]{\includegraphics[scale=0.13]{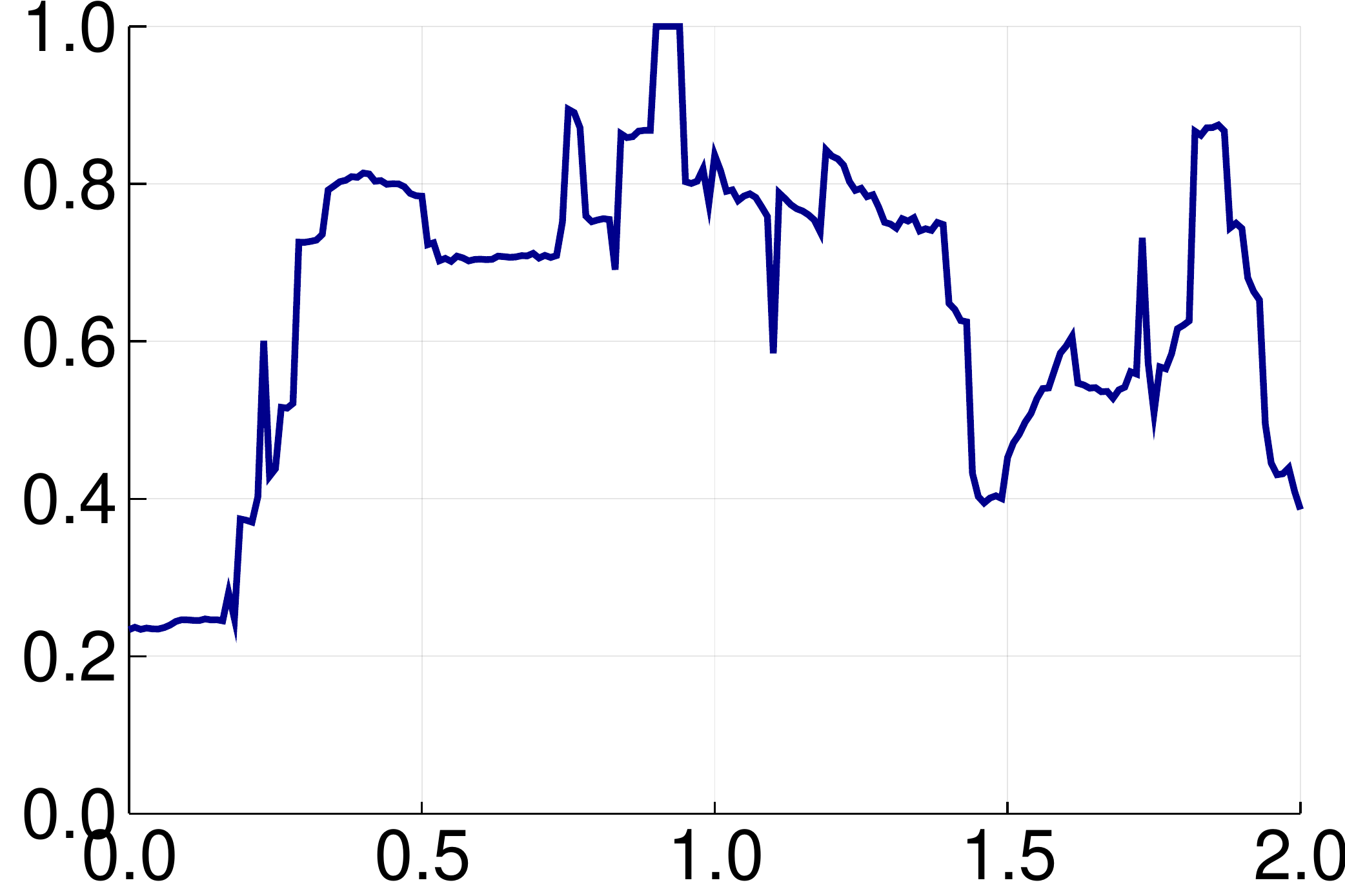}}
	\end{tabular}
	\caption{Success rates in terms of $\lambda$
		for the many sets with many points constellation.}
	\label{fig:BGM_ManySets_ManyPts_LambdaCurves}
\end{figure}

In this section, we consider our four constellations (see Section~\ref{BGMsec:4con})
and run on each of them the four algorithms (see Section~\ref{BGMsec:4alg})
with the parameter $\lambda$ ranging over $\left]0,2\right[$. 
The curves shown in Fig. \ref{fig:BGM_FewSets_FewPts_LambdaCurves}, \ref{fig:BGM_FewSets_ManyPts_LambdaCurves}, \ref{fig:BGM_ManySets_FewPts_LambdaCurves}, and \ref{fig:BGM_ManySets_ManyPts_LambdaCurves} give an estimate of the success rate of each algorithm, evaluated for 200 evenly-spaced values of $\lambda$. For each value of $\lambda$, 5000 starting points are drawn from $[-10,10]\times[-10,10]$ using Quasi-Monte Carlo sampling, and the success rate is estimated by dividing the number of times the algorithm is successful by this number of starting points. Thus,  a ``best'' parameter $\BGMlbest$ is 
determined. It is this parameter that we will use to compare with
the default parameter $\BGMldefault$, which is $1$ in all cases. 
 
\vspace{-0.4cm}
\vspace{-1cm}
\vspace{-1cm}
\vspace{-1cm}

\runinhead{Discussion}

For each of the four constellations considered above,
we visually inspected the $\lambda$-curves indicating
success rates. We then picked for each algorithm a parameter called
$\BGMlbest$ which improved performance over the default parameter
$\BGMldefault=1$. 
The results are recorded in the following table.

	\begin{figure}[H]
	\begin{svgraybox}
		\begin{center}
\begin{tabular}{@{}r r r r r@{}}
\toprule
Algorithm & ~~~~~~~~~CycP &~~ExParP &~~~~~~~~DR &~~CycDR\\ \midrule
$\BGMldefault$ & 1.0 & 1.0 & 1.0 & 1.0\\
$\BGMlbest$ & 1.5 &0.8 &1.6 &1.2\\
\bottomrule
\end{tabular}
\end{center}
\end{svgraybox}
\caption{Best parameters $\BGMlbest$ chosen by inspecting the success rates curves}
\label{BGMfig:lambdatable}
\end{figure}

We will use the parameters $\BGMlbest$ for the experiments in subsequent sections.

\section{Tracking orbits}
\label{BGMsec:track}

In this section, we consider our four given constellations
(see Section~\ref{BGMsec:4con}).
For each constellation, which is organized in a separate subsection, 
the same starting point is used.
We then consider each of our four fixed algorithms 
(see Section~\ref{BGMsec:4alg})
and show orbits for $\BGMldefault=1$ and for $\BGMlbest$ 
(see Section~\ref{BGMsec:det}),
and the corresponding feasibility measure $d$ (see Section \ref{BGMsec:setup}).

\subsection{Few sets with few points}

\BGMvspace{0.5}

\begin{figure}[H]
\begin{tabular}{cccc}
\large \textbf{CycP} & & & \\
\subfloat[$\BGMldefault$ orbit]{\includegraphics[scale=0.15]{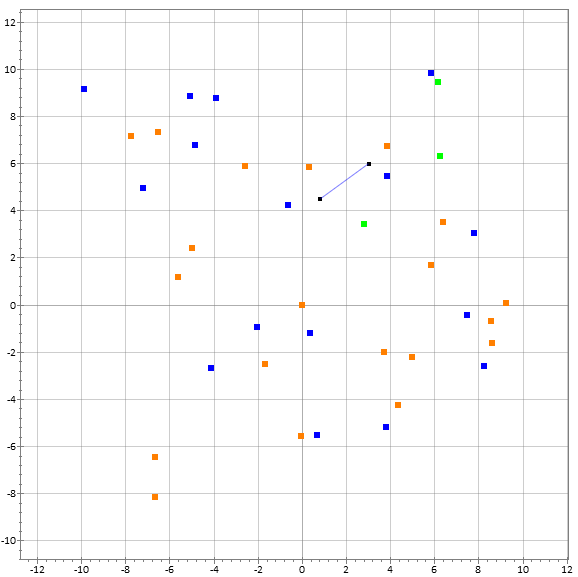}} & 
\subfloat[$\BGMldefault$ error]{\includegraphics[scale=0.15]{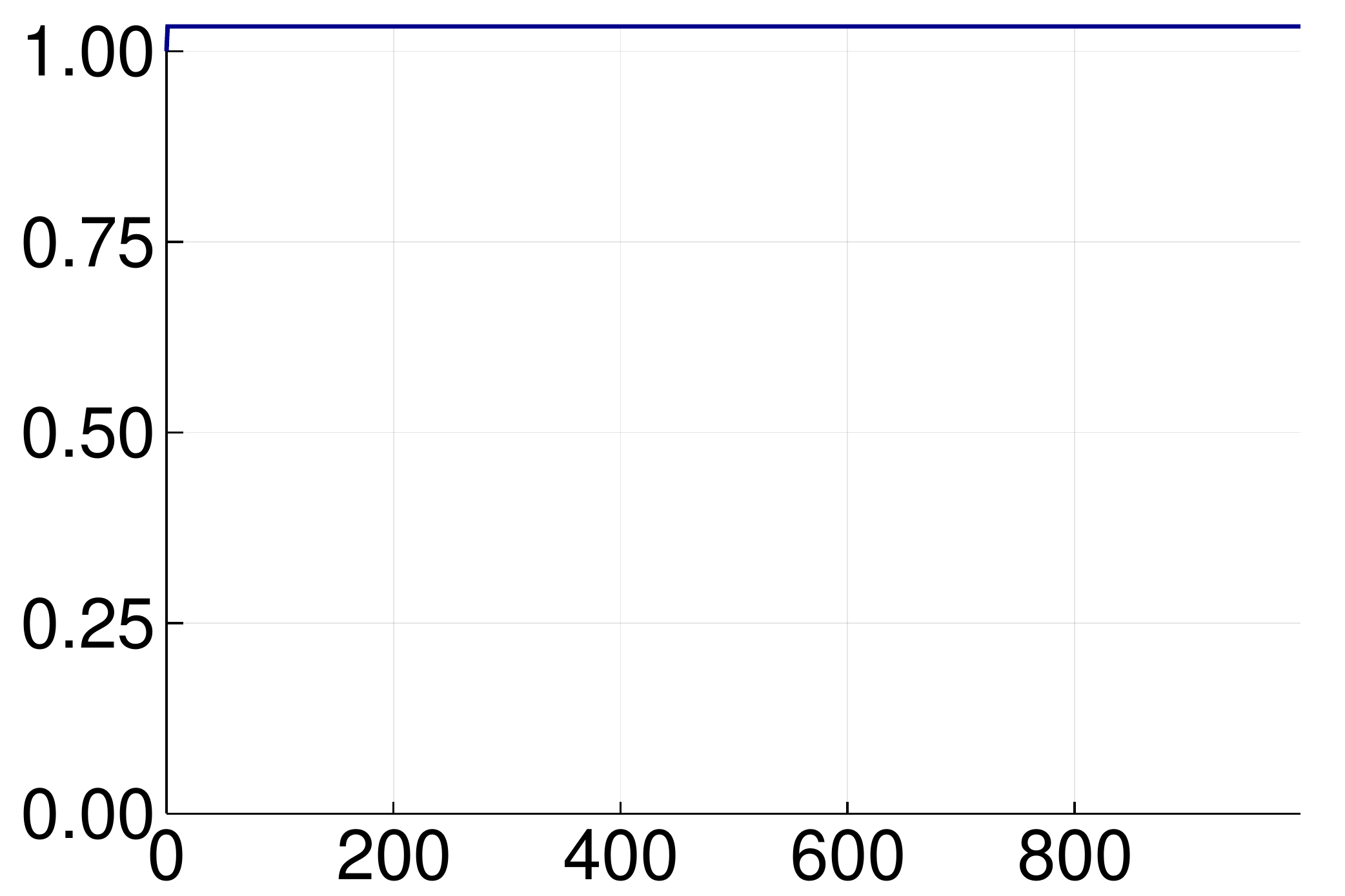}} & 
\subfloat[$\BGMlbest$ orbit]{\includegraphics[scale=0.15]{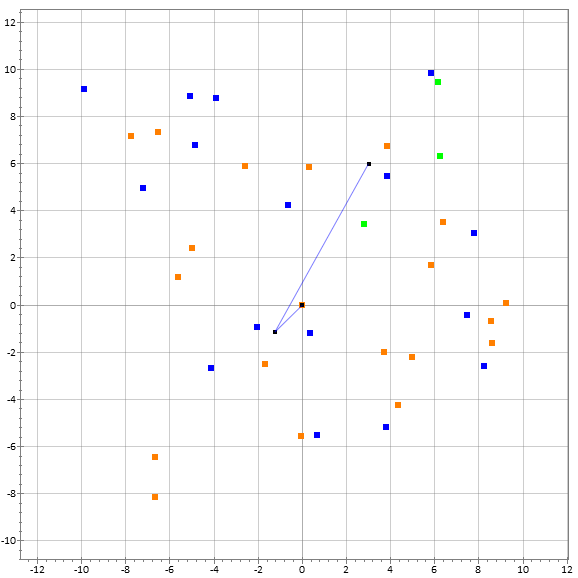}} &
\subfloat[$\BGMlbest$ error]{\includegraphics[scale=0.15]{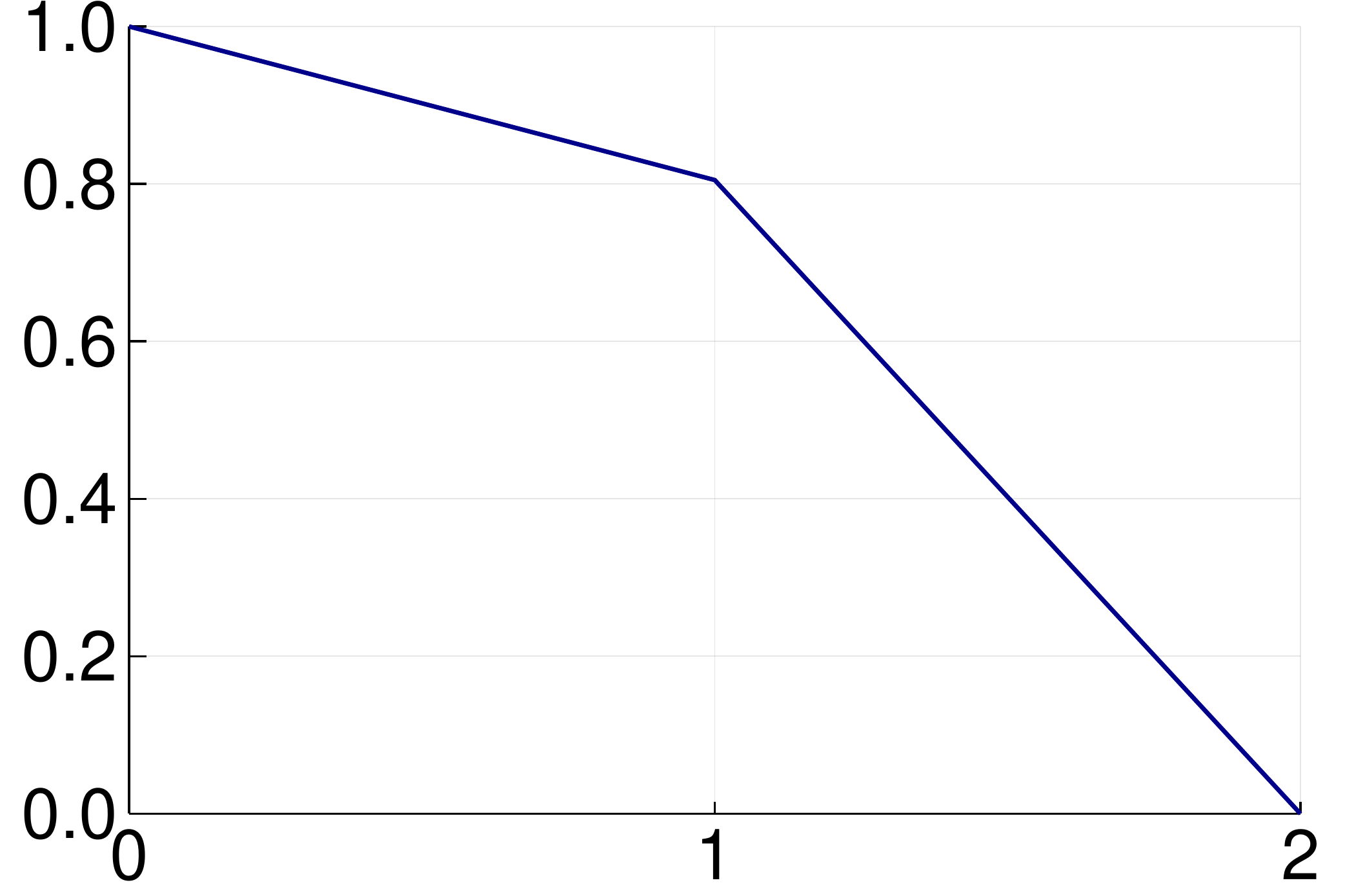}}
\\
 & & & \\ 
\large \textbf{ExParP} & & & \\
\subfloat[$\BGMldefault$ orbit]{\includegraphics[scale=0.15]{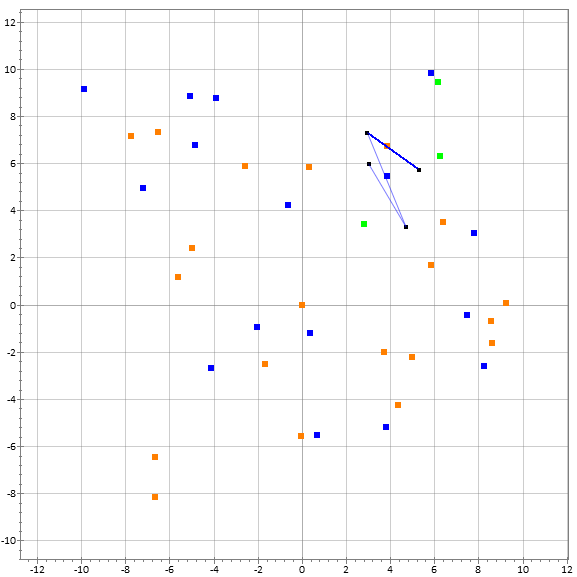}} & 
\subfloat[$\BGMldefault$ error]{\includegraphics[scale=0.15]{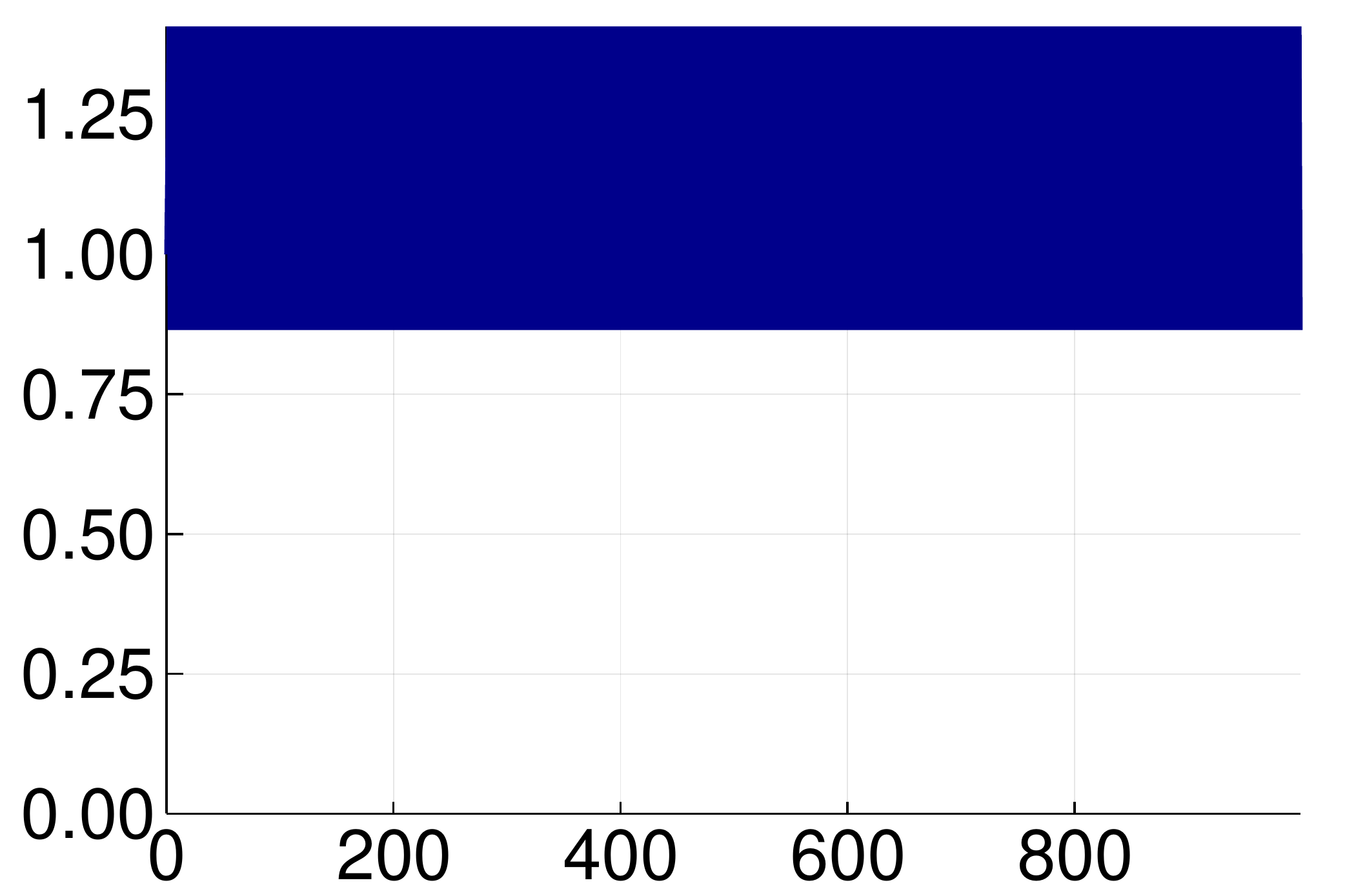}} & 
\subfloat[$\BGMlbest$ orbit]{\includegraphics[scale=0.15]{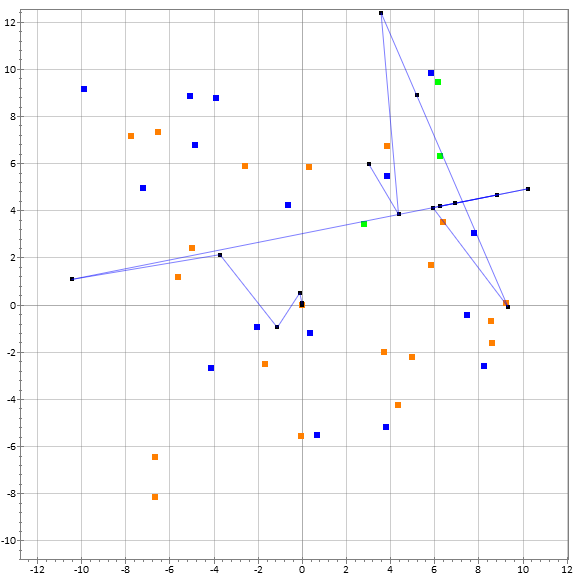}} &
\subfloat[$\BGMlbest$ error]{\includegraphics[scale=0.15]{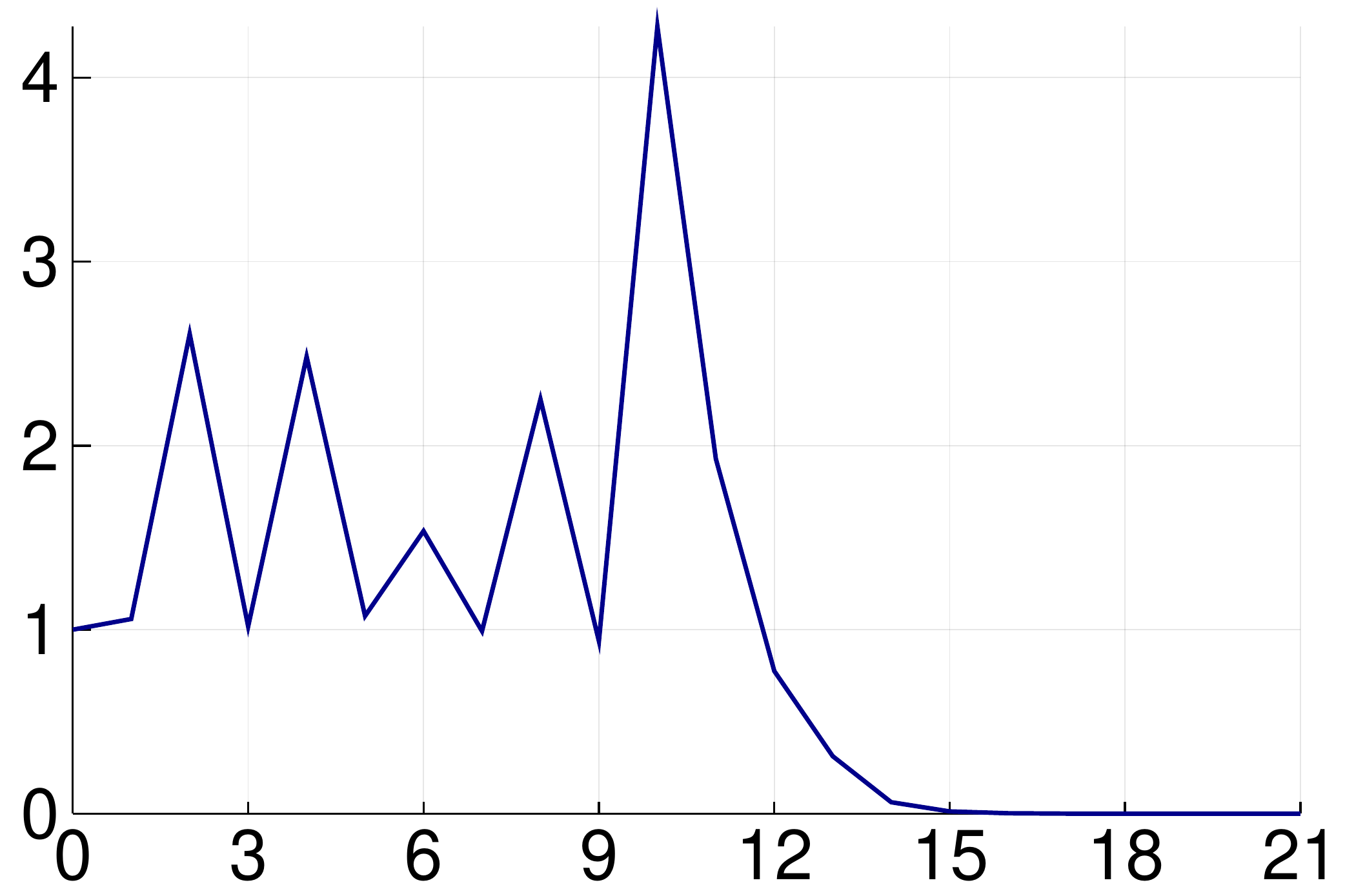}}
\\
& & & \\ 
\large \textbf{DR} & & & \\
\subfloat[$\BGMldefault$ orbit]{\includegraphics[scale=0.15]{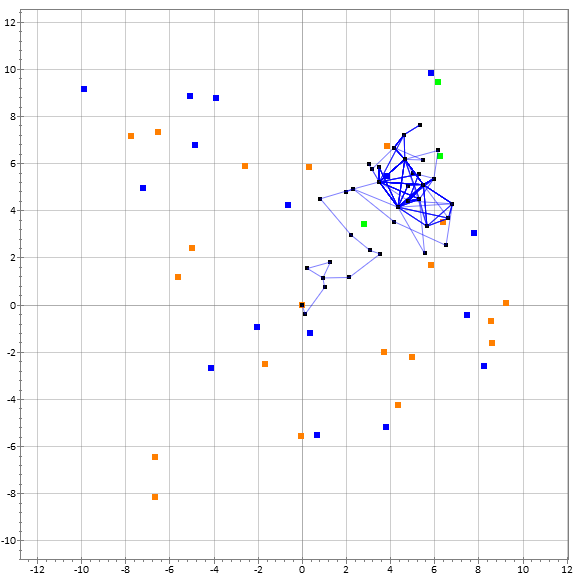}} & 
\subfloat[$\BGMldefault$ error]{\includegraphics[scale=0.15]{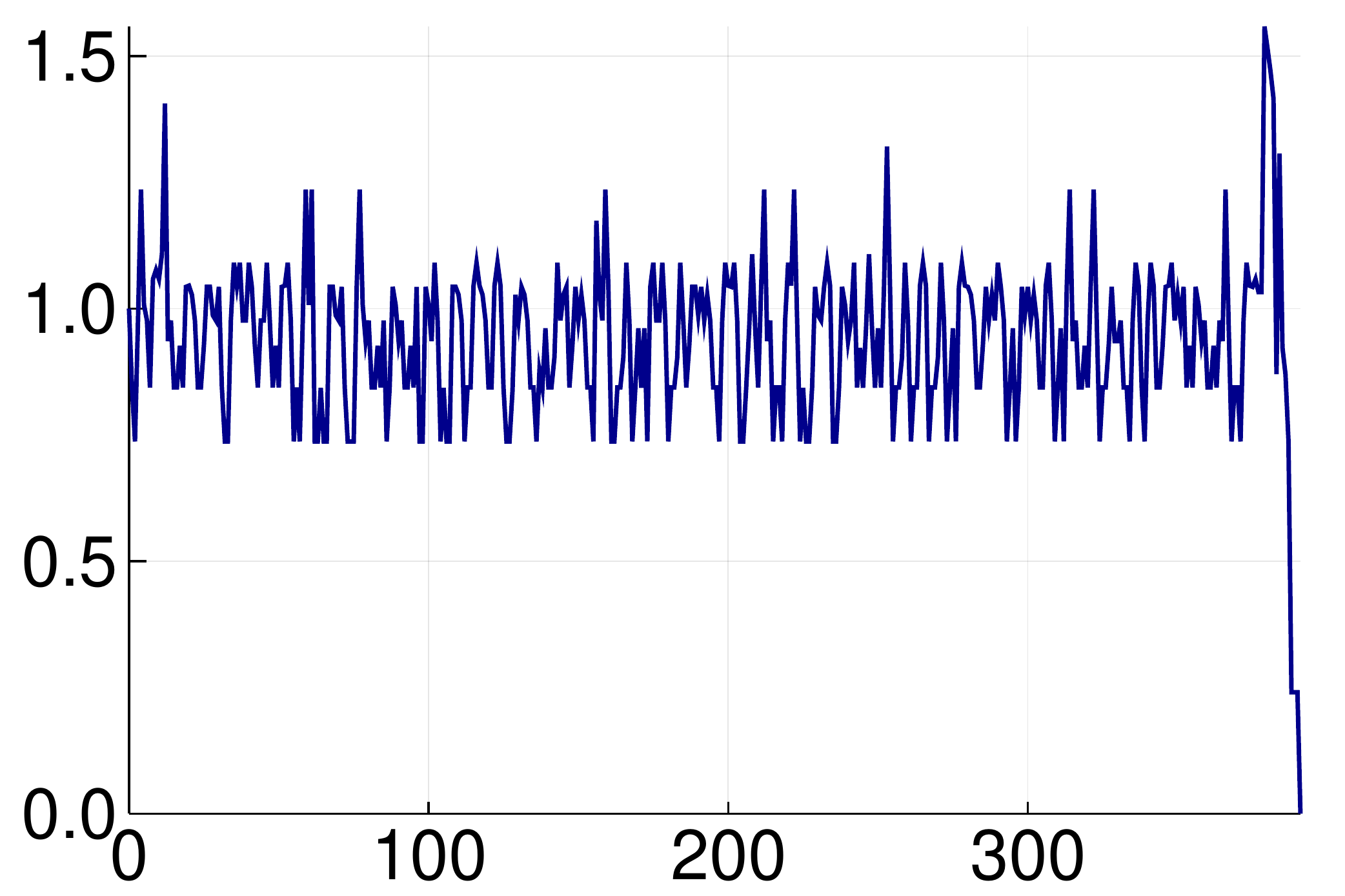}} & 
\subfloat[$\BGMlbest$ orbit]{\includegraphics[scale=0.15]{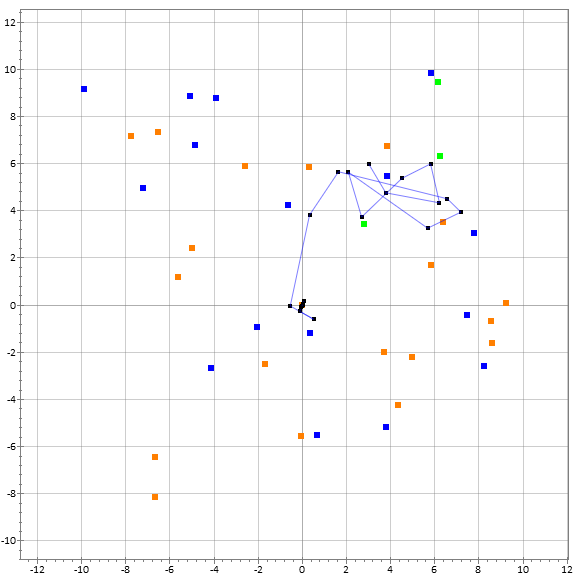}} &
\subfloat[$\BGMlbest$ error]{\includegraphics[scale=0.15]{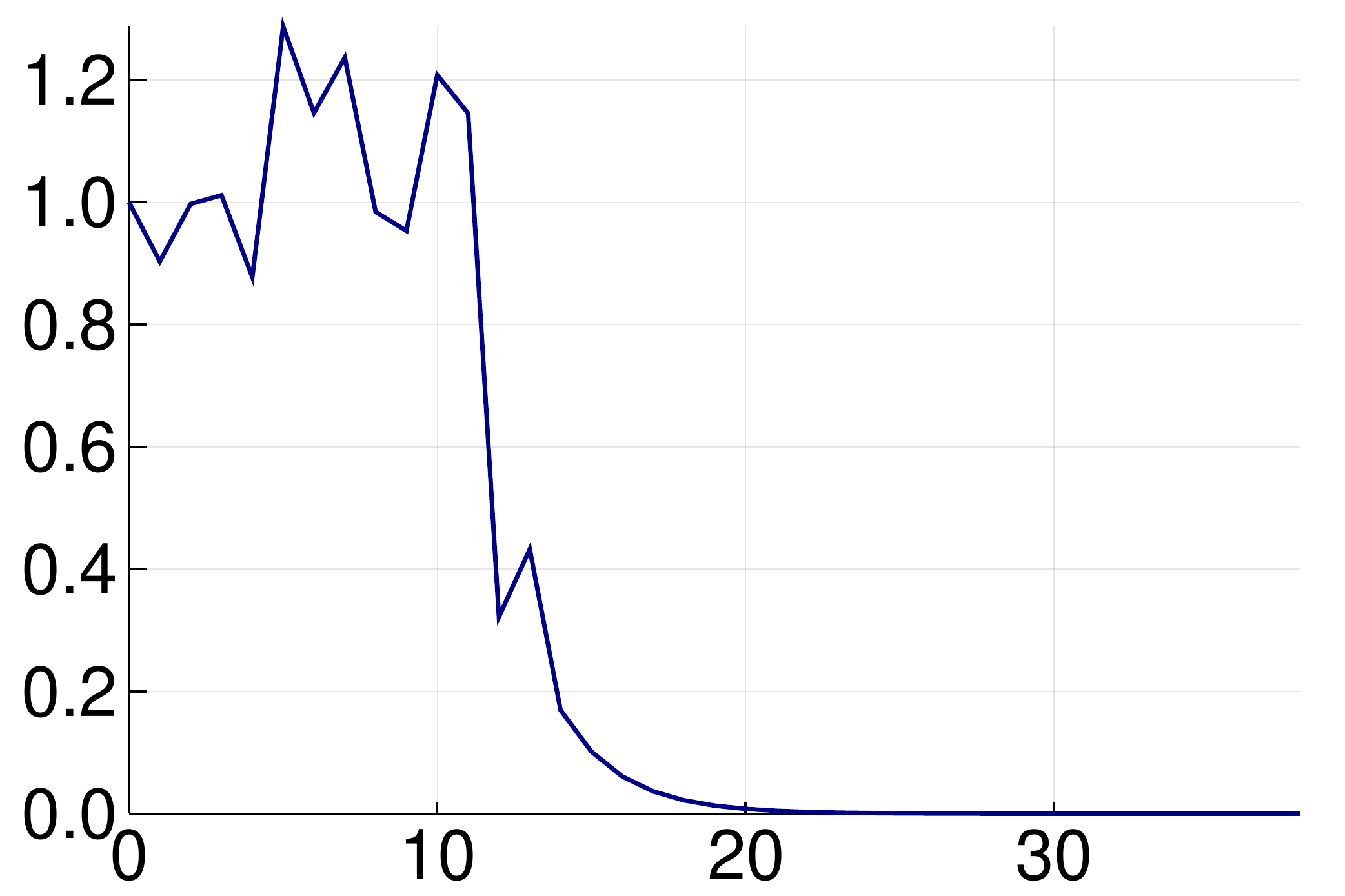}}
\\
& & & \\ 
\large \textbf{CycDR} & & & \\
\subfloat[$\BGMldefault$ orbit]{\includegraphics[scale=0.15]{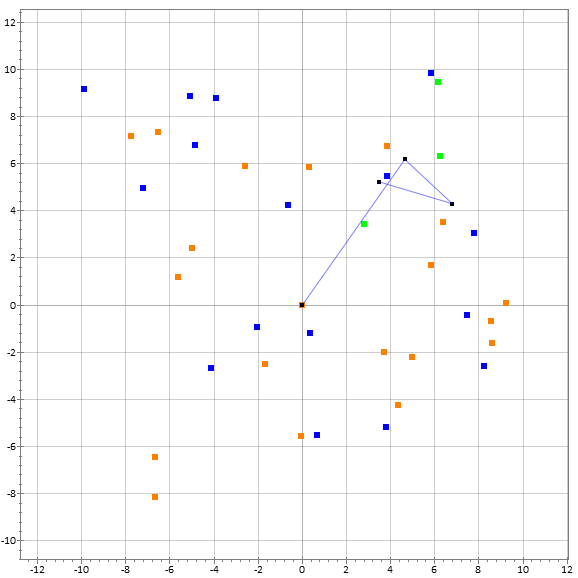}} & 
\subfloat[$\BGMldefault$ error]{\includegraphics[scale=0.15]{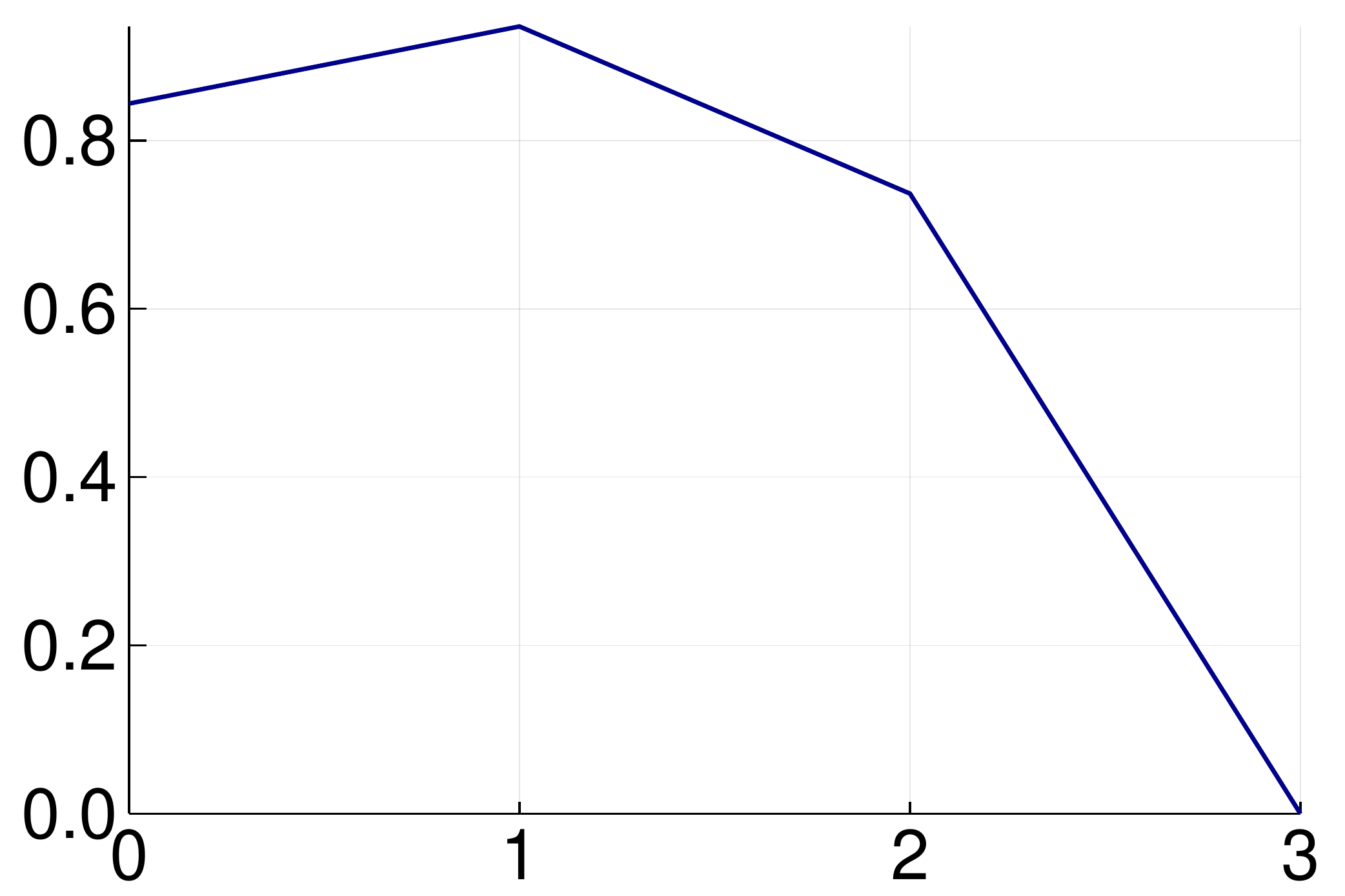}} & 
\subfloat[$\BGMlbest$ orbit]{\includegraphics[scale=0.15]{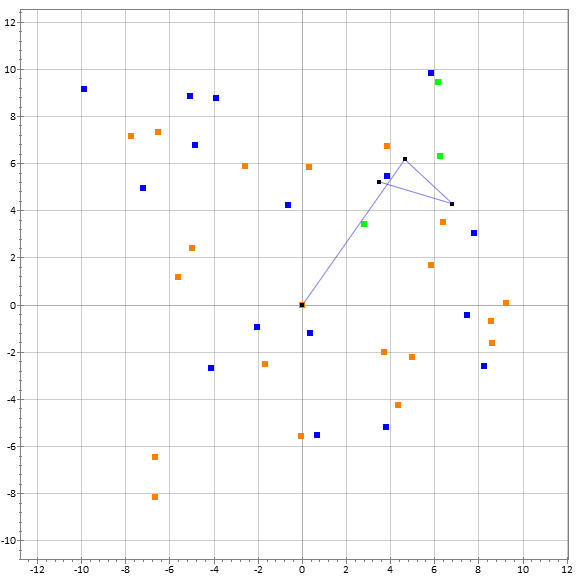}} &
\subfloat[$\BGMlbest$ error]{\includegraphics[scale=0.15]{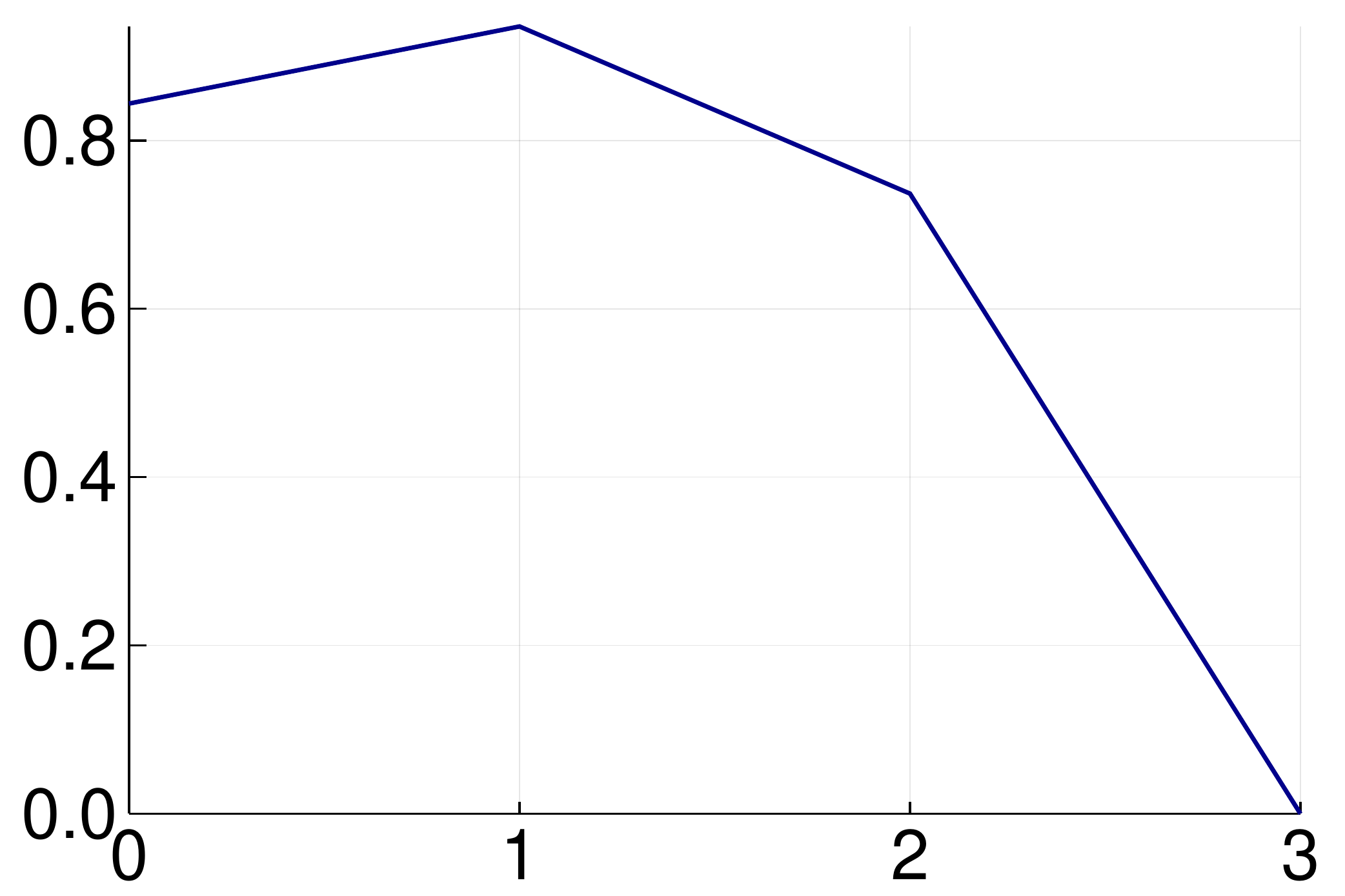}}
\end{tabular}
\caption{Orbits and errors for CycP, ExParP, DR, and CycDR in the few sets with few points constellation}
\end{figure}

\subsection{Few sets with many points}

\BGMvspace{0.5}

\begin{figure}[H]
\begin{tabular}{cccc}
\large \textbf{CycP} & & & \\
\subfloat[$\BGMldefault$ orbit]{\includegraphics[scale=0.15]{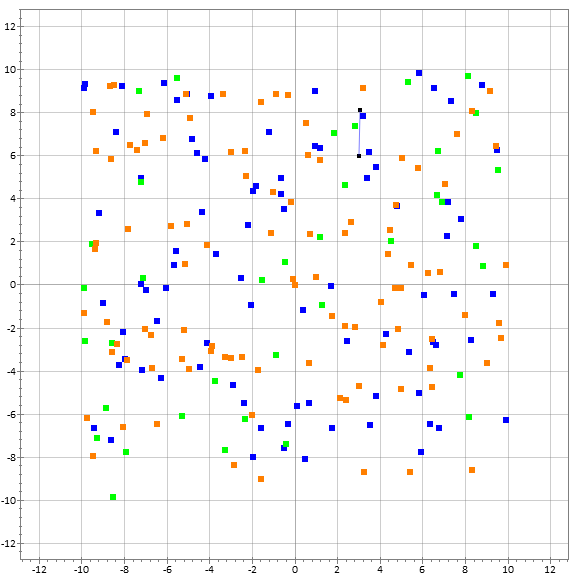}} & 
\subfloat[$\BGMldefault$ error]{\includegraphics[scale=0.15]{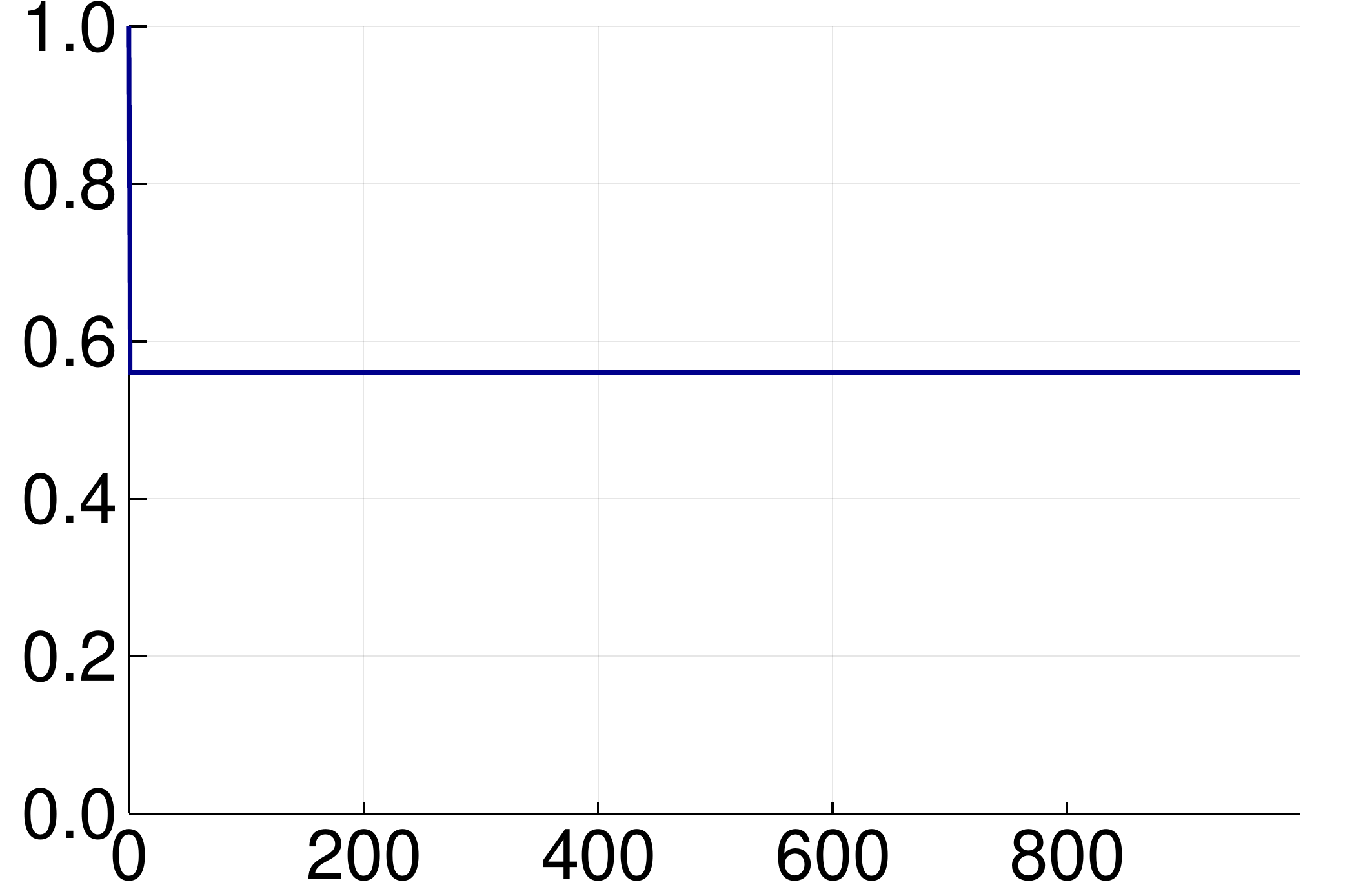}} & 
\subfloat[$\BGMlbest$ orbit]{\includegraphics[scale=0.15]{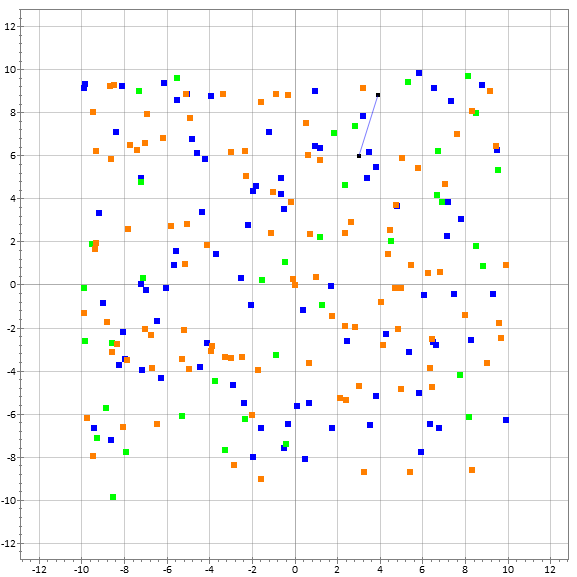}} &
\subfloat[$\BGMlbest$ error]{\includegraphics[scale=0.15]{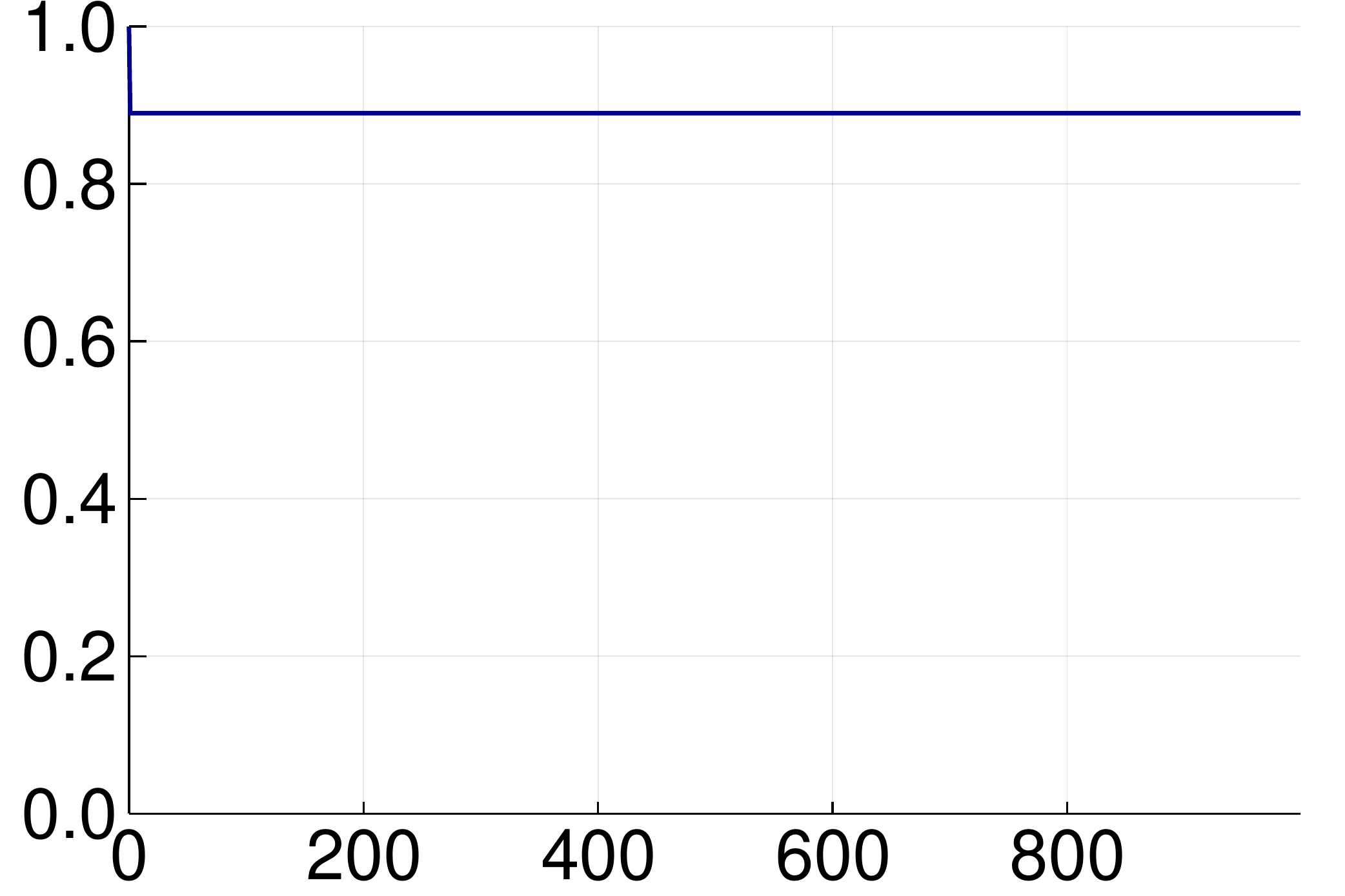}}
\\
& & & \\ 
\large \textbf{ExParP} & & & \\
\subfloat[$\BGMldefault$ orbit]{\includegraphics[scale=0.15]{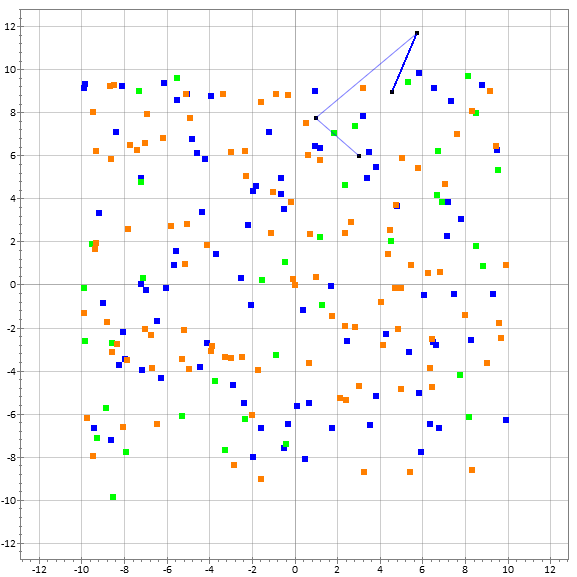}} & 
\subfloat[$\BGMldefault$ error]{\includegraphics[scale=0.15]{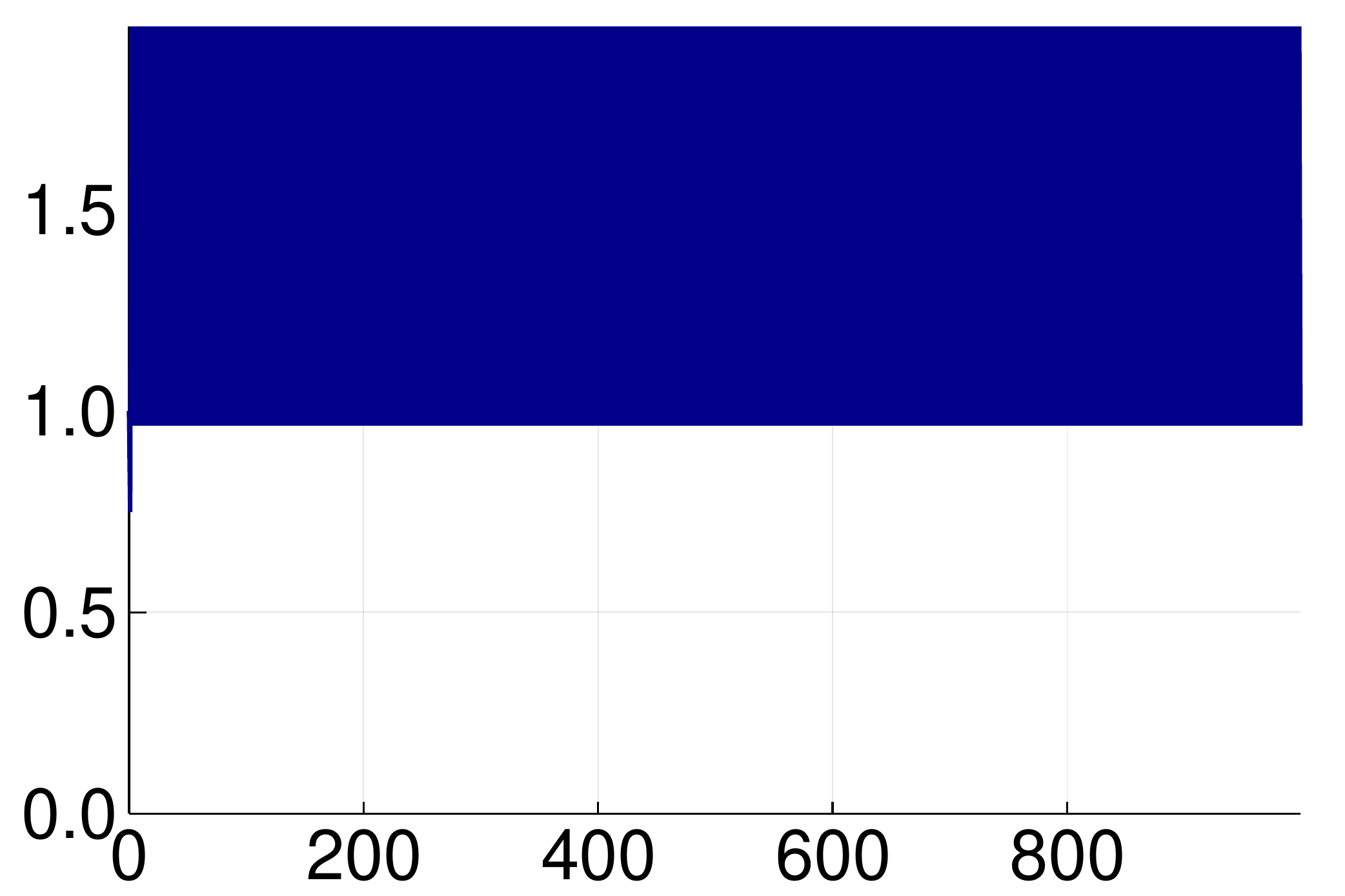}} & 
\subfloat[$\BGMlbest$ orbit]{\includegraphics[scale=0.15]{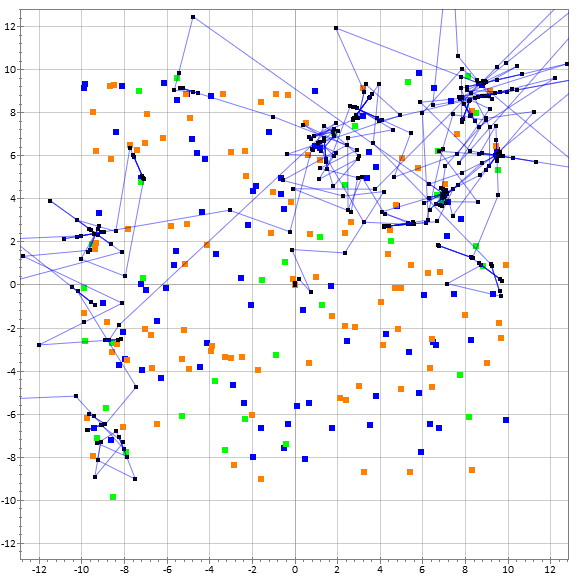}} &
\subfloat[$\BGMlbest$ error]{\includegraphics[scale=0.15]{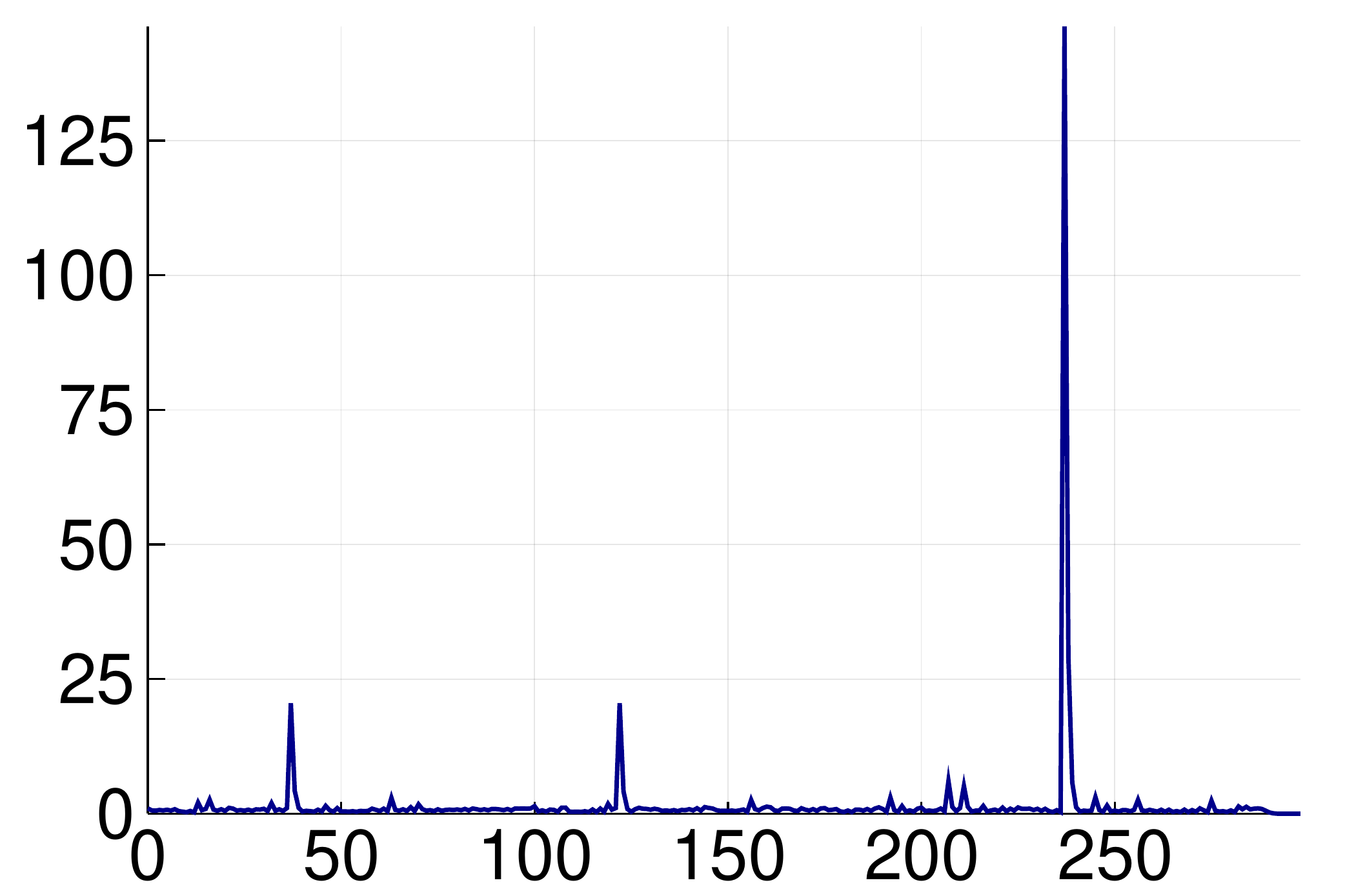}}
\\
& & & \\ 
\large \textbf{DR} & & & \\
\subfloat[$\BGMldefault$ orbit]{\includegraphics[scale=0.15]{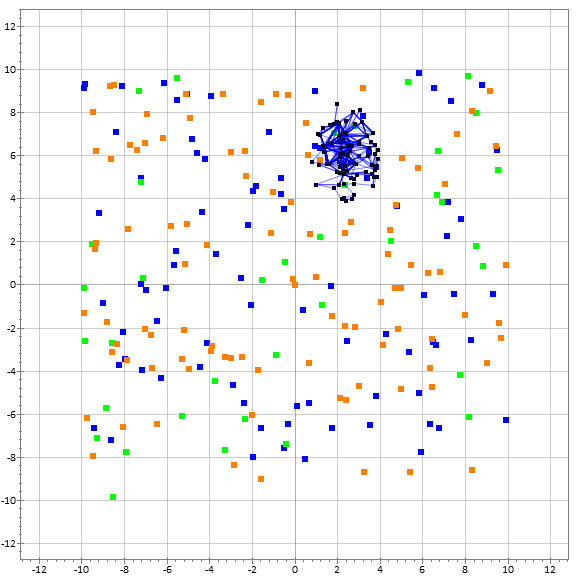}} & 
\subfloat[$\BGMldefault$ error]{\includegraphics[scale=0.15]{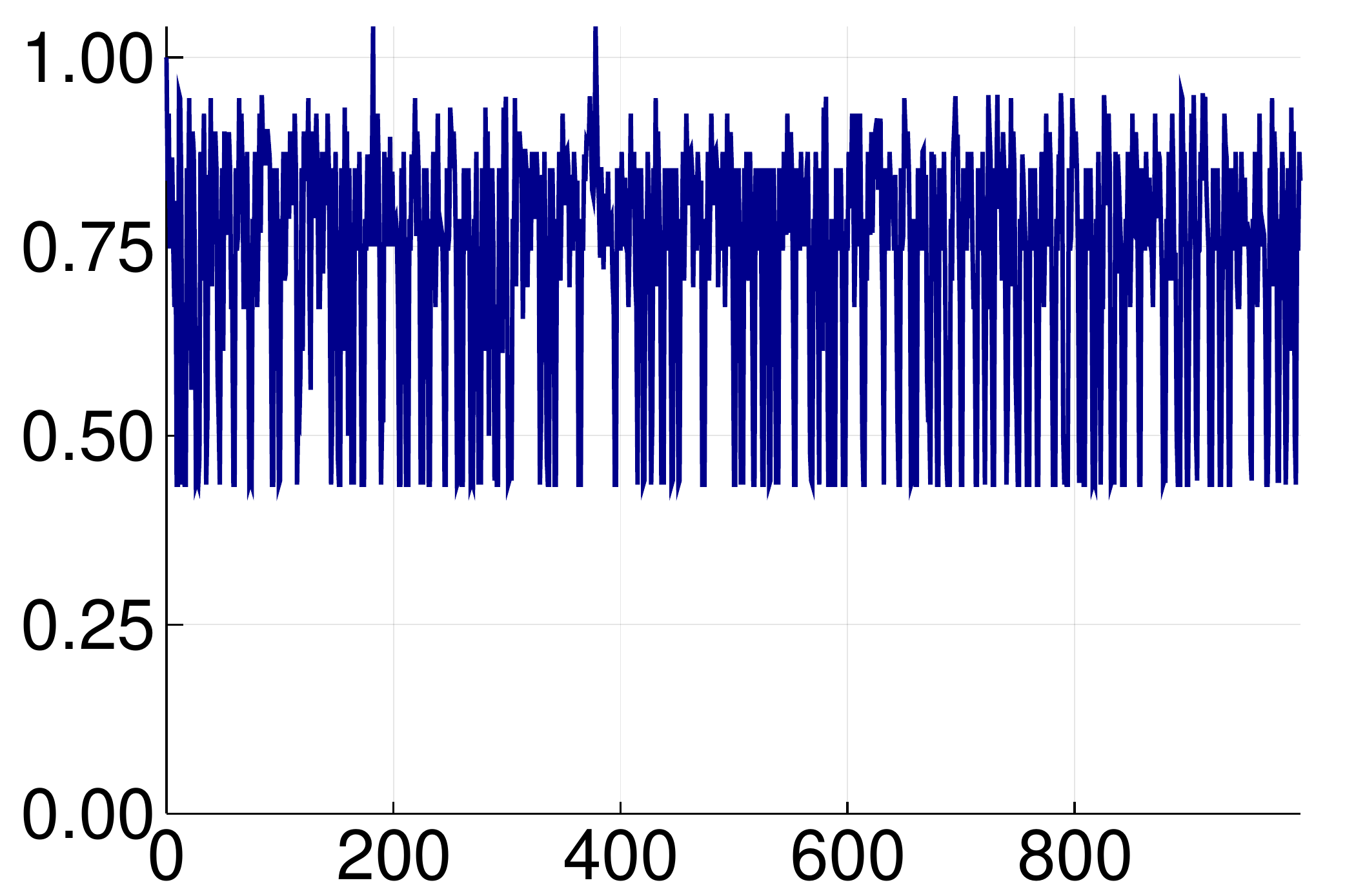}} & 
\subfloat[$\BGMlbest$ orbit]{\includegraphics[scale=0.15]{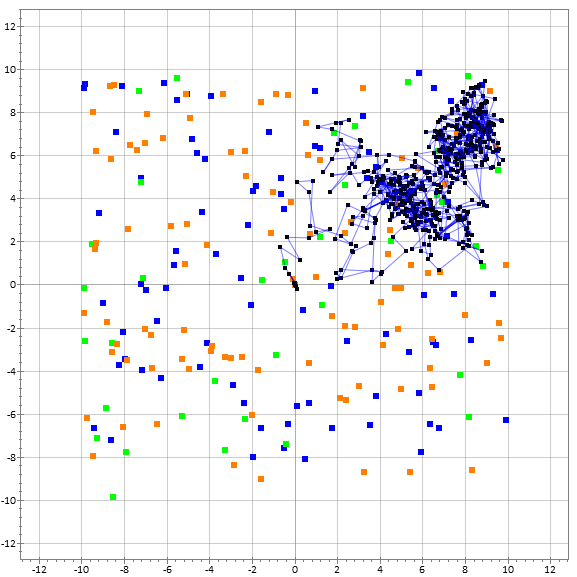}} &
\subfloat[$\BGMlbest$ error]{\includegraphics[scale=0.15]{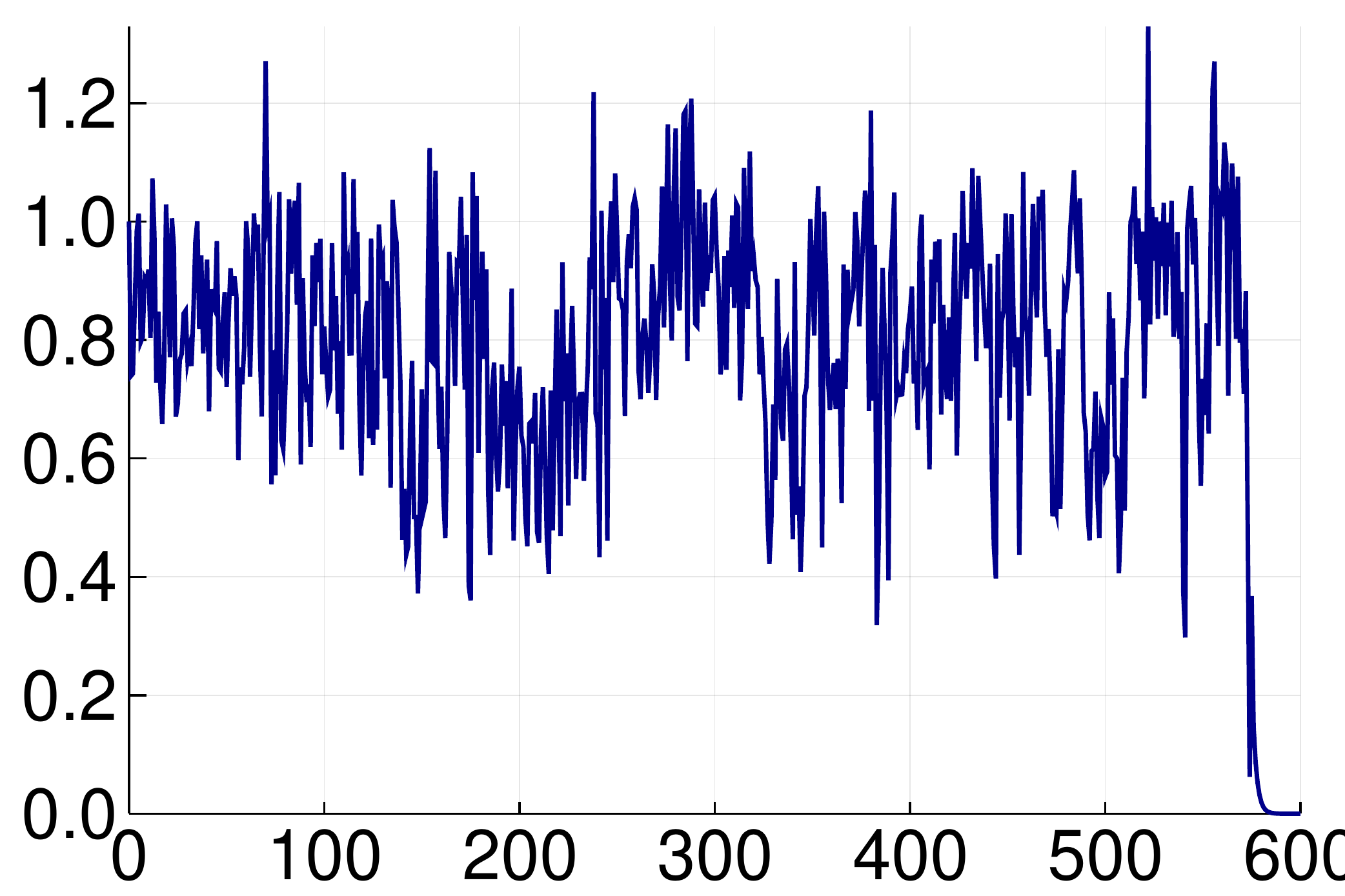}}
\\
& & & \\ 
\large \textbf{CycDR} & & & \\
\subfloat[$\BGMldefault$ orbit]{\includegraphics[scale=0.15]{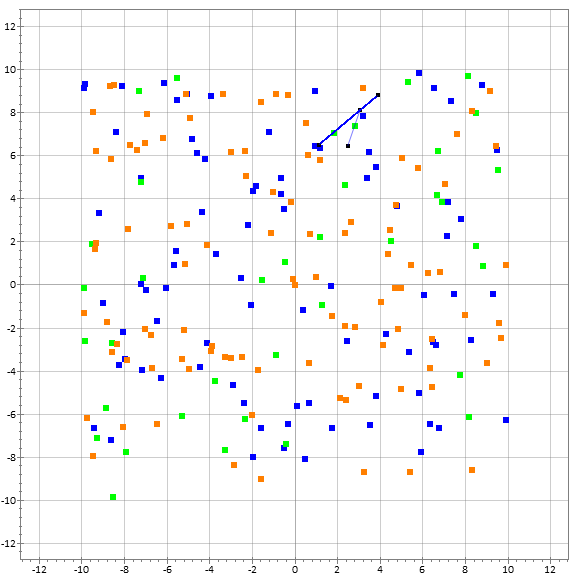}} & 
\subfloat[$\BGMldefault$ error]{\includegraphics[scale=0.15]{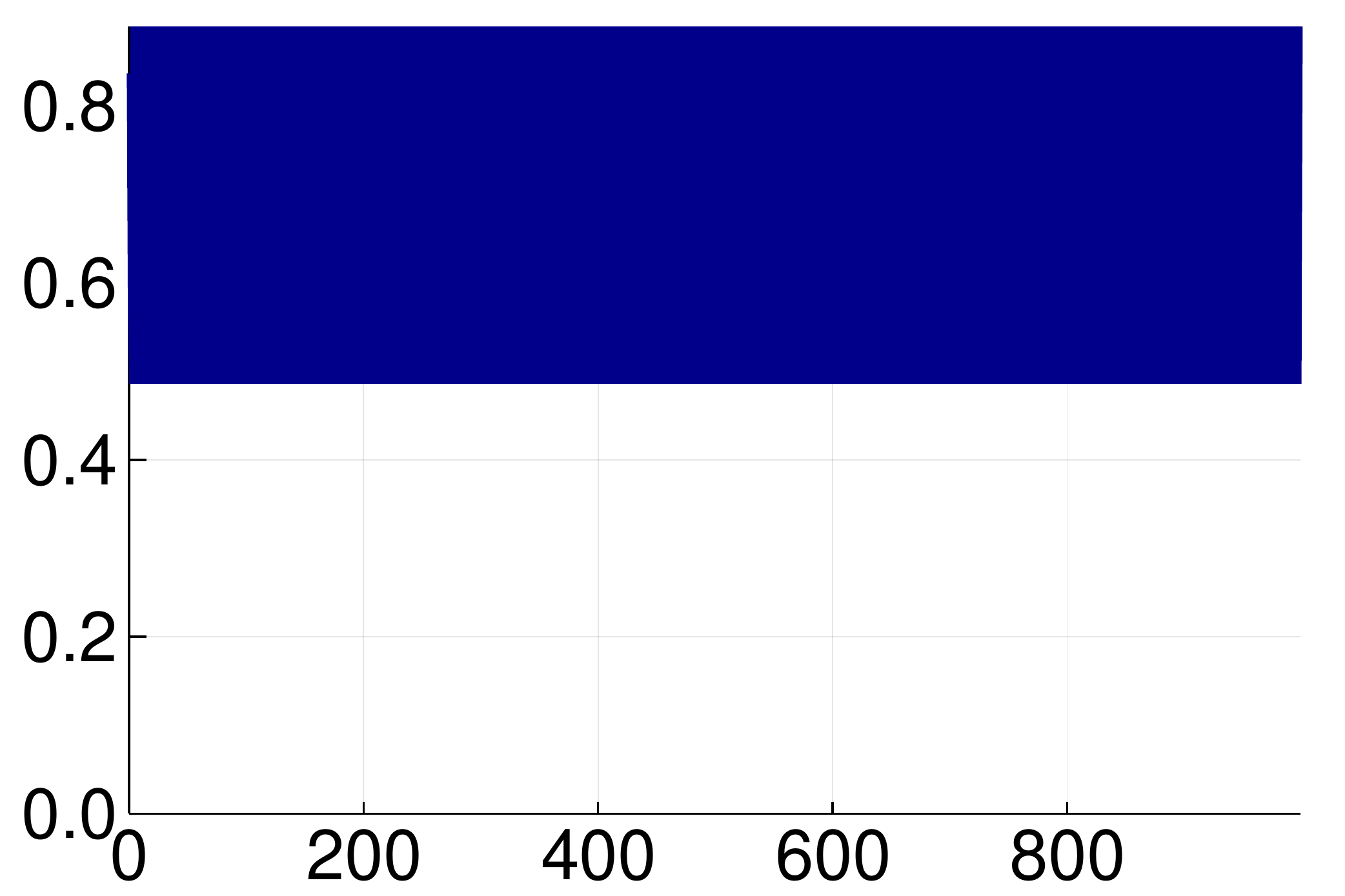}} & 
\subfloat[$\BGMlbest$ orbit]{\includegraphics[scale=0.15]{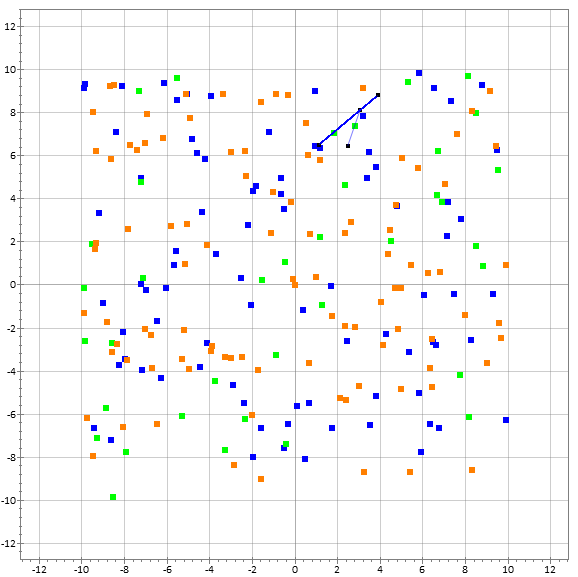}} &
\subfloat[$\BGMlbest$ error]{\includegraphics[scale=0.15]{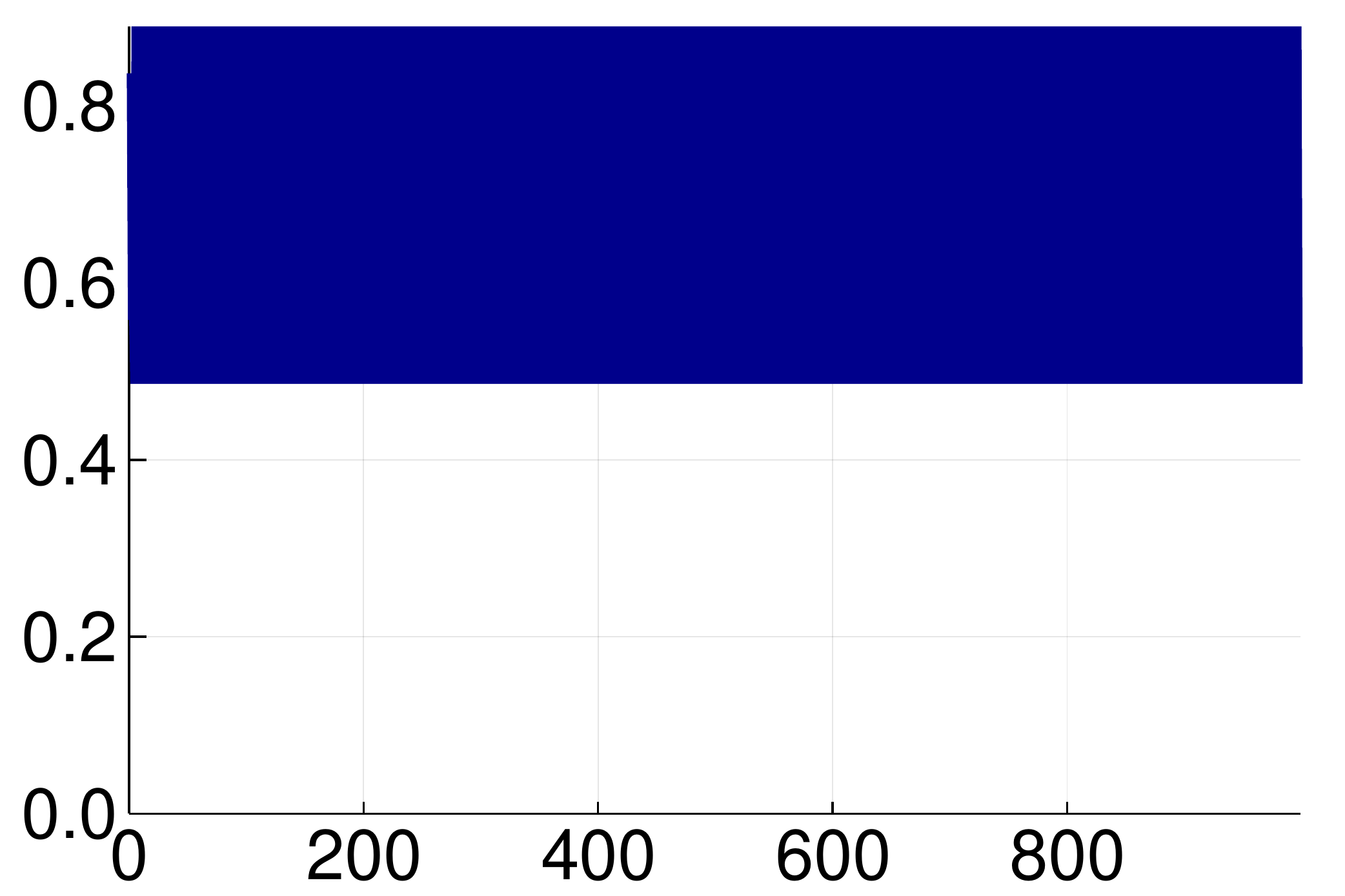}}
\end{tabular}
\caption{Orbits and errors for CycP, ExParP, DR, and CycDR in the few sets with many points constellation}
\end{figure}

\subsection{Many sets with few points}

\BGMvspace{0.5}

\begin{figure}[H]
\begin{tabular}{cccc}
\large \textbf{CycP} & & & \\
\subfloat[$\BGMldefault$ orbit]{\includegraphics[scale=0.15]{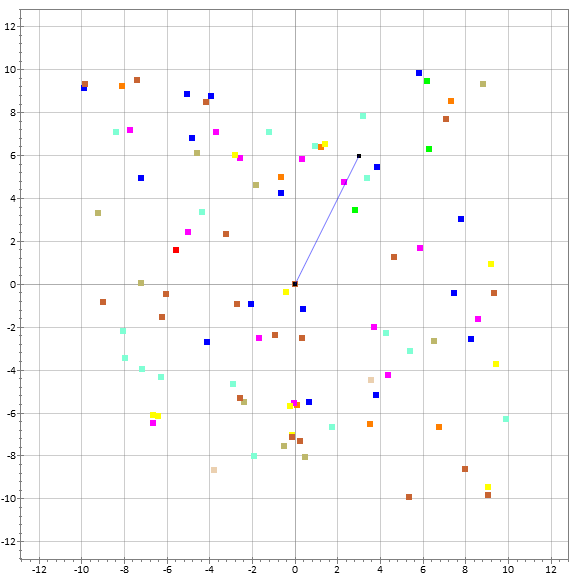}} & 
\subfloat[$\BGMldefault$ error]{\includegraphics[scale=0.15]{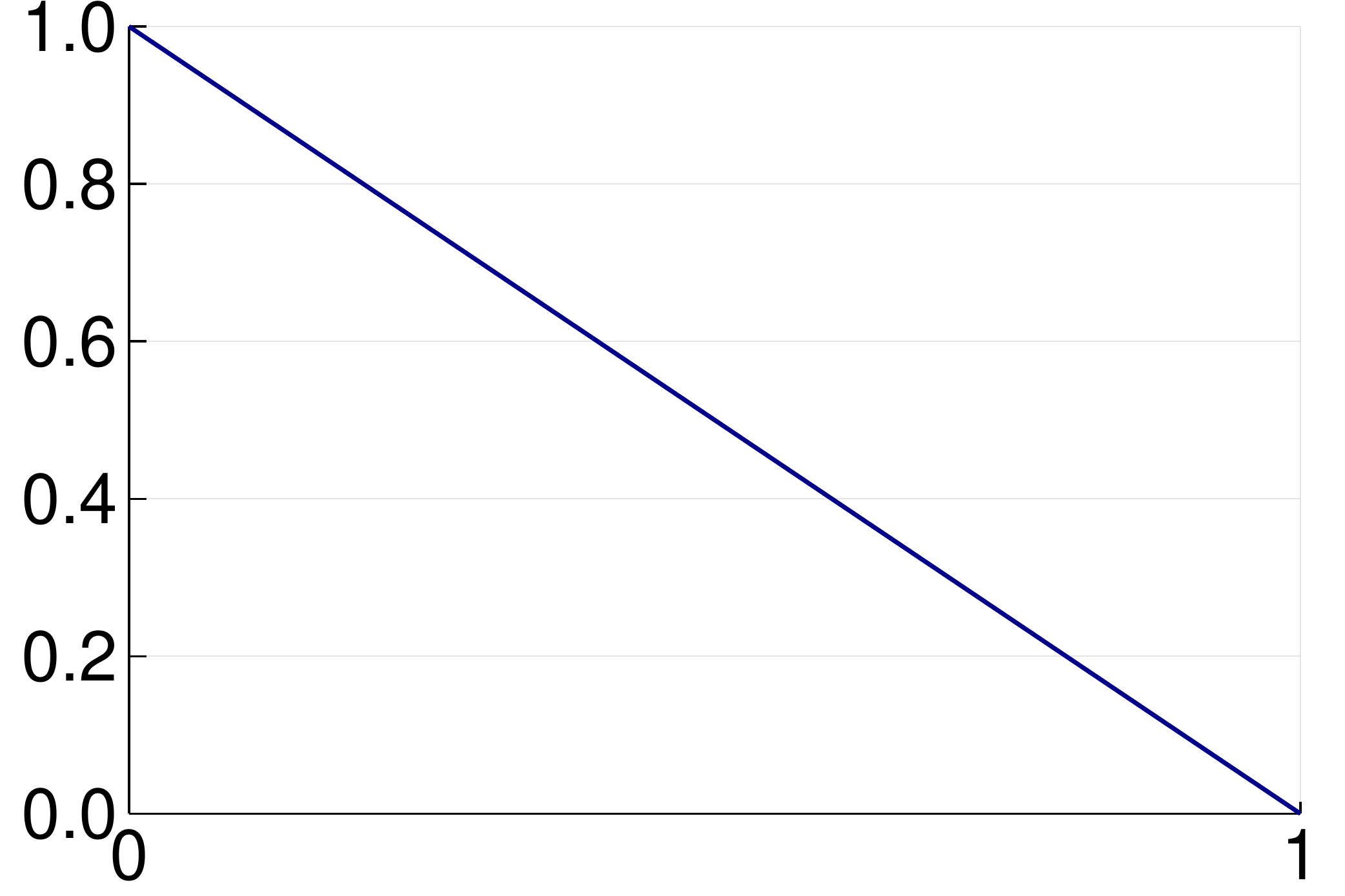}} & 
\subfloat[$\BGMlbest$ orbit]{\includegraphics[scale=0.15]{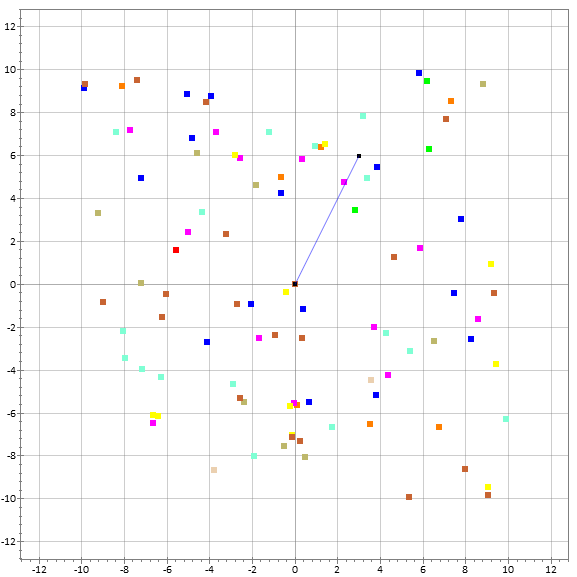}} &
\subfloat[$\BGMlbest$ error]{\includegraphics[scale=0.15]{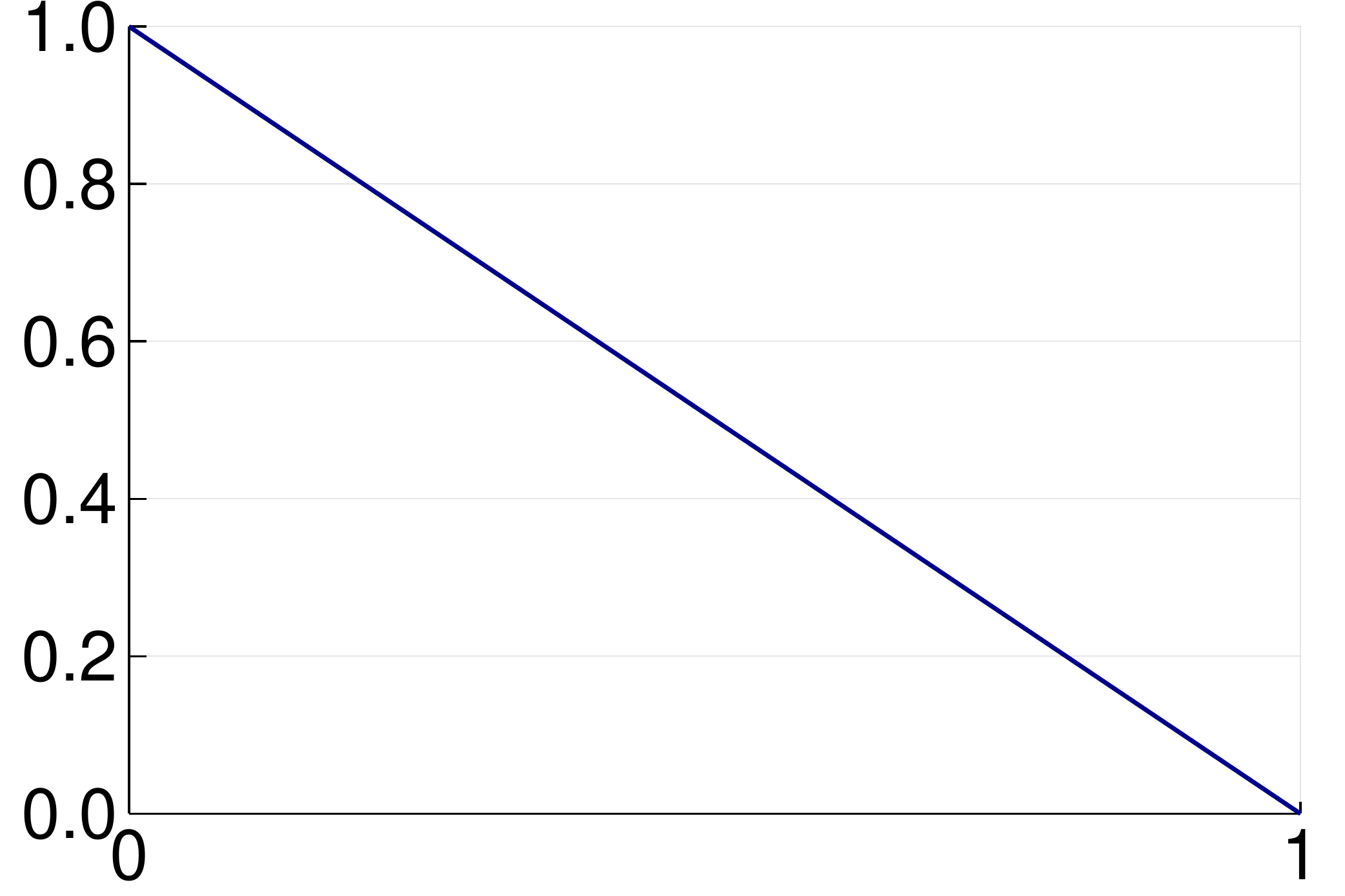}}
\\
& & & \\ 
\large \textbf{ExParP} & & & \\
\subfloat[$\BGMldefault$ orbit]{\includegraphics[scale=0.15]{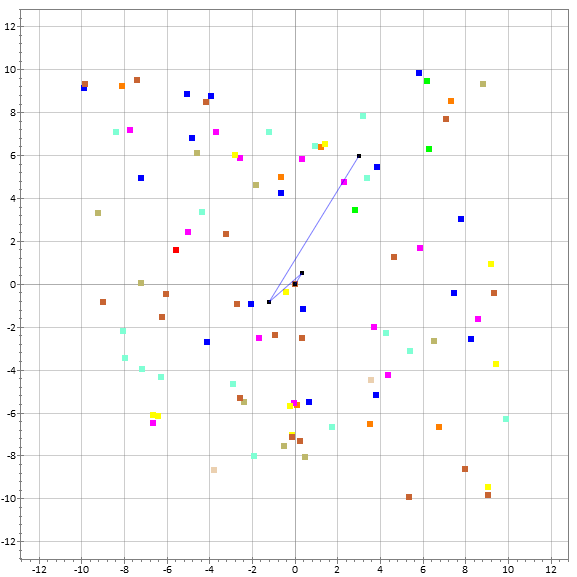}} & 
\subfloat[$\BGMldefault$ error]{\includegraphics[scale=0.15]{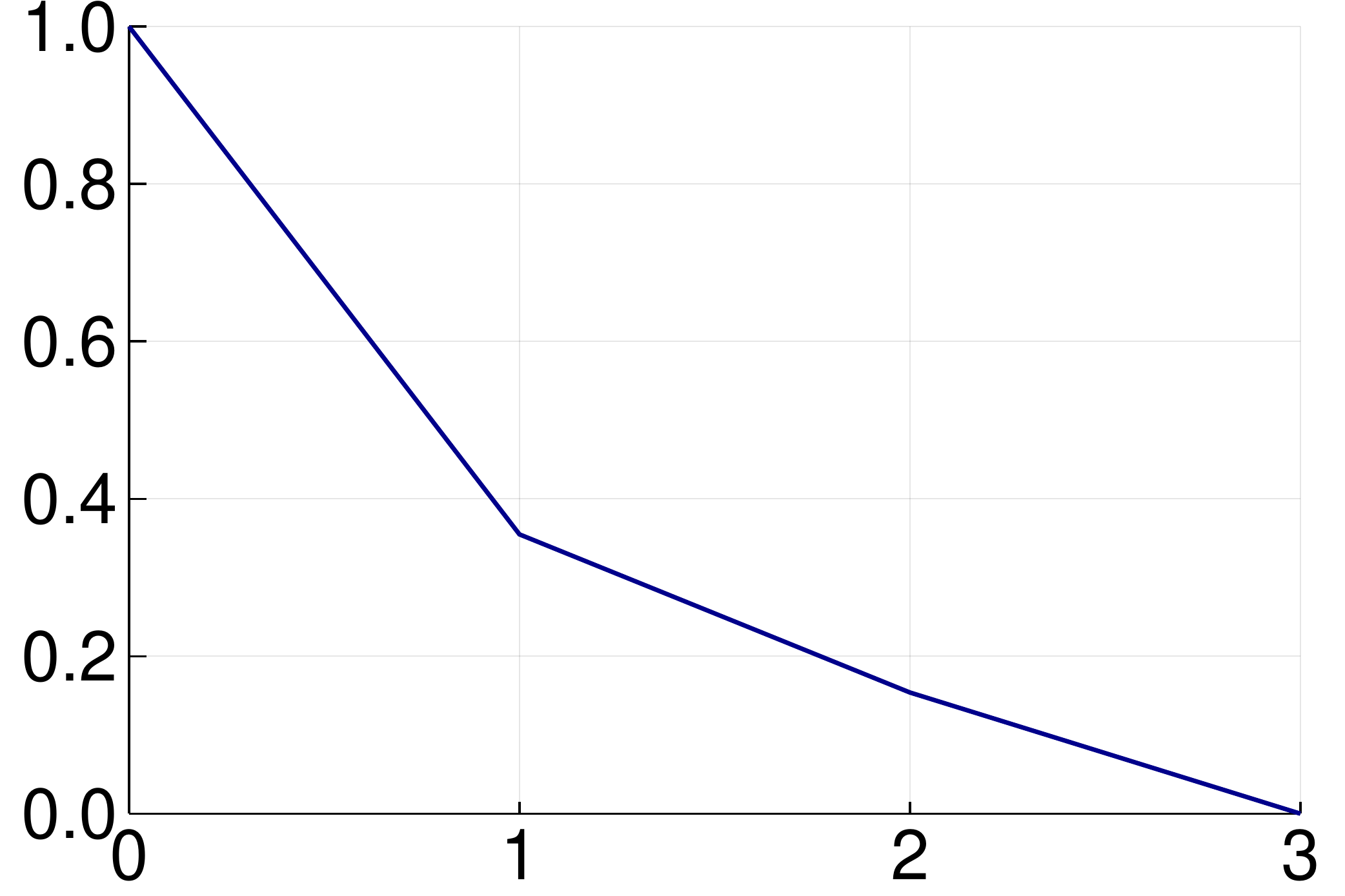}} & 
\subfloat[$\BGMlbest$ orbit]{\includegraphics[scale=0.15]{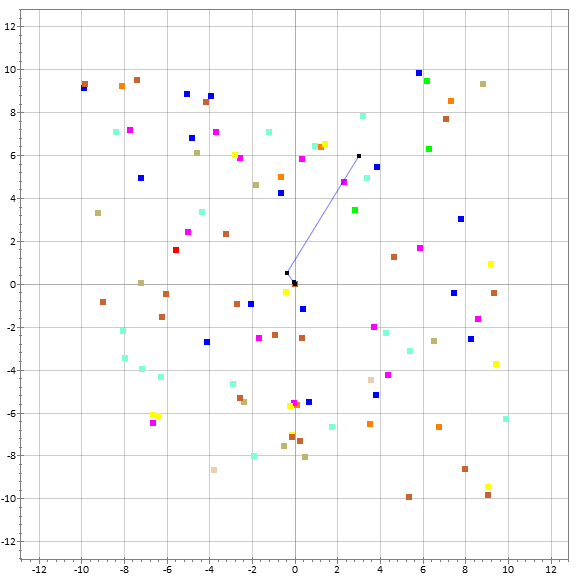}} &
\subfloat[$\BGMlbest$ error]{\includegraphics[scale=0.15]{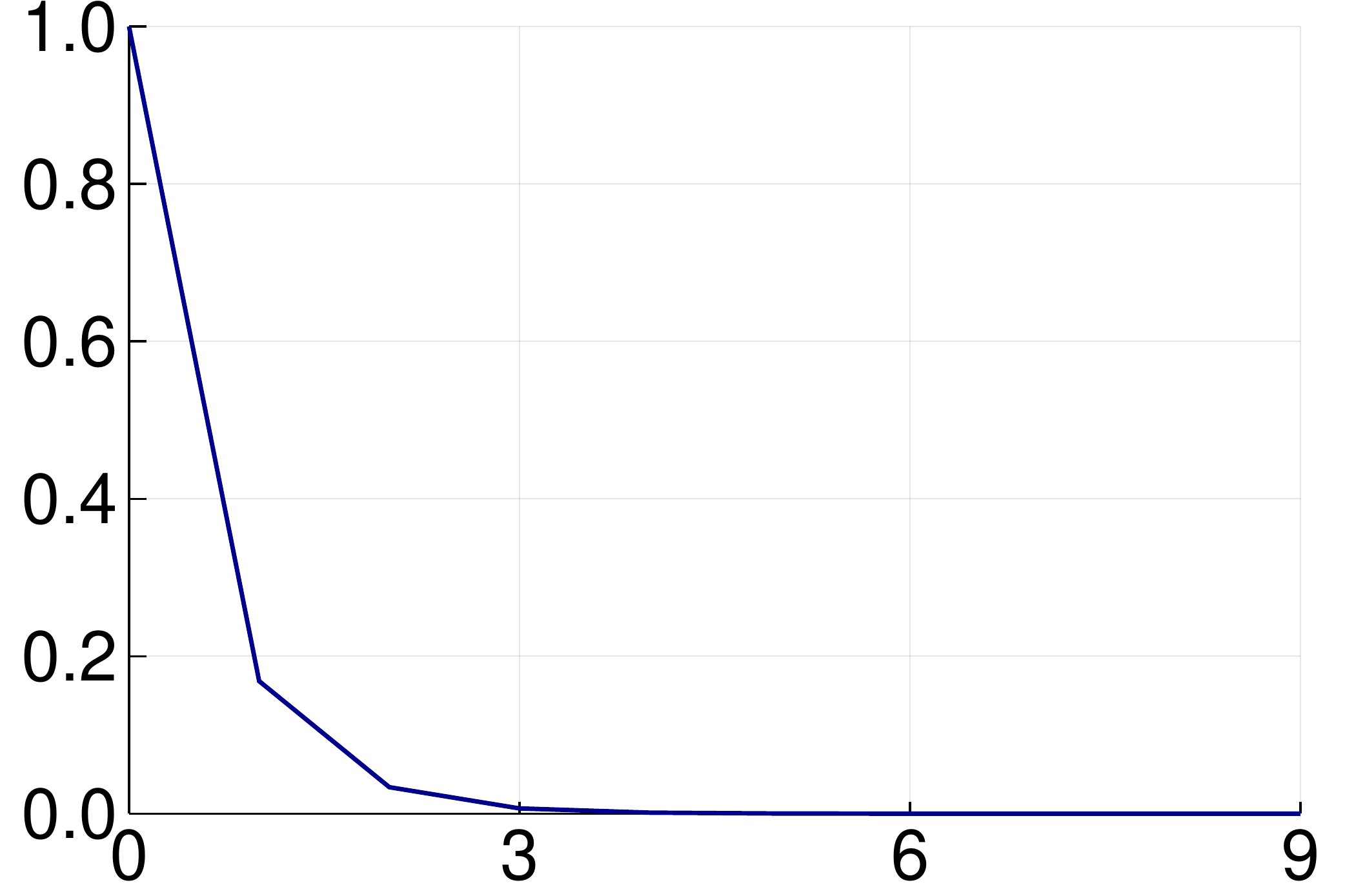}}
\\
& & & \\ 
\large \textbf{DR} & & & \\
\subfloat[$\BGMldefault$ orbit]{\includegraphics[scale=0.15]{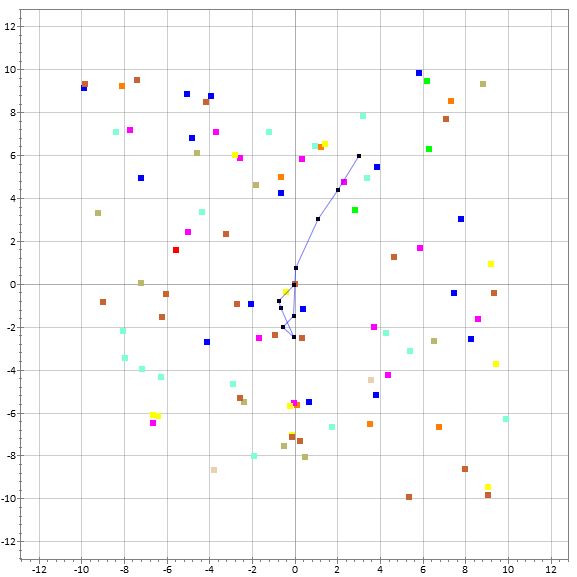}} & 
\subfloat[$\BGMldefault$ error]{\includegraphics[scale=0.15]{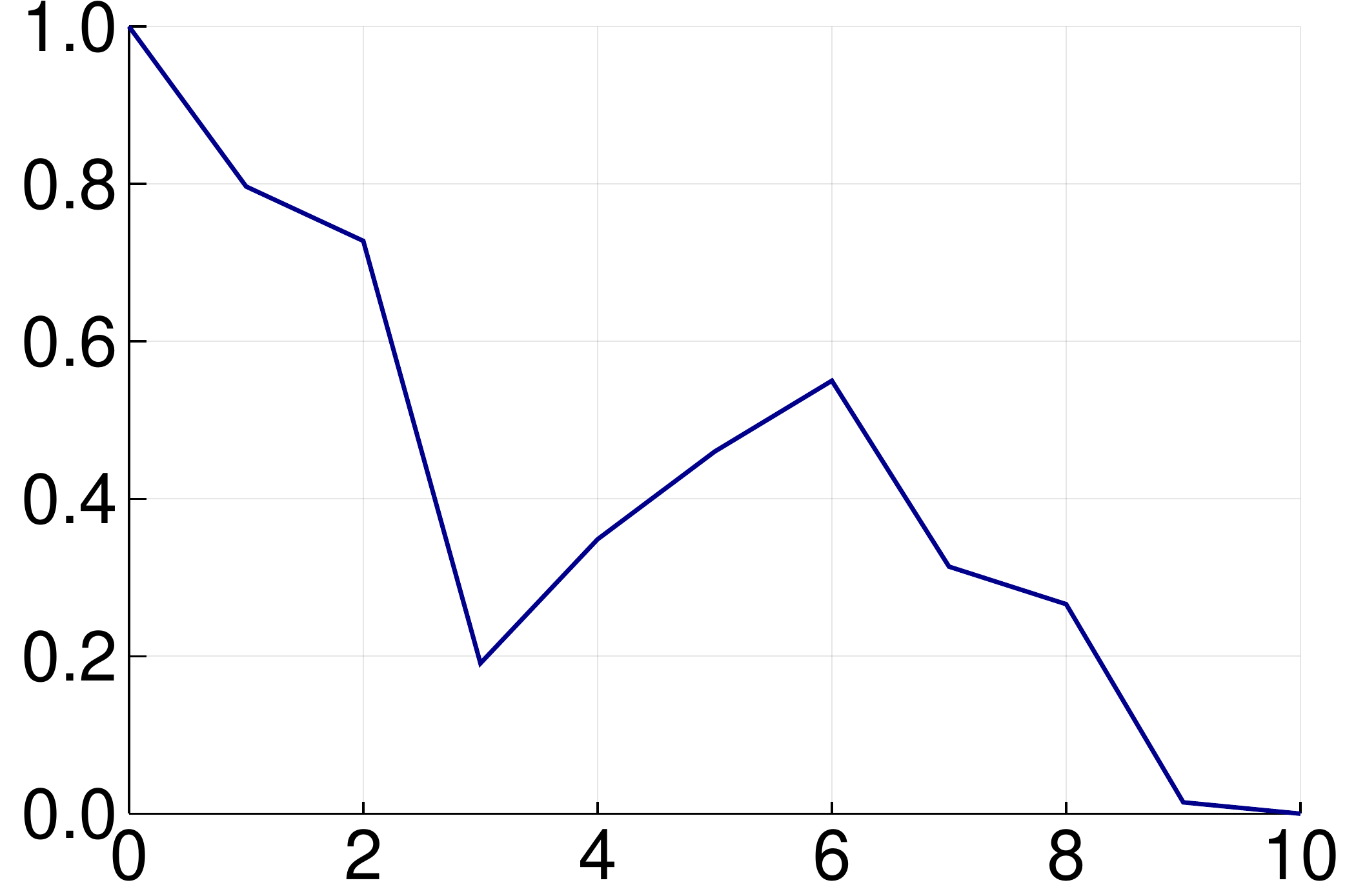}} & 
\subfloat[$\BGMlbest$ orbit]{\includegraphics[scale=0.15]{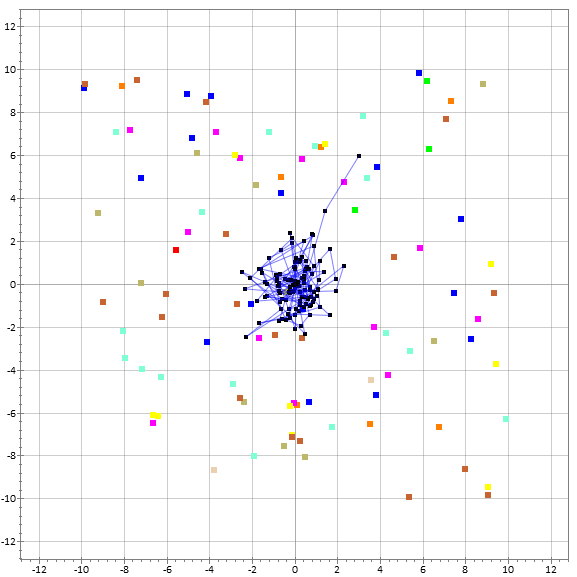}} &
\subfloat[$\BGMlbest$ error]{\includegraphics[scale=0.15]{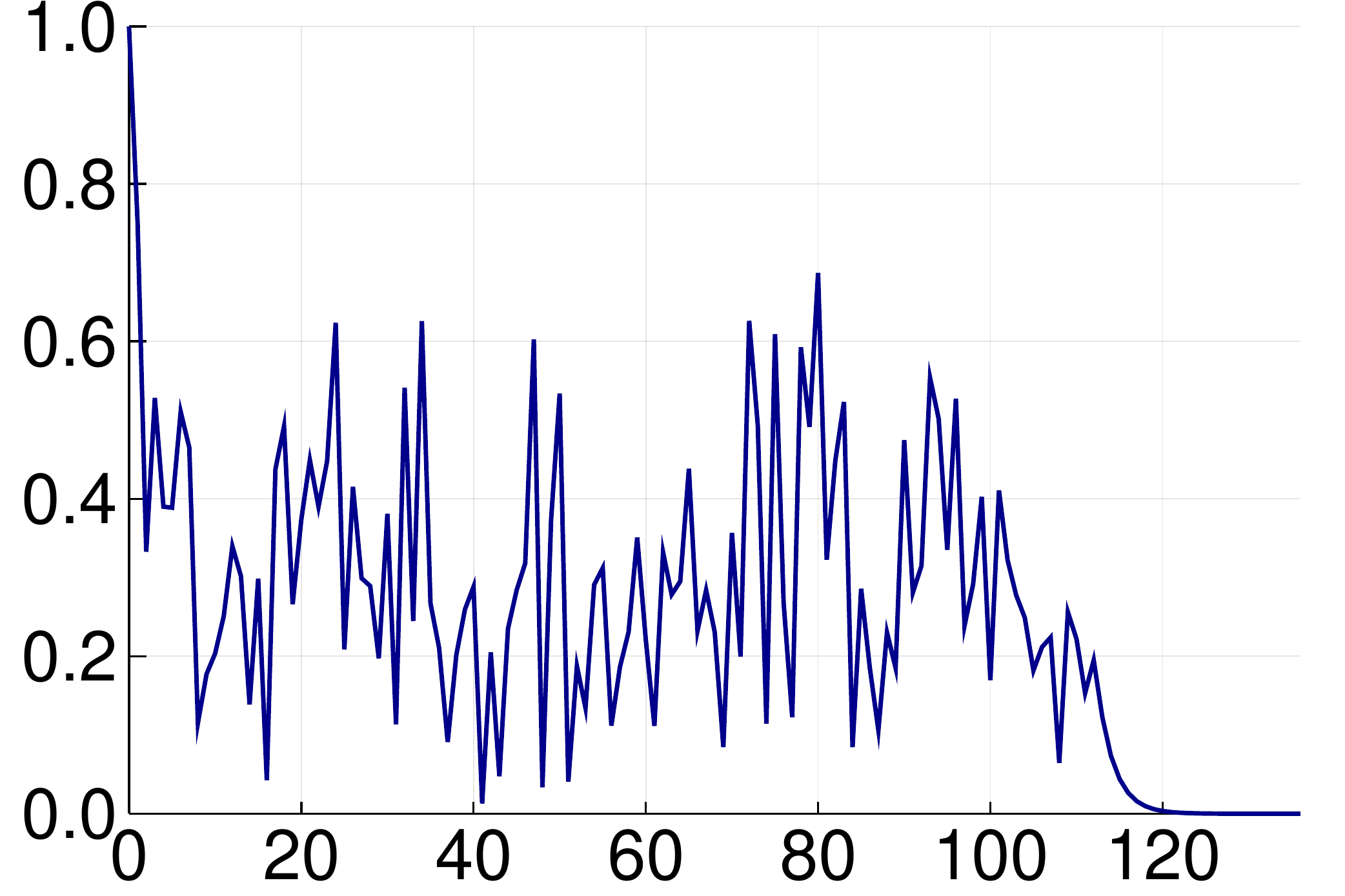}}
\\
& & & \\ 
\large \textbf{CycDR} & & & \\
\subfloat[$\BGMldefault$ orbit]{\includegraphics[scale=0.15]{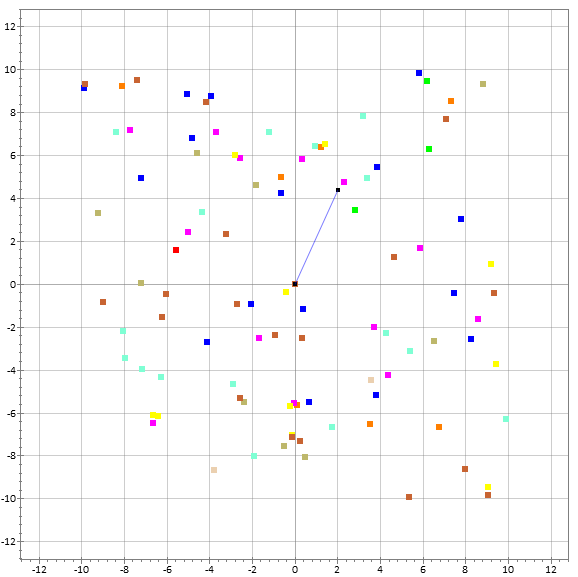}} & 
\subfloat[$\BGMldefault$ error]{\includegraphics[scale=0.15]{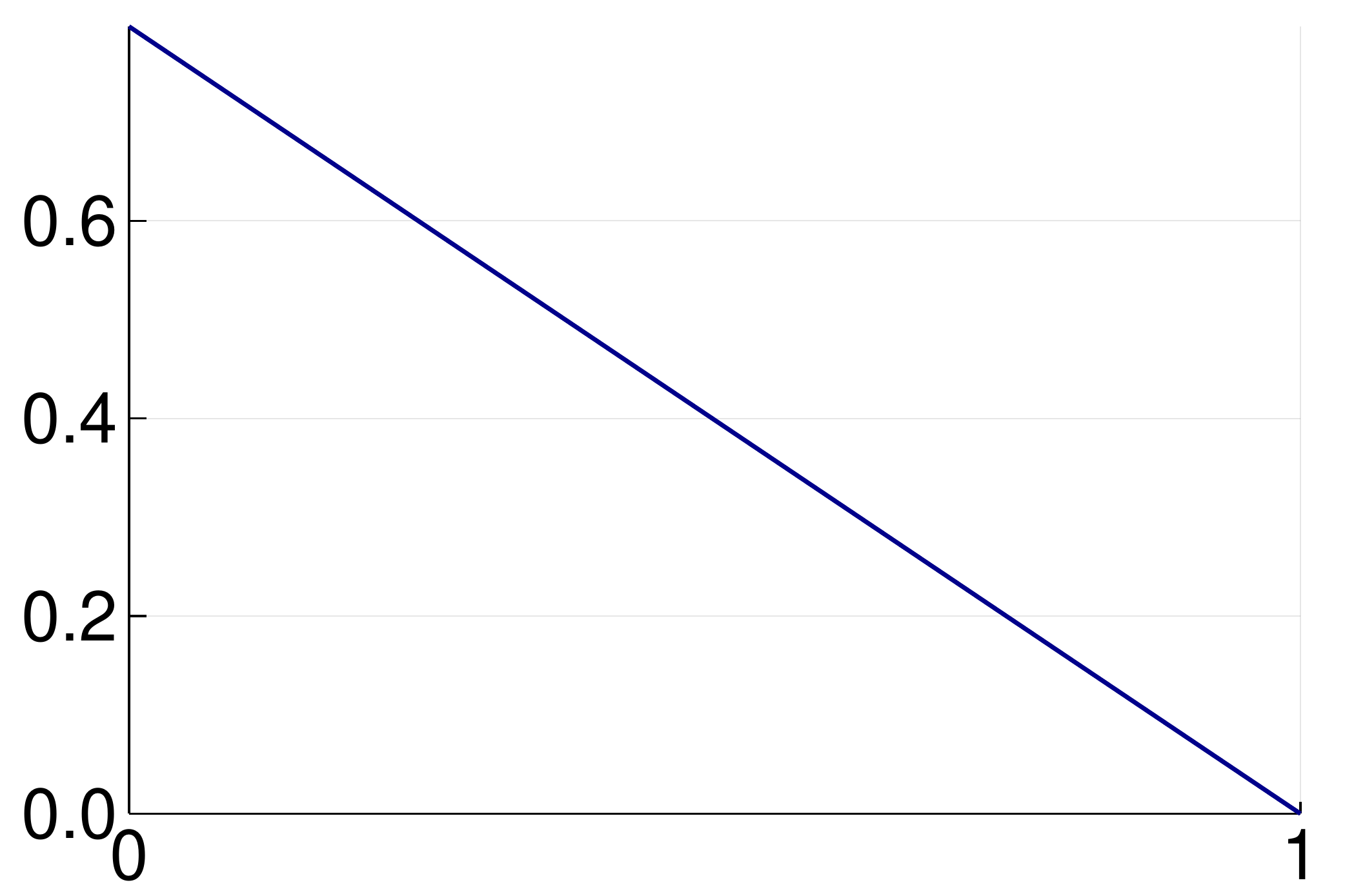}} & 
\subfloat[$\BGMlbest$ orbit]{\includegraphics[scale=0.15]{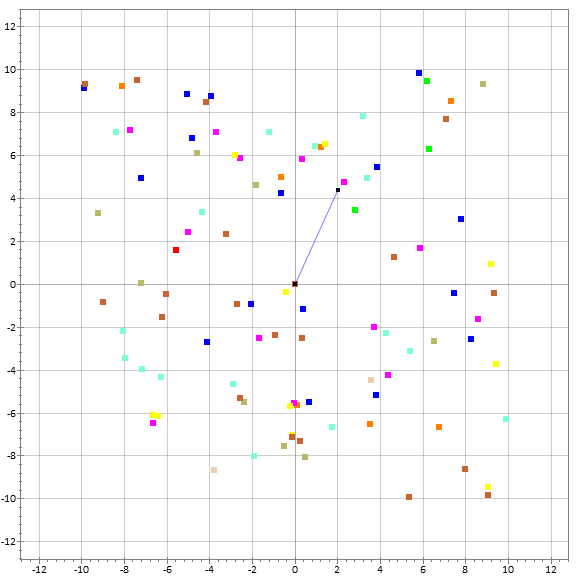}} &
\subfloat[$\BGMlbest$ error]{\includegraphics[scale=0.15]{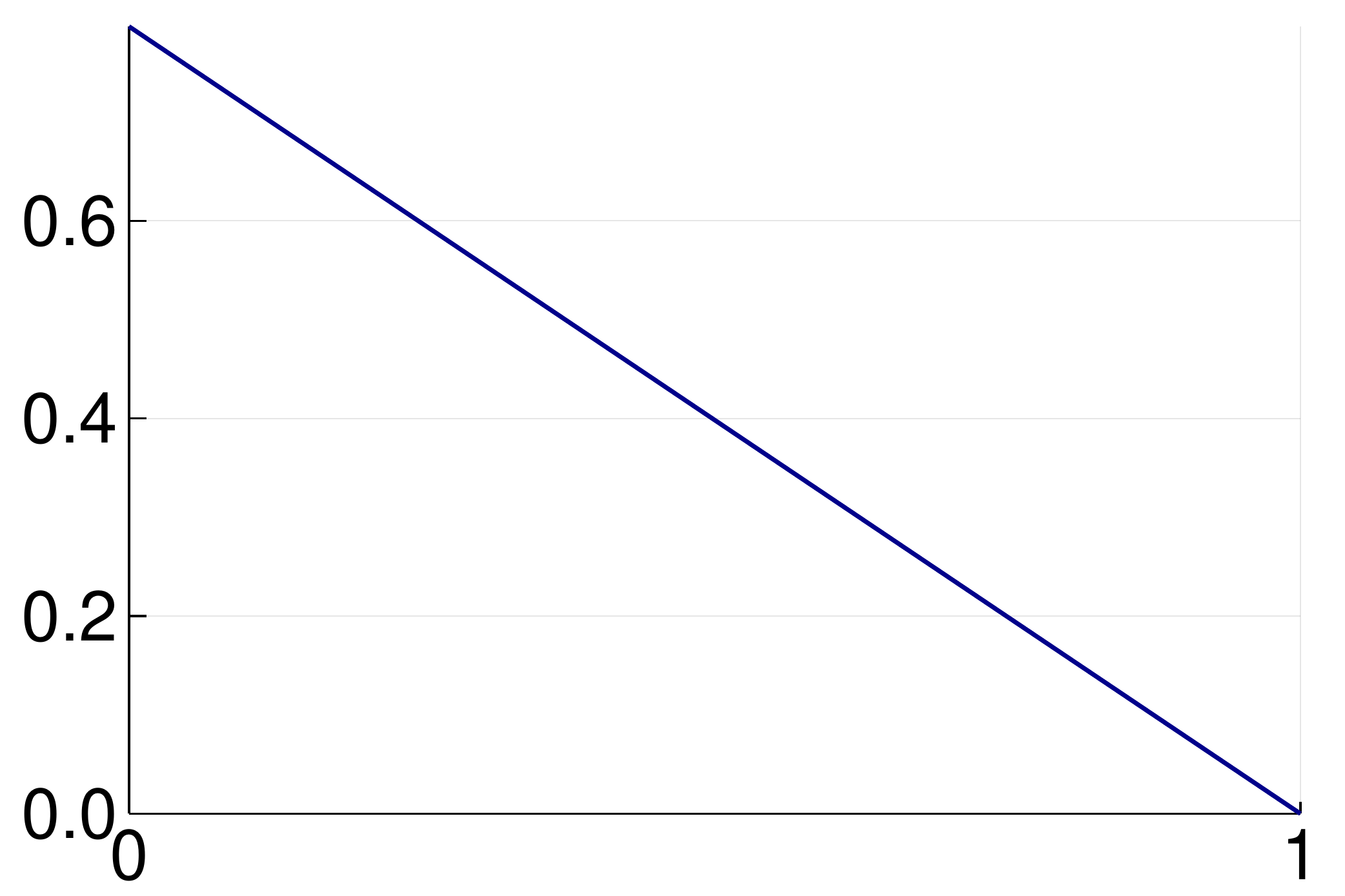}}
\end{tabular}
\caption{Orbits and errors for CycP, ExParP, DR, and CycDR in the many sets with few points constellation}
\end{figure}

\subsection{Many sets with many points}

\BGMvspace{0.5}

\begin{figure}[H]
\begin{tabular}{cccc}
\large \textbf{CycP} & & & \\
\subfloat[$\BGMldefault$ orbit]{\includegraphics[scale=0.15]{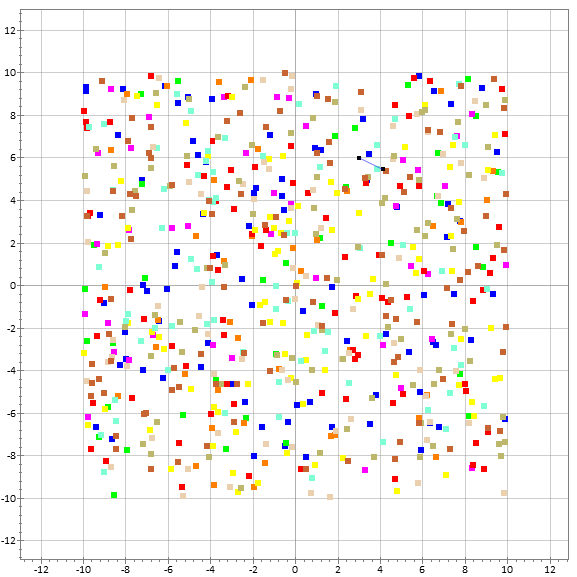}} & 
\subfloat[$\BGMldefault$ error]{\includegraphics[scale=0.15]{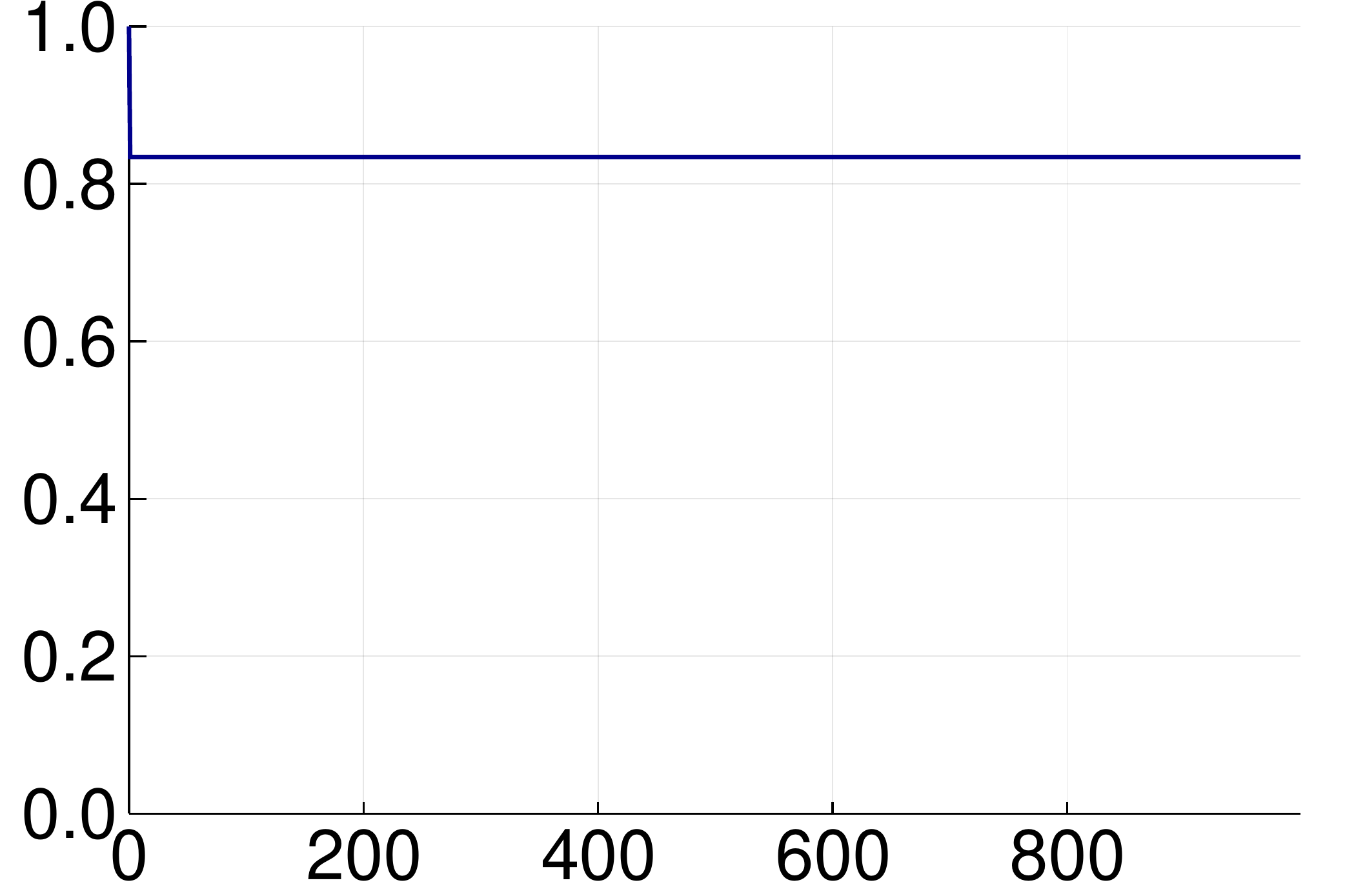}} & 
\subfloat[$\BGMlbest$ orbit]{\includegraphics[scale=0.15]{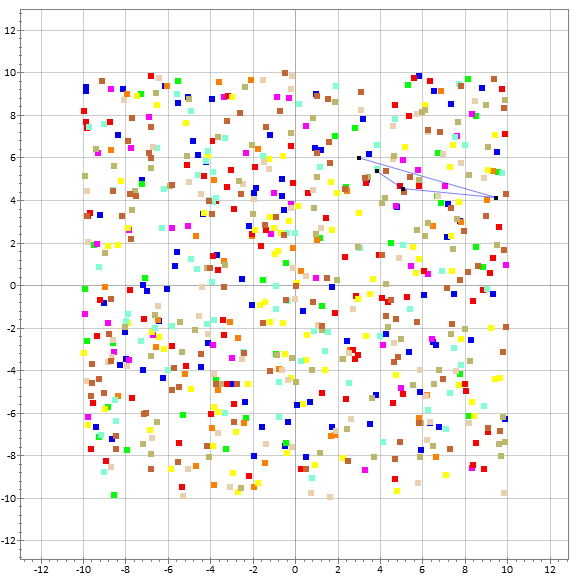}} &
\subfloat[$\BGMlbest$ error]{\includegraphics[scale=0.15]{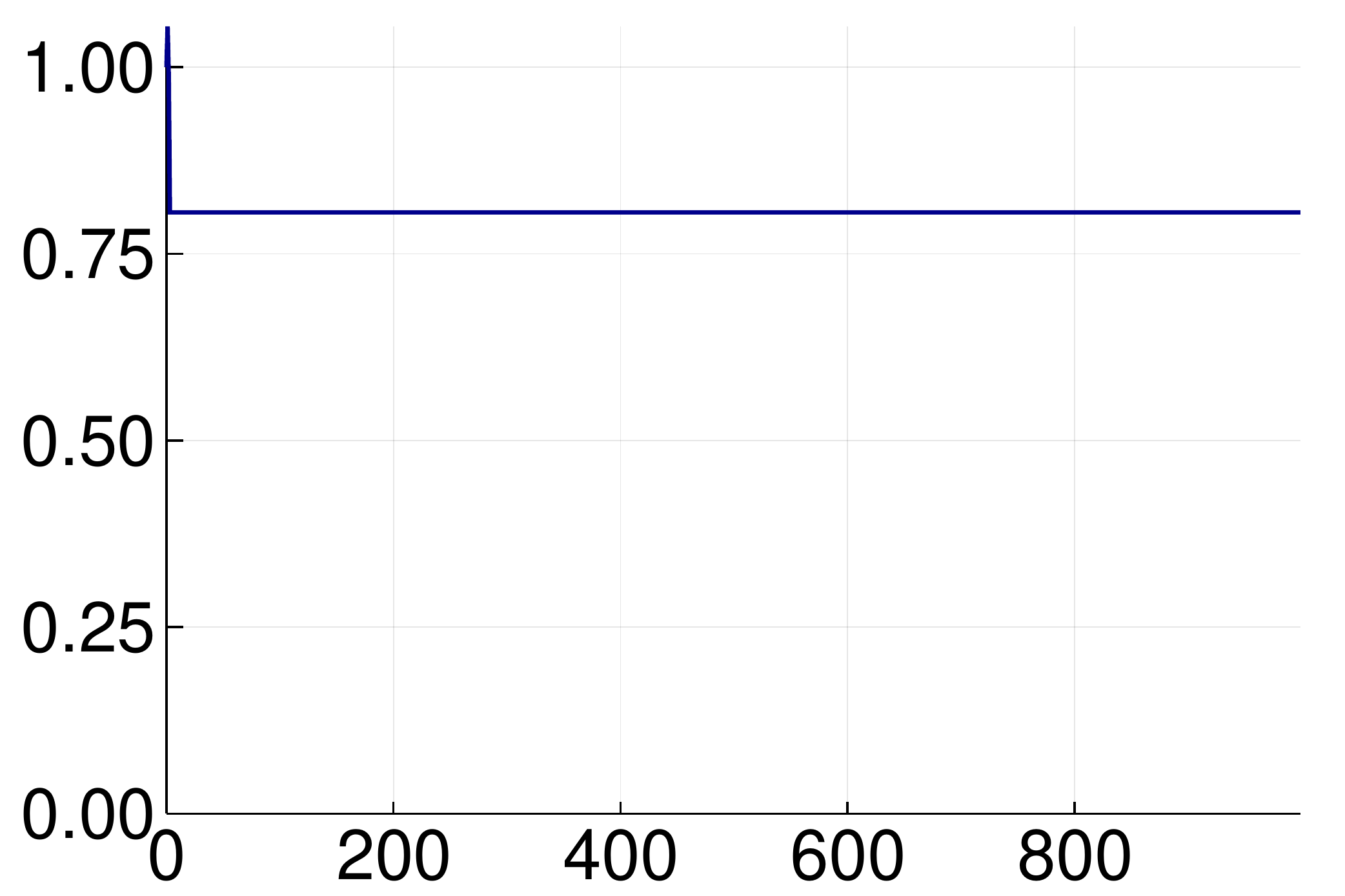}}
\\
& & & \\ 
\large \textbf{ExParP} & & & \\
\subfloat[$\BGMldefault$ orbit]{\includegraphics[scale=0.15]{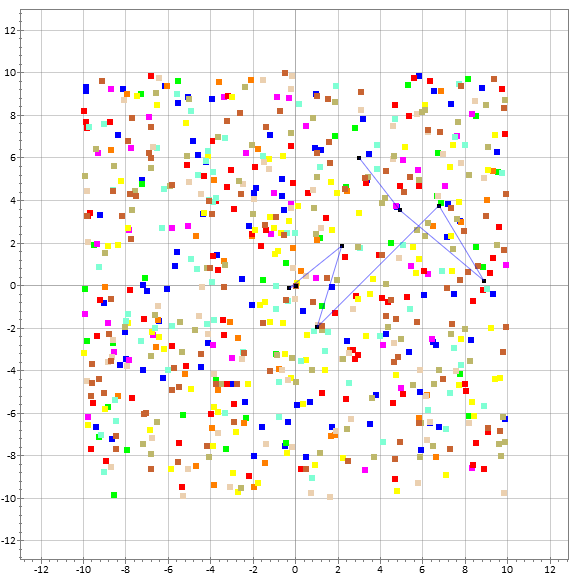}} & 
M\subfloat[$\BGMldefault$ error]{\includegraphics[scale=0.15]{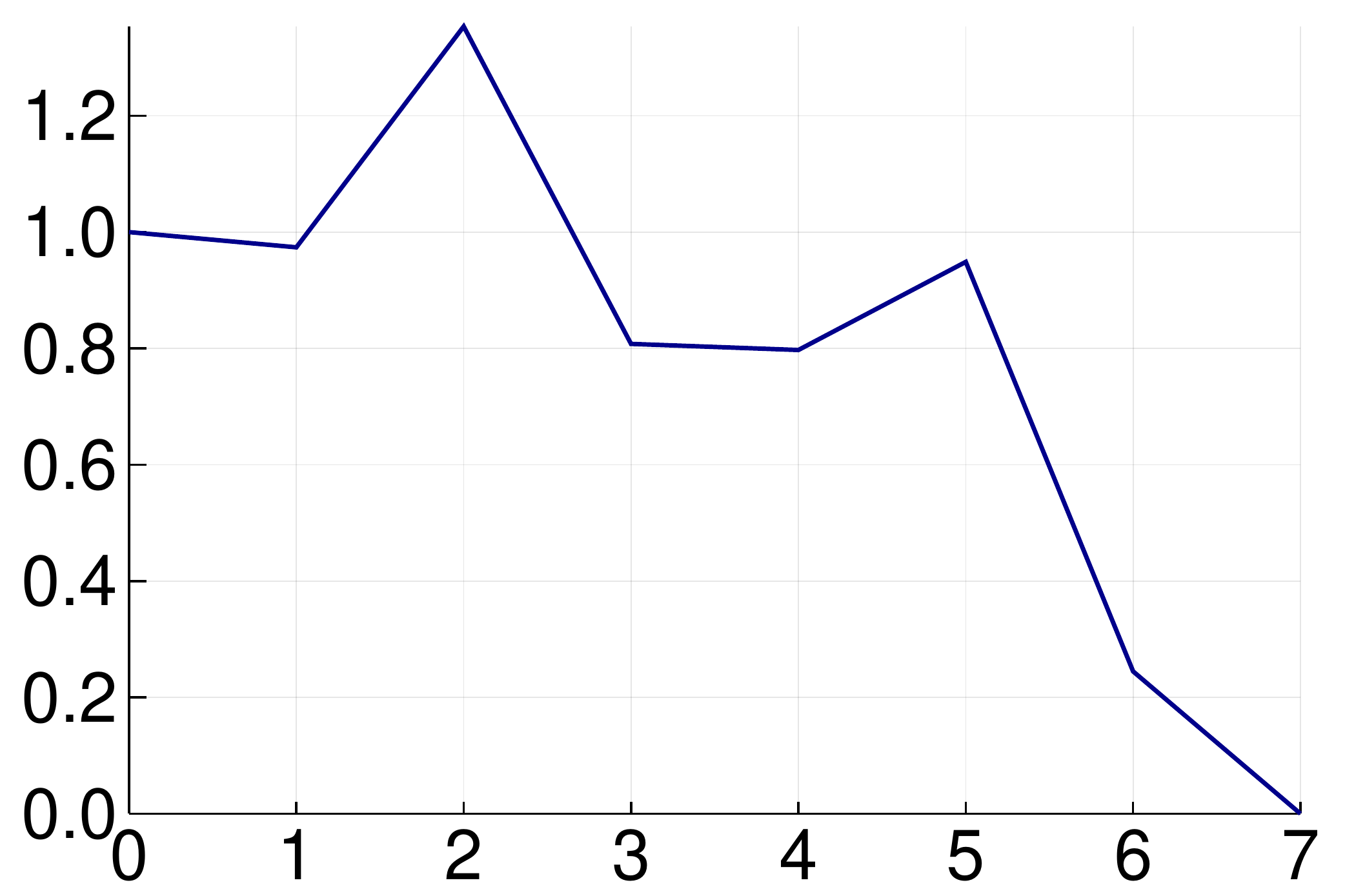}} & 
\subfloat[$\BGMlbest$ orbit]{\includegraphics[scale=0.15]{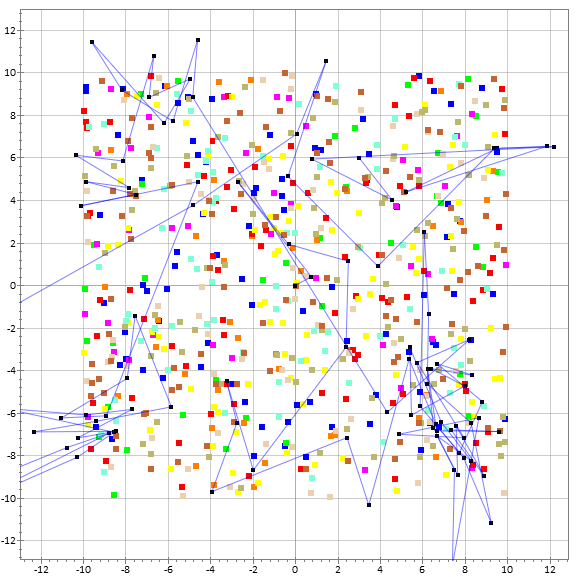}} &
\subfloat[$\BGMlbest$ error]{\includegraphics[scale=0.15]{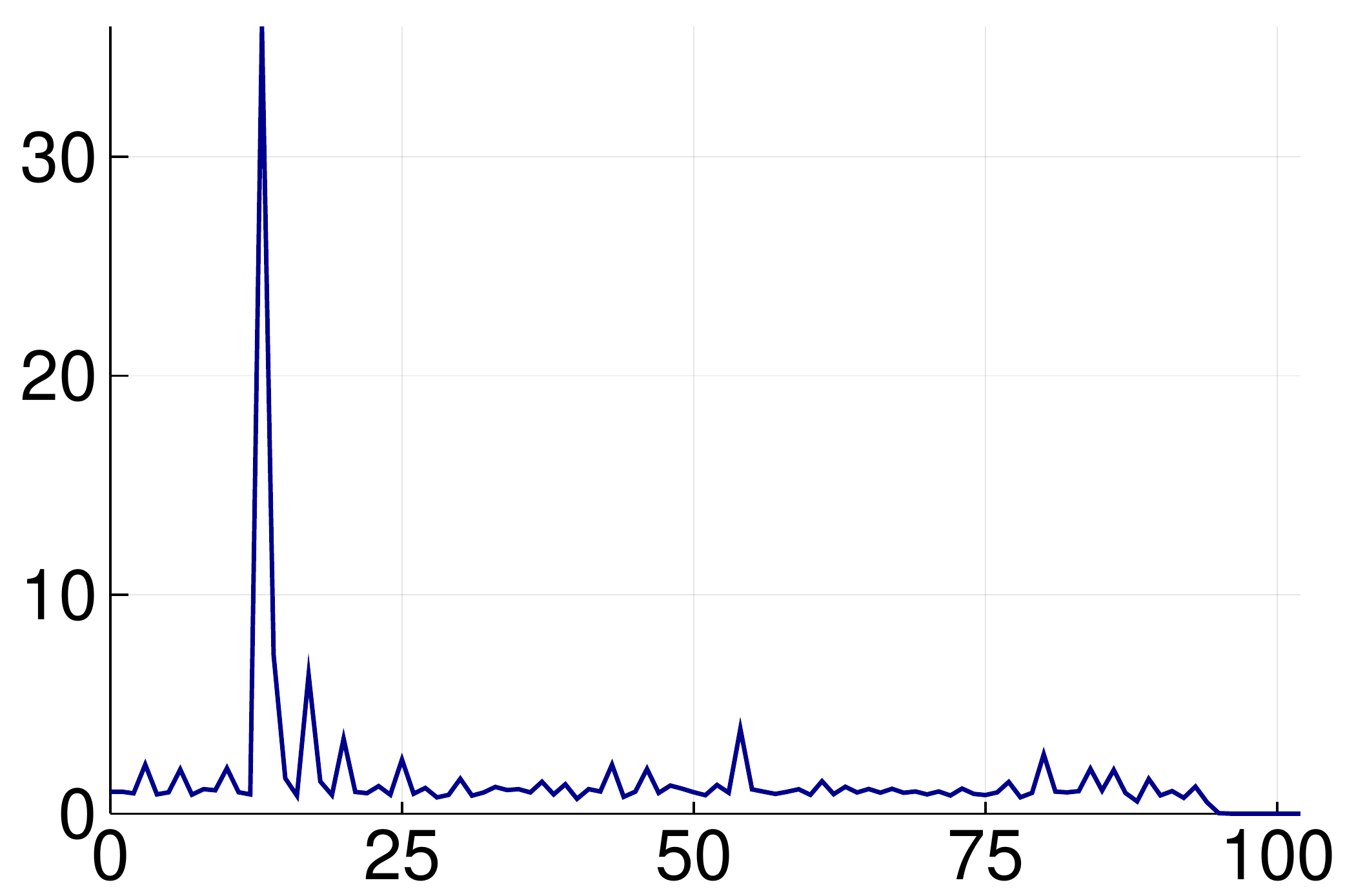}}
\\
& & & \\
\large \textbf{DR} & & & \\
\subfloat[$\BGMldefault$ orbit]{\includegraphics[scale=0.15]{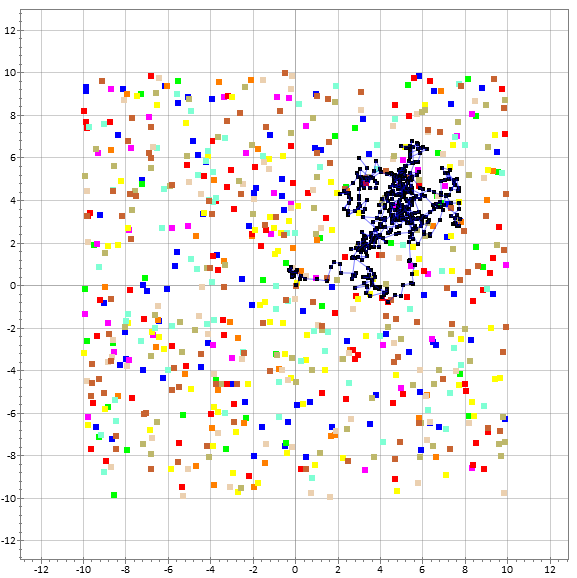}} & 
\subfloat[$\BGMldefault$ error]{\includegraphics[scale=0.15]{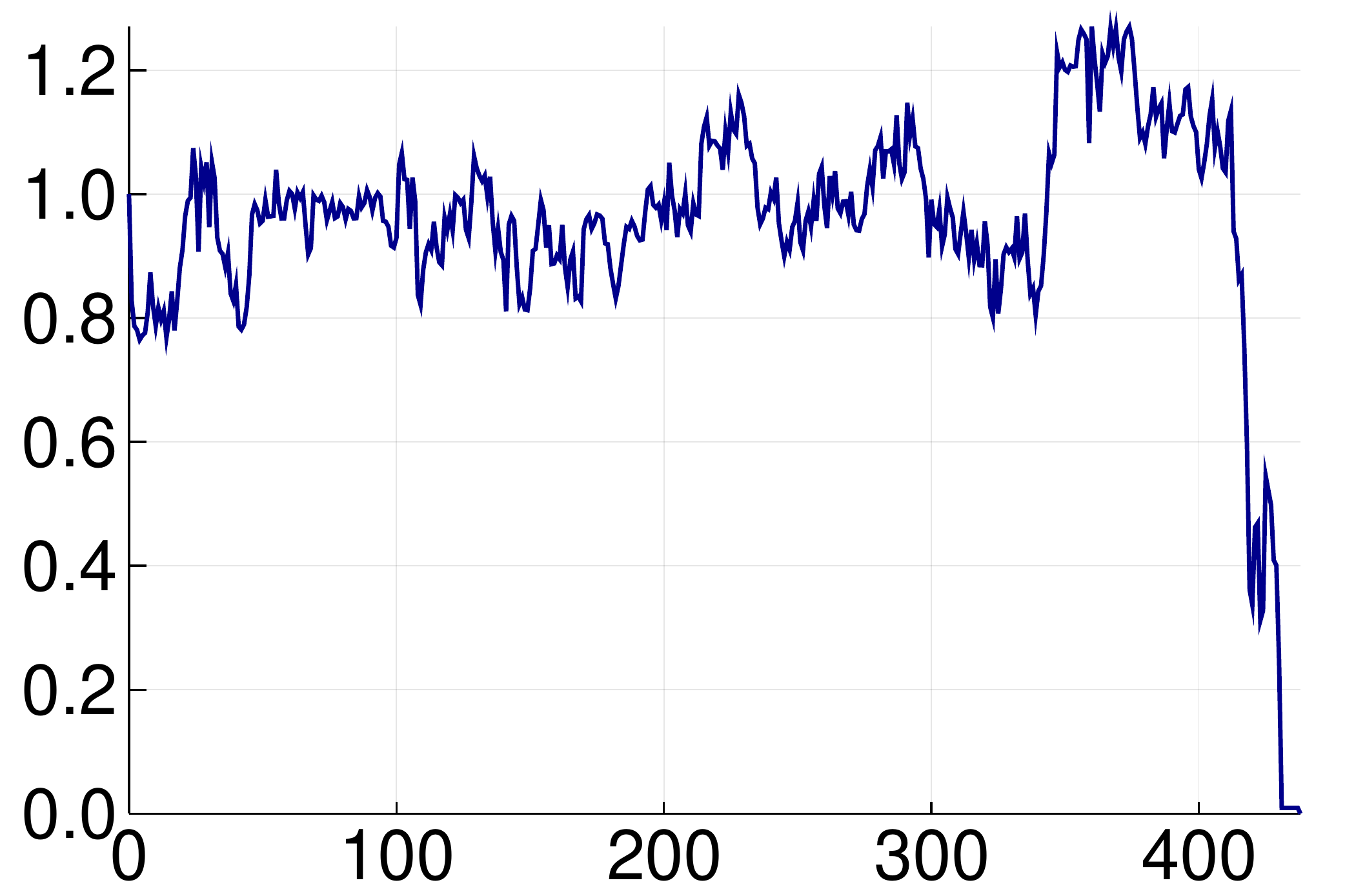}} & 
\subfloat[$\BGMlbest$ orbit]{\includegraphics[scale=0.15]{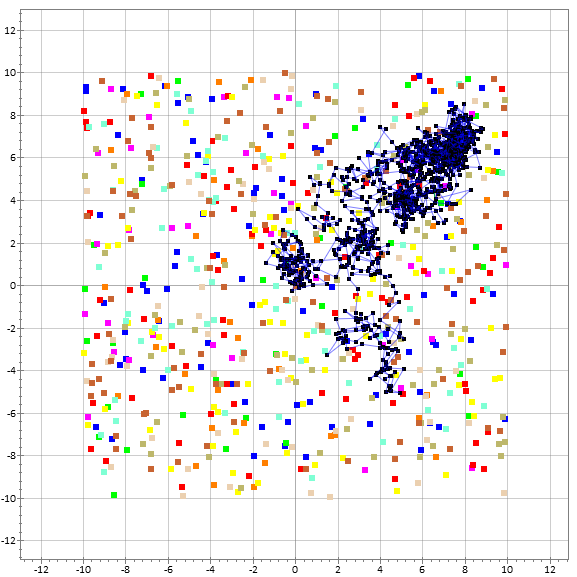}} &
\subfloat[$\BGMlbest$ error]{\includegraphics[scale=0.15]{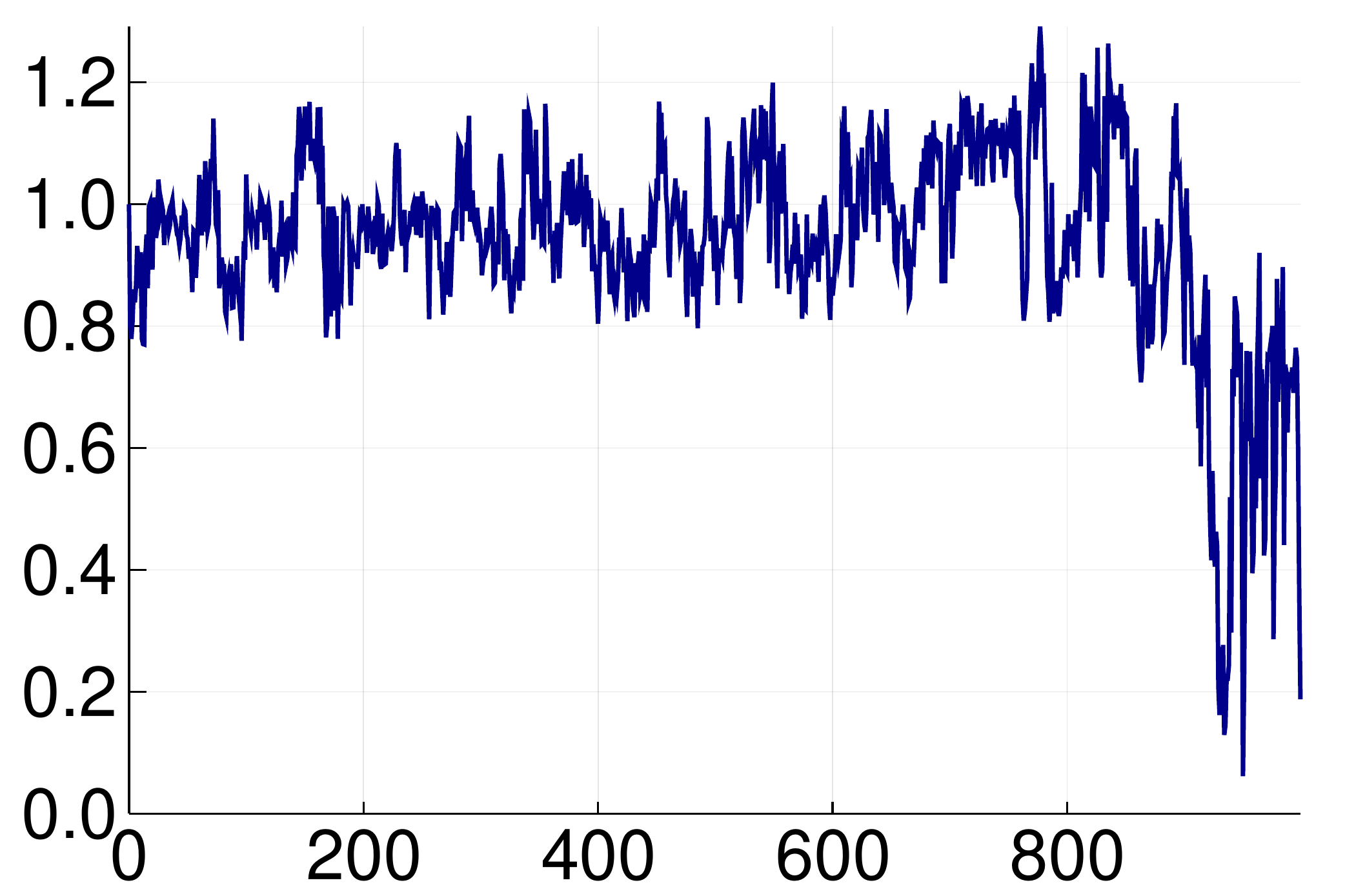}}
\\
& & & \\
\large \textbf{CycDR} & & & \\
\subfloat[$\BGMldefault$ orbit]{\includegraphics[scale=0.15]{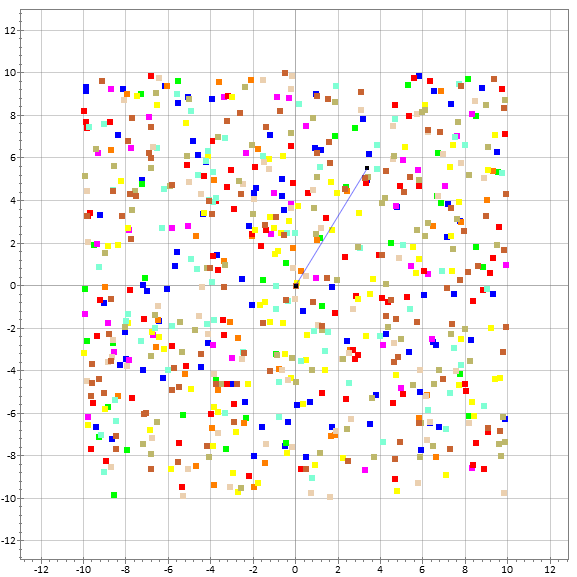}} & 
\subfloat[$\BGMldefault$ error]{\includegraphics[scale=0.15]{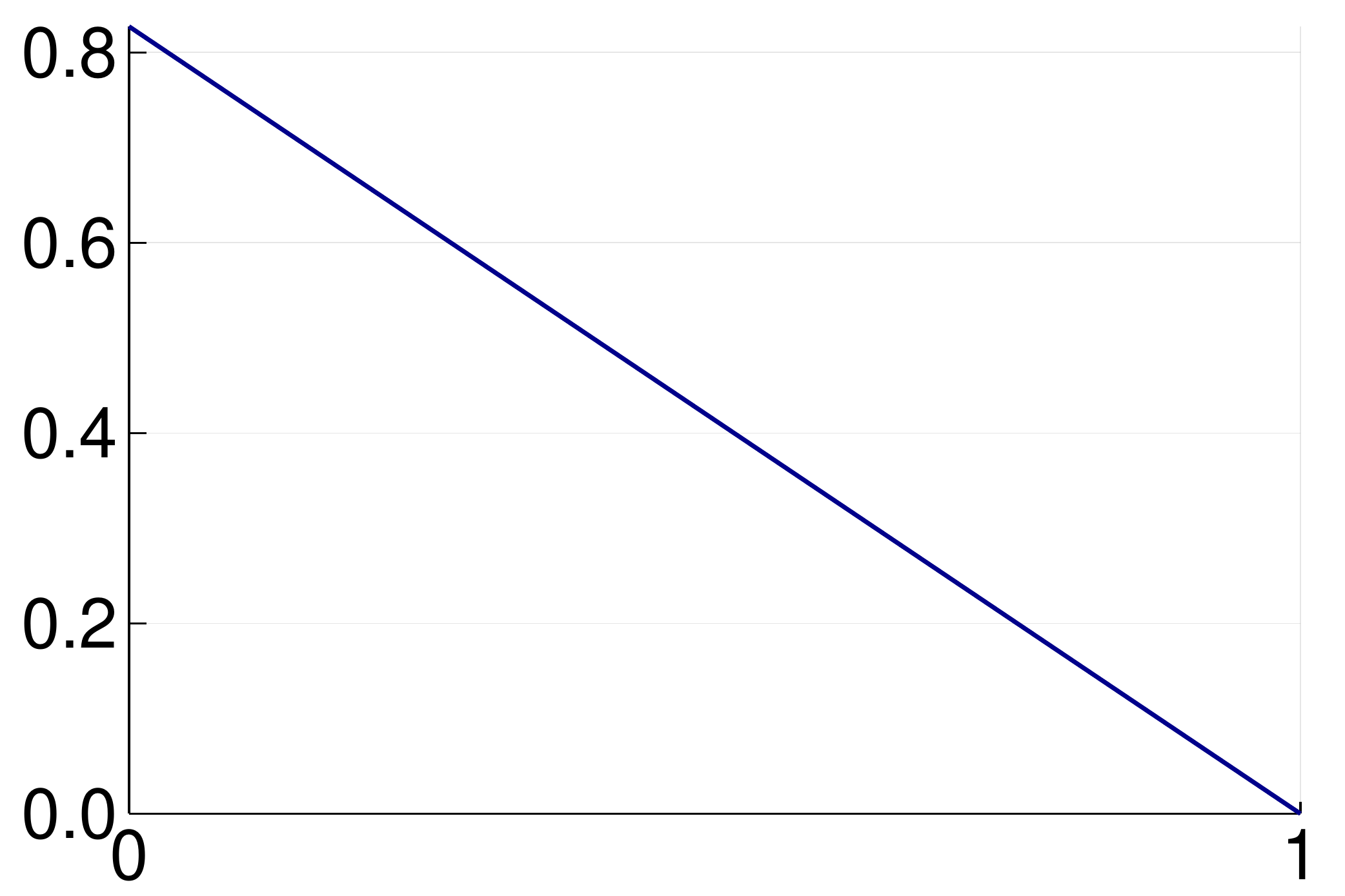}} & 
\subfloat[$\BGMlbest$ orbit]{\includegraphics[scale=0.15]{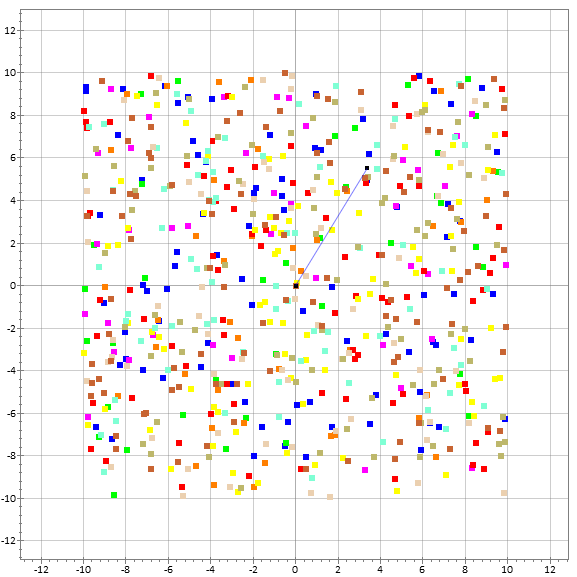}} &
\subfloat[$\BGMlbest$ error]{\includegraphics[scale=0.15]{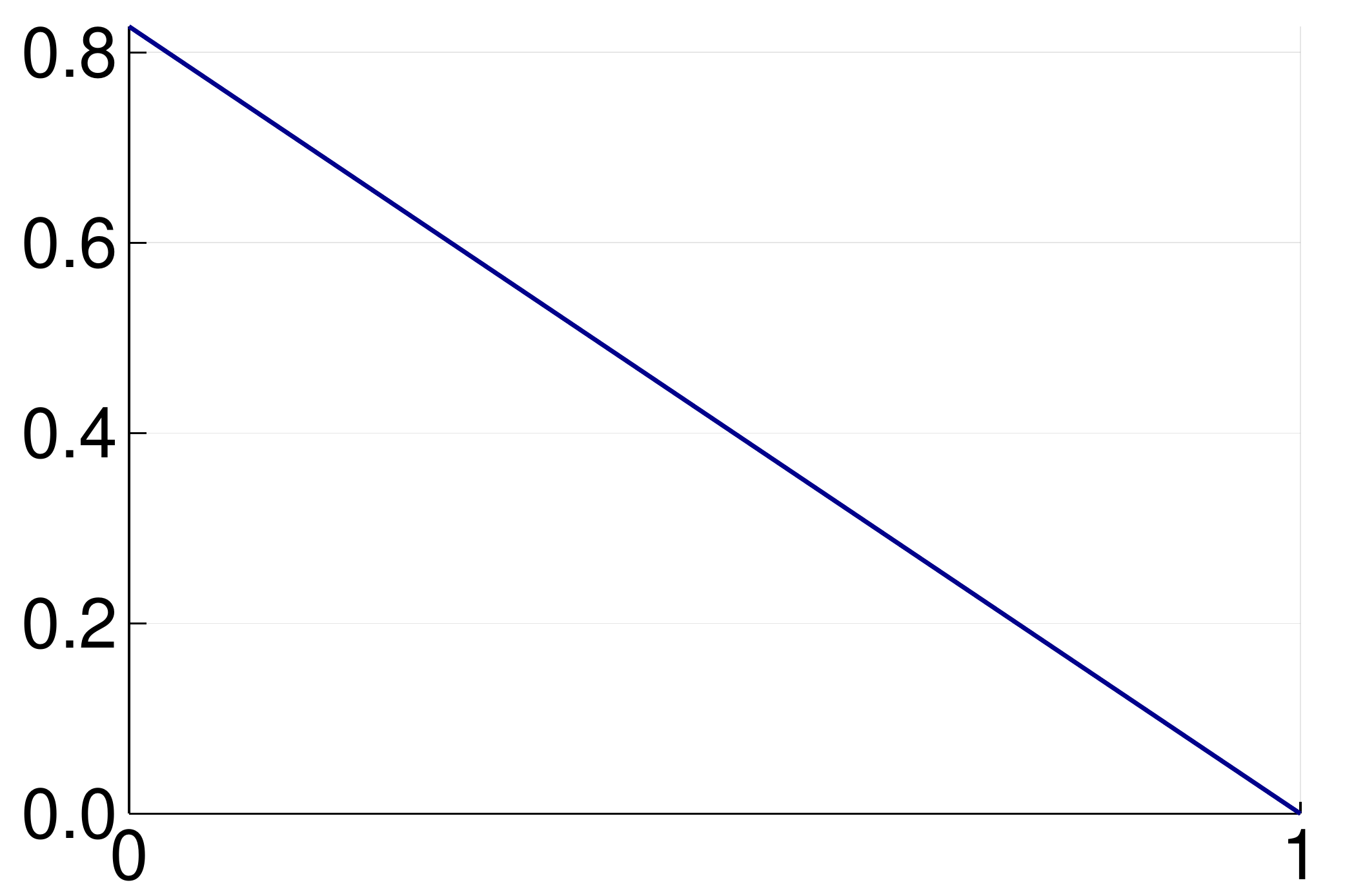}}
\end{tabular}
\caption{Orbits and errors for CycP, ExParP, DR, and CycDR in the many sets with many points constellation}
\end{figure} 

\subsection{Discussion}

The numerical results in this subsection suggest the following:
The most challenging constellation is the one with few sets and many points.
The least challenging constellation is the one with many sets and few points
for which all algorithms are successful.

The worst algorithm is CycP.
ExParP with $\BGMlbest$ solves all four constellations. 
DR solves all constellations but $\lambda$ has to be chosen appropriately.
CycDR works well with $\BGMldefault$ in terms of number of iterations
required; however, it was not able to solve the constellation with
few sets and many points. 

The experiments in this section suggest that 
(i) ExParP, DR, and CycDR are algorithms worthwhile exploring and that 
(ii) experimenting with $\lambda$ may lead to improved convergence.

Because the results in this section feature a \emph{fixed} starting point,
we will explore in the next section the four constellations for a \emph{multitude} of
starting points.

\section{Local and global behaviour}

\label{BGMsec:locglob}

In this section, we continue to consider our four constellations (see Section~\ref{BGMsec:4con}) which our four algorithms attempt to solve
(see Section~\ref{BGMsec:4alg}).

In contrast to Section~\ref{BGMsec:track} where we tracked a single orbit,
we here illustrate 
\emph{local} and \emph{global} behaviour of the algorithms for a multitude 
of starting points, sampled from $[-10,10]\times[-10,10]$ and
$[-100,100]\times[-100,100]$, respectively.
We do this for $\BGMldefault=1$ and for $\BGMlbest$
(see the table in Figure~\ref{BGMfig:lambdatable} in Section~\ref{BGMsec:det});

For each starting point in the given range, these plots display as gray levels the number of iterations the algorithm needed in its attempt to solve the problem represented by the given constellation. Black corresponds to the minimum number of iterations (zero), and white to the maximum number of iterations (1000). The latter is obtained when the algorithm is unsuccessful. Therefore, the darker the image, the better the performance.

To quantitatively assess the performance of each algorithm, 
success rates are also provided. 
These are obtained by dividing the number of times the algorithm is 
successful by the number of starting points used.

Each of these images was generated using at least 15 million starting points. 
Depending on the constellation, the time required to generate these pictures ranged between a few minutes to about 3 hours using a quad-core computer.

\subsection{Few sets with few points}

\BGMvspace{0.5}

\begin{figure}[H]
\begin{tabular}{cccc}
\large \textbf{CycP} & & & \\
\subfloat[\scriptsize $\BGMldefault$ local (57\%)]{\includegraphics[scale=0.165]{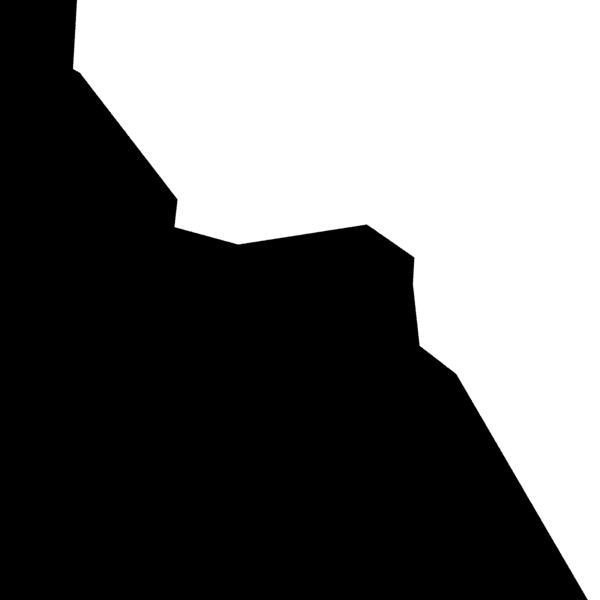}} & 
\subfloat[\scriptsize $\BGMldefault$ global (57\%)]{\includegraphics[scale=0.165]{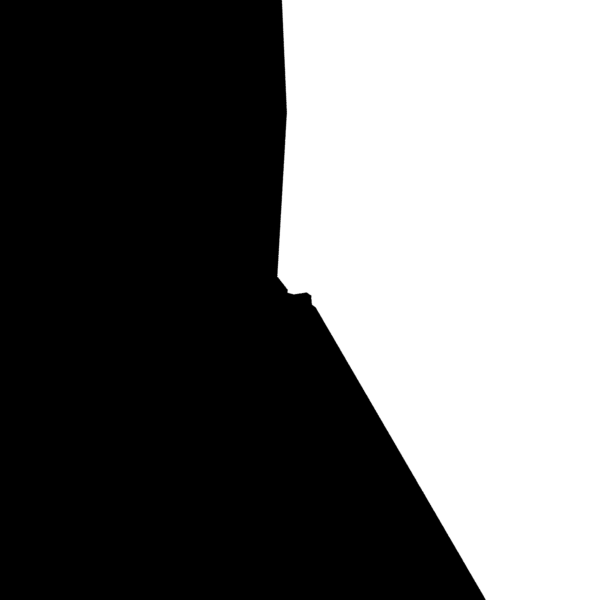}} & 
\subfloat[\scriptsize $\BGMlbest$ local (95\%)]{\includegraphics[scale=0.165]{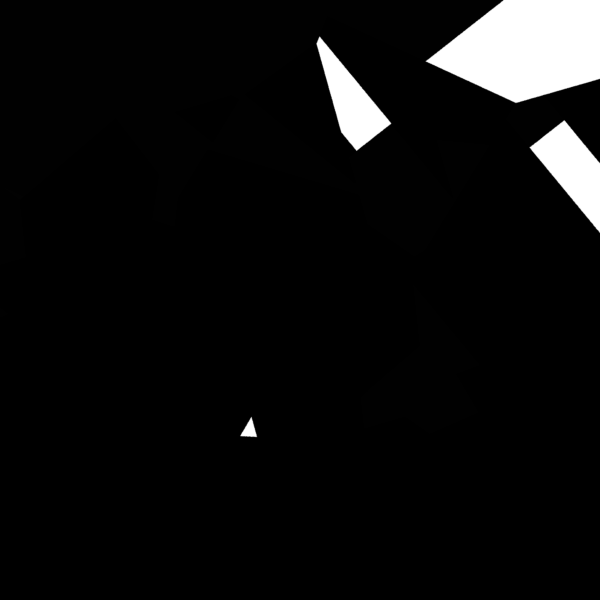}} &
\subfloat[\scriptsize $\BGMlbest$ global (98\%)]{\includegraphics[scale=0.165]{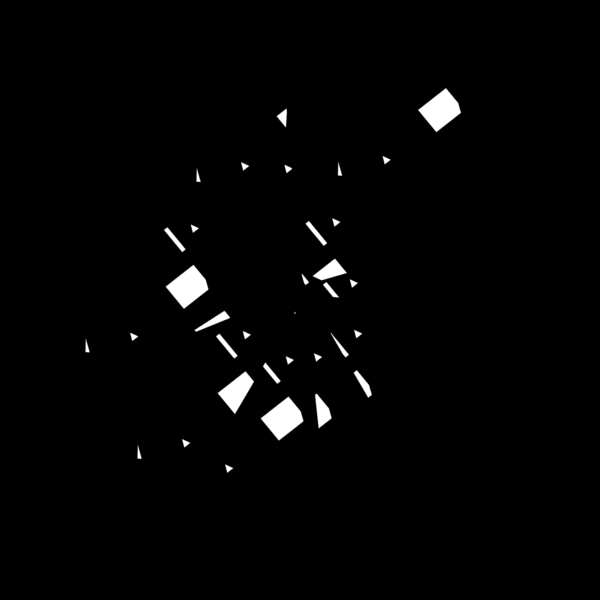}}
\\
 & & & \\ 
\large \textbf{ExParP} & & & \\
\subfloat[\scriptsize $\BGMldefault$ local (68\%)]{\includegraphics[scale=0.165]{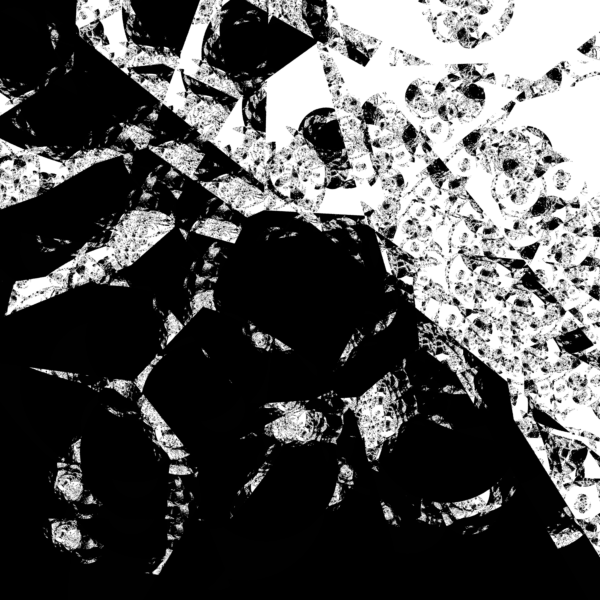}} & 
\subfloat[\scriptsize $\BGMldefault$ global (52\%)]{\includegraphics[scale=0.165]{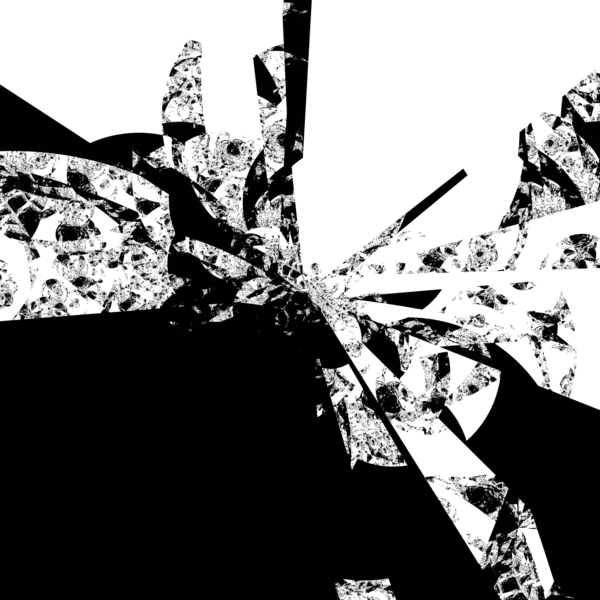}} & 
\subfloat[\scriptsize $\BGMlbest$ local (100\%)]{\includegraphics[scale=0.165]{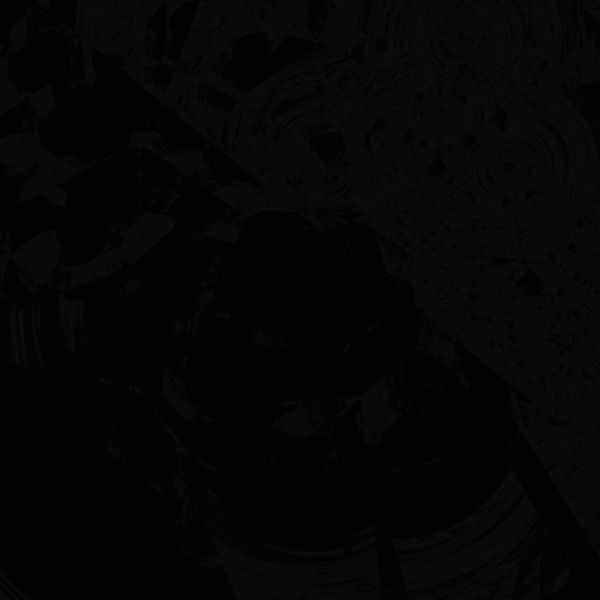}} &
\subfloat[\scriptsize $\BGMlbest$ global (100\%)]{\includegraphics[scale=0.165]{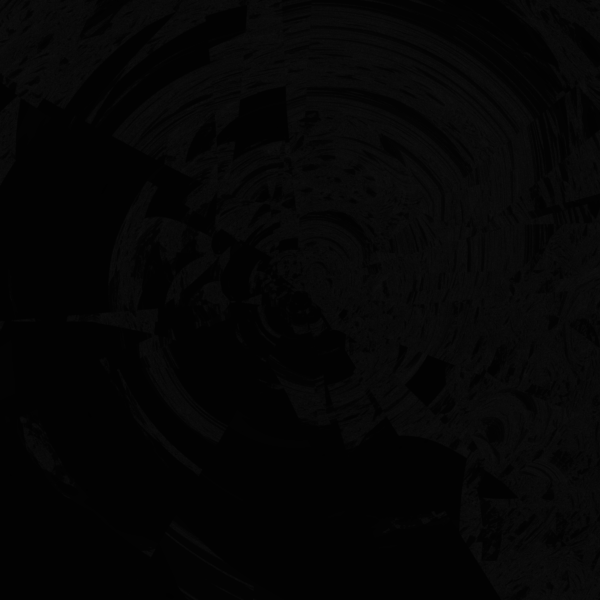}}
\\
 & & & \\
\large \textbf{DR} & & & \\
\subfloat[\scriptsize $\BGMldefault$ local (96\%)]{\includegraphics[scale=0.165]{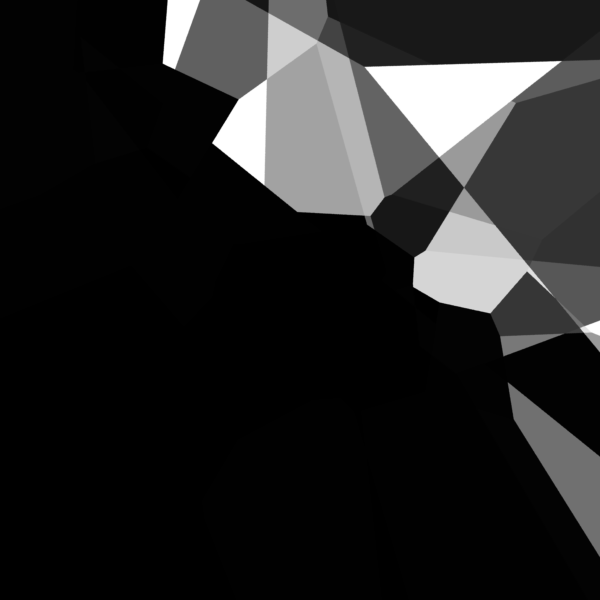}} & 
\subfloat[\scriptsize $\BGMldefault$ global (94\%)]{\includegraphics[scale=0.165]{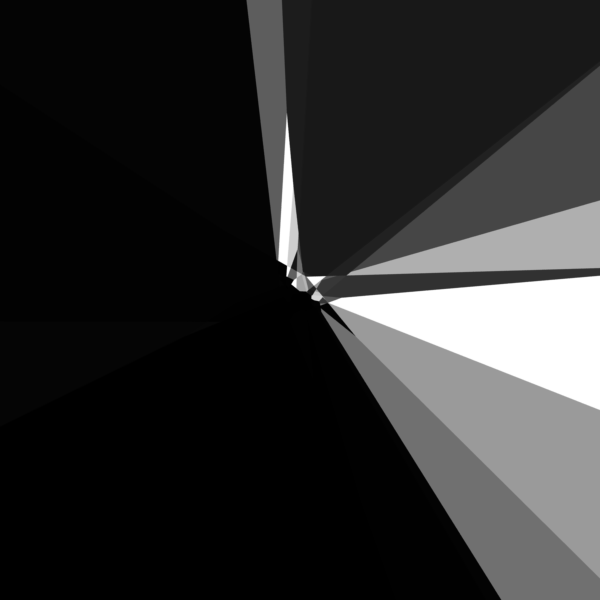}} & 
\subfloat[\scriptsize $\BGMlbest$ local (100\%)]{\includegraphics[scale=0.165]{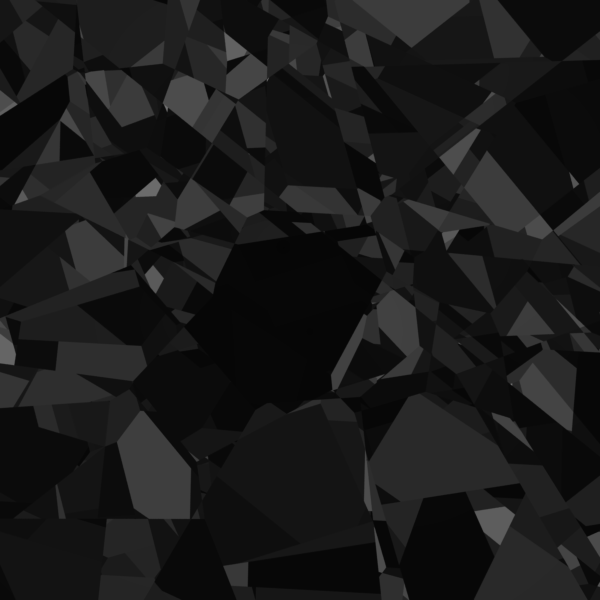}} &
\subfloat[\scriptsize $\BGMlbest$ global (100\%)]{\includegraphics[scale=0.165]{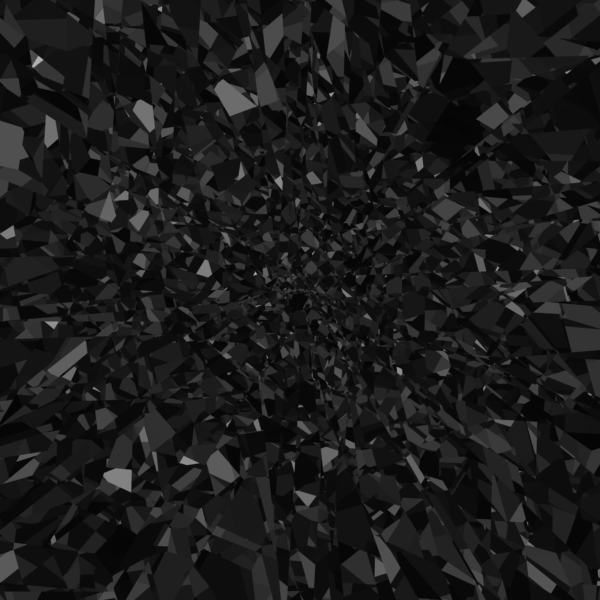}}
\\
 & & & \\
\large \textbf{CycDR} & & & \\
\subfloat[\scriptsize $\BGMldefault$ local (93\%)]{\includegraphics[scale=0.165]{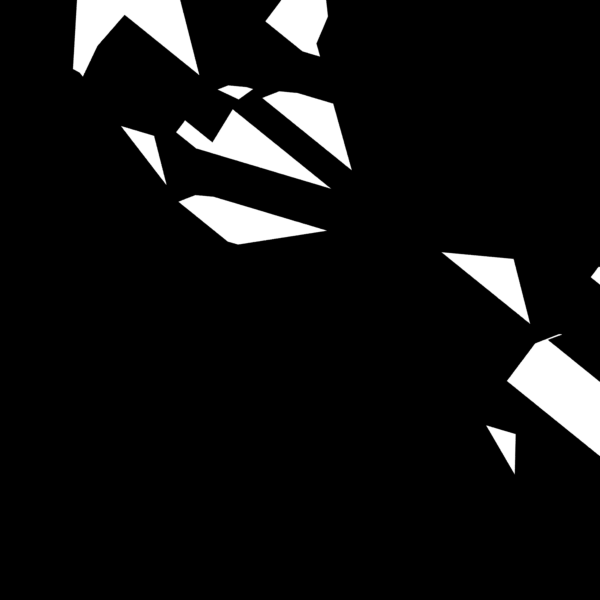}} & 
\subfloat[\scriptsize $\BGMldefault$ global (91\%)]{\includegraphics[scale=0.165]{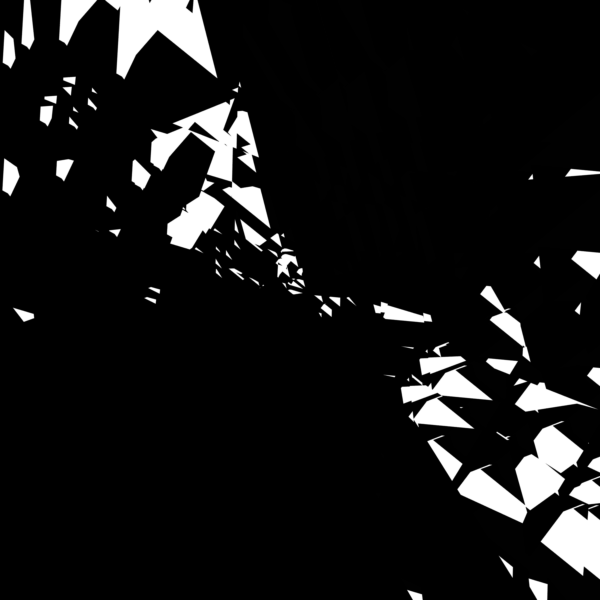}} & 
\subfloat[\scriptsize $\BGMlbest$ local (100\%)]{\includegraphics[scale=0.165]{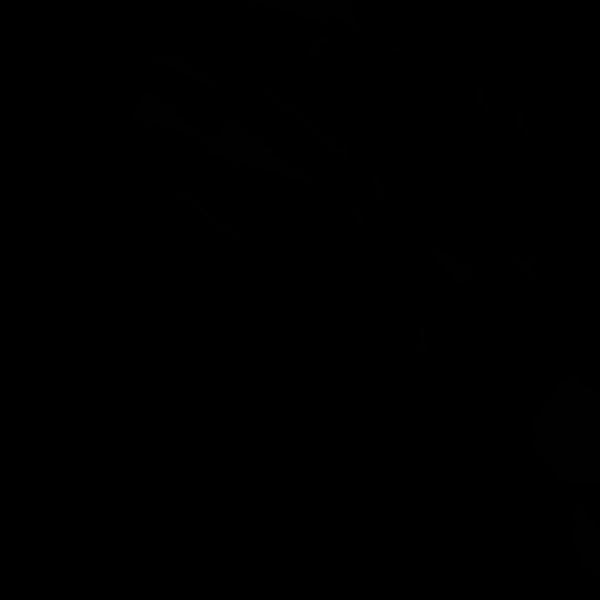}} &
\subfloat[\scriptsize $\BGMlbest$ global (100\%)]{\includegraphics[scale=0.165]{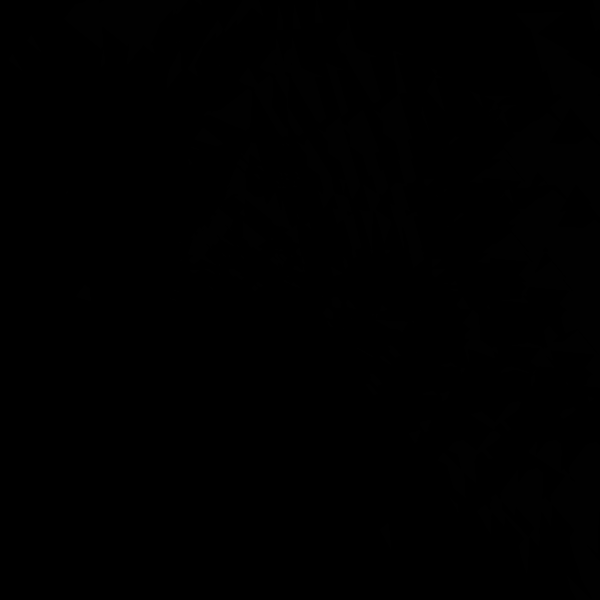}}
\end{tabular}
\caption{Behaviour of CycP, ExParP, DR, and CycDR for the few sets with few points 
constellation (success rates indicated in parentheses)}
\end{figure}

\subsection{Few sets with many points}

\BGMvspace{0.5}

\begin{figure}[H]
\begin{tabular}{cccc}
\large \textbf{CycP} & & & \\
\subfloat[\scriptsize $\BGMldefault$ local (6.8\%)]{\includegraphics[scale=0.165]{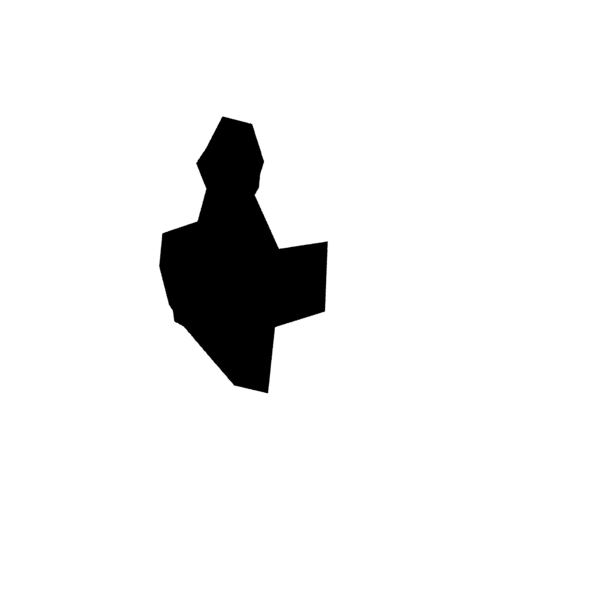}} & 
\subfloat[\scriptsize $\BGMldefault$ global (0.1\%)]{\includegraphics[scale=0.165]{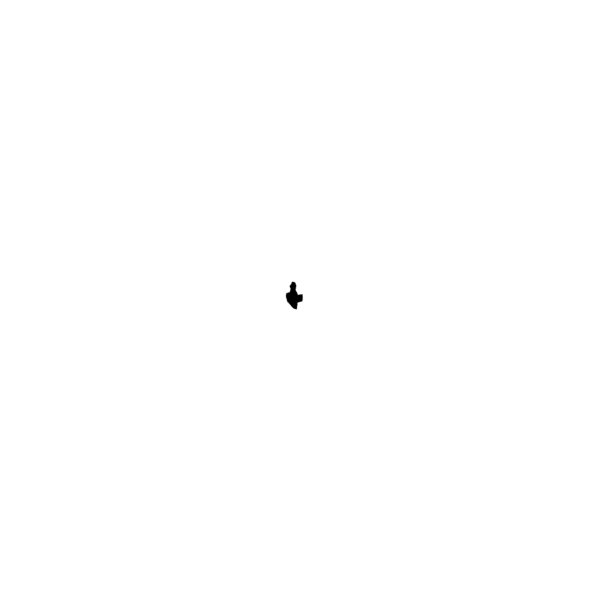}} & 
\subfloat[\scriptsize $\BGMlbest$ local (11\%)]{\includegraphics[scale=0.165]{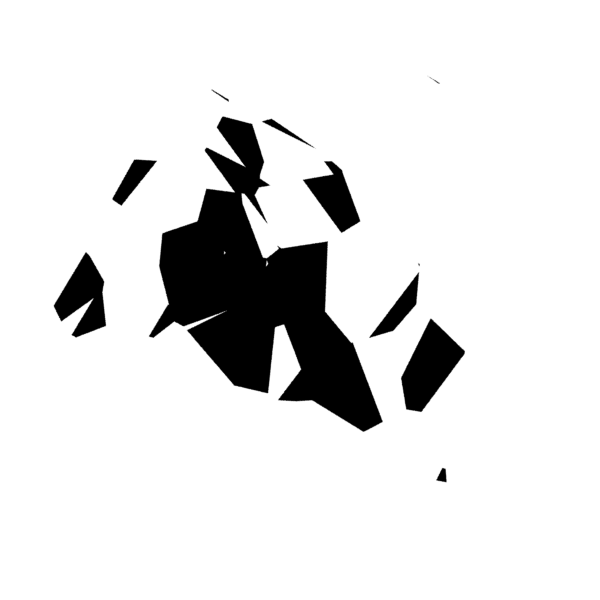}} &
\subfloat[\scriptsize $\BGMlbest$ global (12\%)]{\includegraphics[scale=0.165]{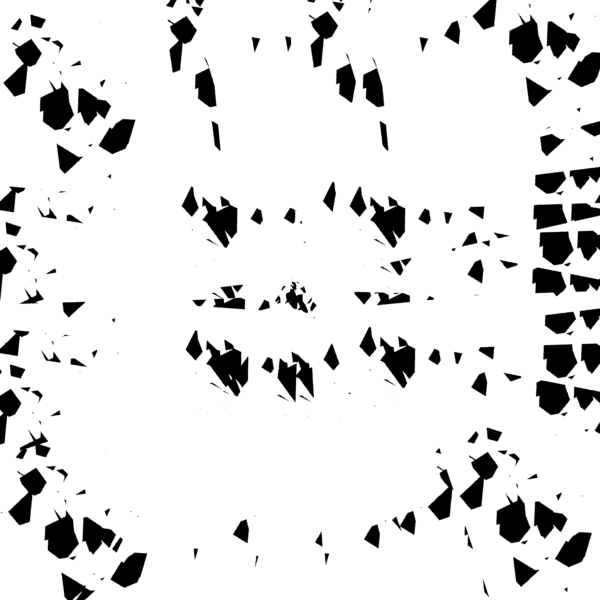}}
\\
 & & & \\ 
\large \textbf{ExParP} & & & \\
\subfloat[\scriptsize $\BGMldefault$ local (10\%)]{\includegraphics[scale=0.165]{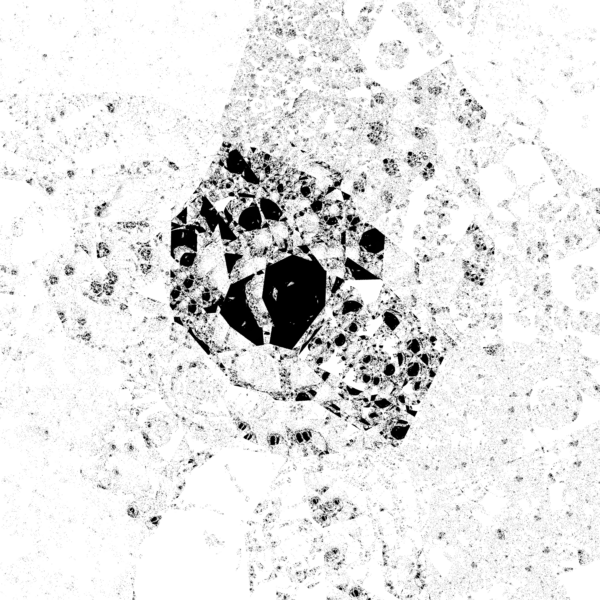}} & 
\subfloat[\scriptsize $\BGMldefault$ global (0.9\%)]{\includegraphics[scale=0.165]{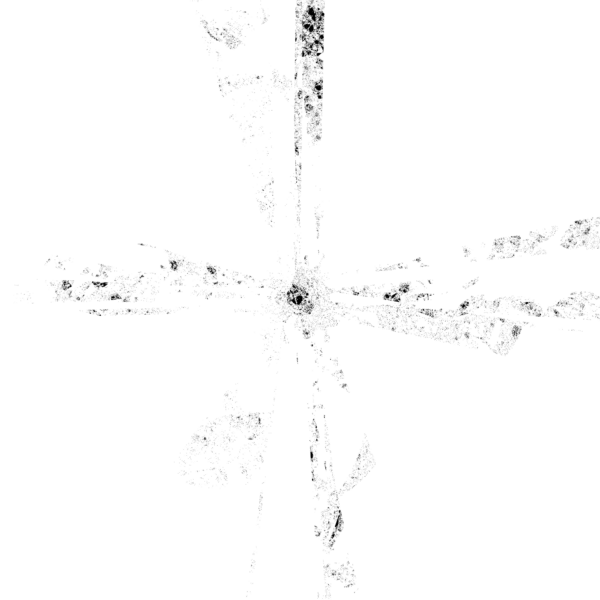}} & 
\subfloat[\scriptsize $\BGMlbest$ local (99\%)]{\includegraphics[scale=0.165]{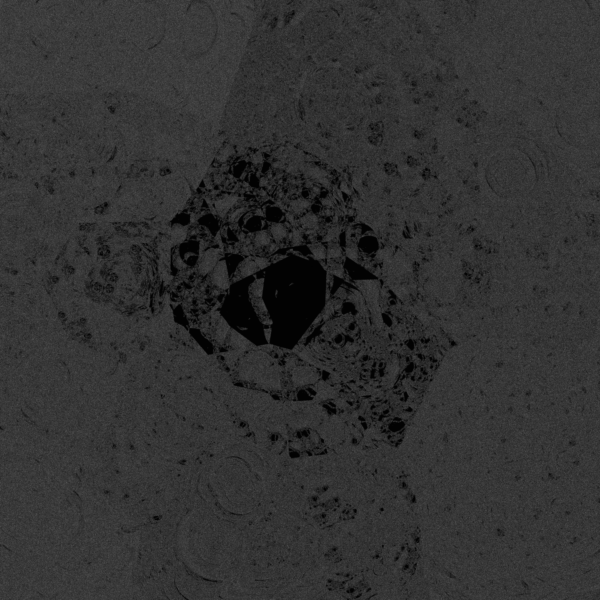}} &
\subfloat[\scriptsize $\BGMlbest$ global (99\%)]{\includegraphics[scale=0.165]{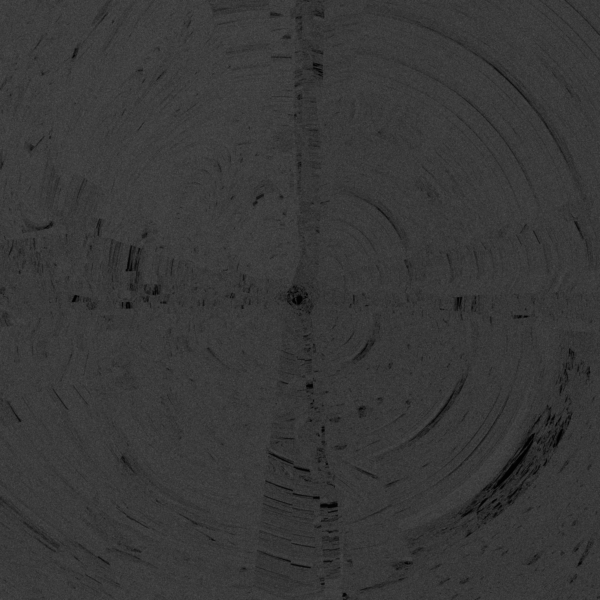}}
\\
 & & & \\
\large \textbf{DR} & & & \\
\subfloat[\scriptsize $\BGMldefault$ local (15\%)]{\includegraphics[scale=0.165]{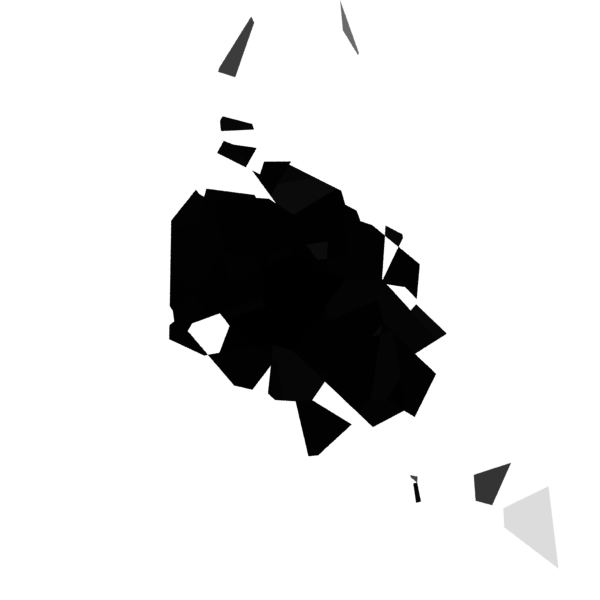}} & 
\subfloat[\scriptsize $\BGMldefault$ global (0.2\%)]{\includegraphics[scale=0.165]{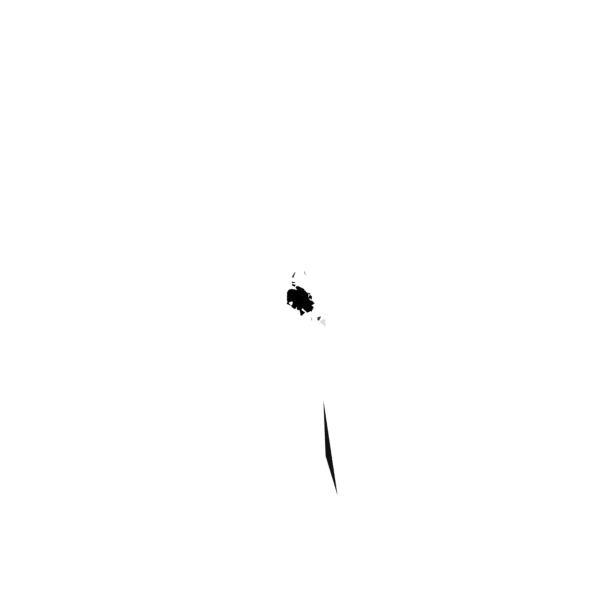}} & 
\subfloat[\scriptsize $\BGMlbest$ local (80\%)]{\includegraphics[scale=0.165]{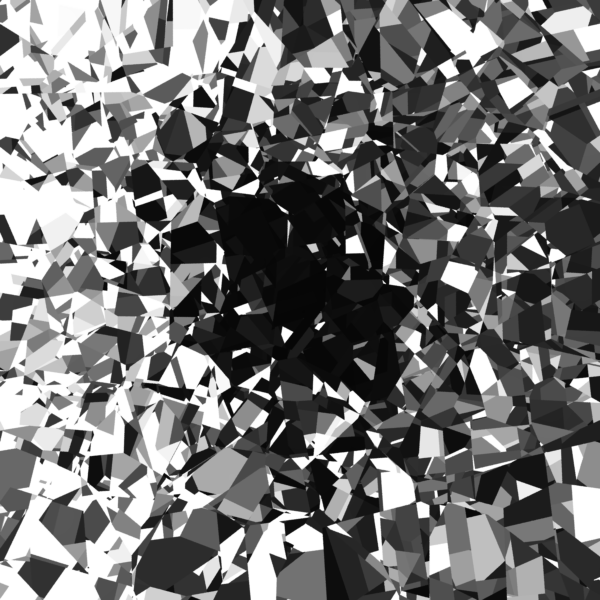}} &
\subfloat[\scriptsize $\BGMlbest$ global (81\%)]{\includegraphics[scale=0.165]{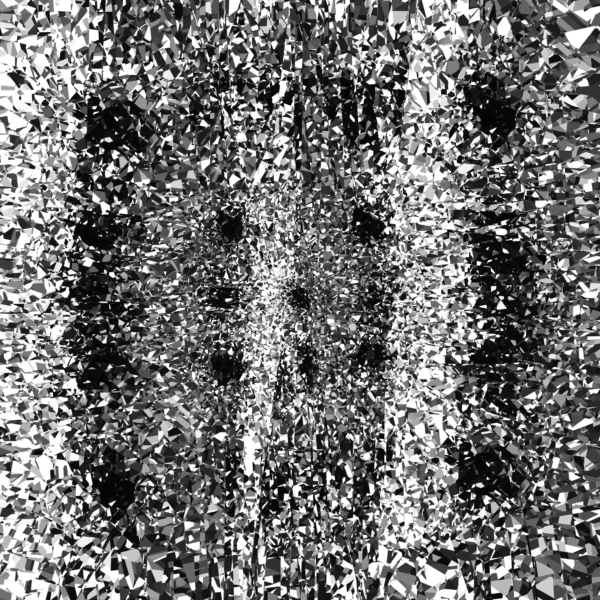}}
\\
 & & & \\
\large \textbf{CycDR} & & & \\
\subfloat[\scriptsize $\BGMldefault$ local (17\%)]{\includegraphics[scale=0.165]{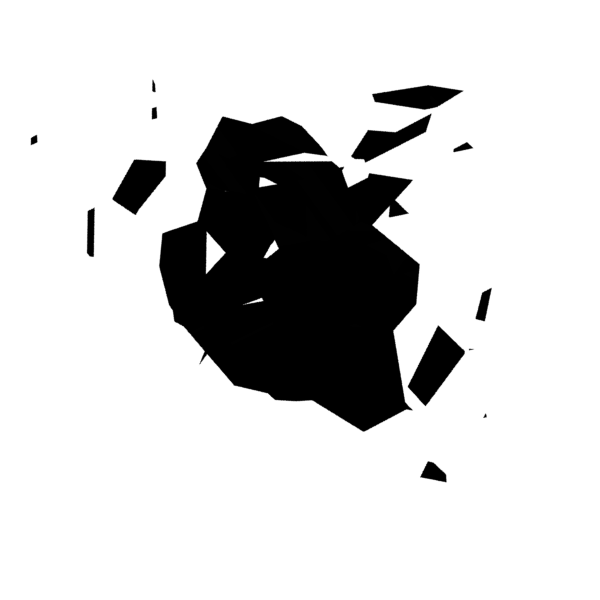}} & 
\subfloat[\scriptsize $\BGMldefault$ global (0.2\%)]{\includegraphics[scale=0.165]{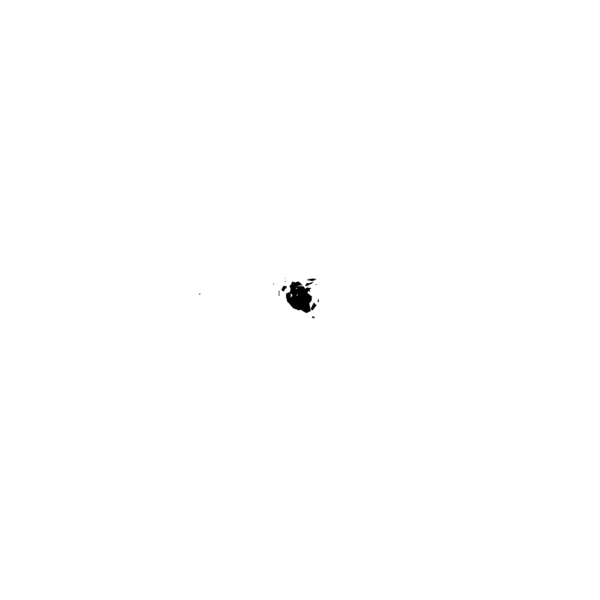}} & 
\subfloat[\scriptsize $\BGMlbest$ local (18\%)]{\includegraphics[scale=0.165]{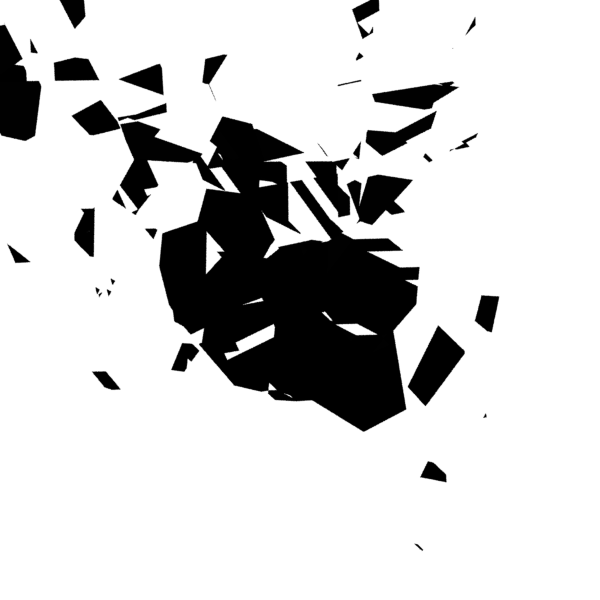}} &
\subfloat[\scriptsize $\BGMlbest$ global (4.9\%)]{\includegraphics[scale=0.165]{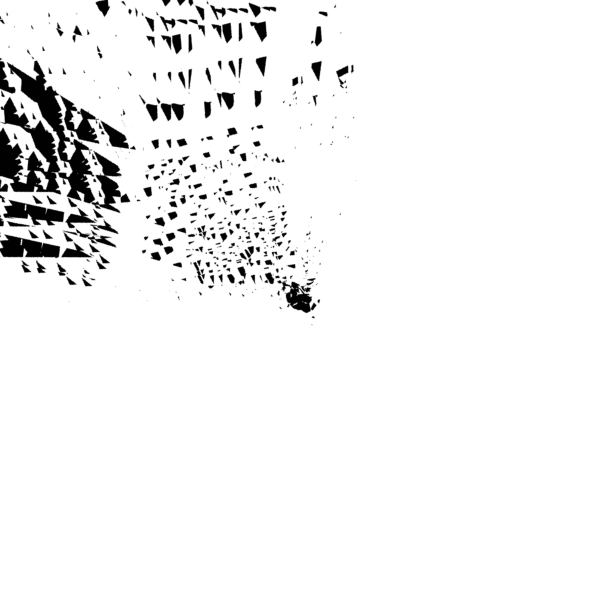}}
\end{tabular}
\caption{Behaviour of CycP, ExParP, DR, and CycDR for the few sets with many points constellation (success 
rates indicated in parentheses)}
\end{figure}

\subsection{Many sets with few points}

\BGMvspace{0.5}

\begin{figure}[H]
\begin{tabular}{cccc}
\large \textbf{CycP} & & & \\
\subfloat[\scriptsize $\BGMldefault$ local (100\%)]{\includegraphics[scale=0.165]{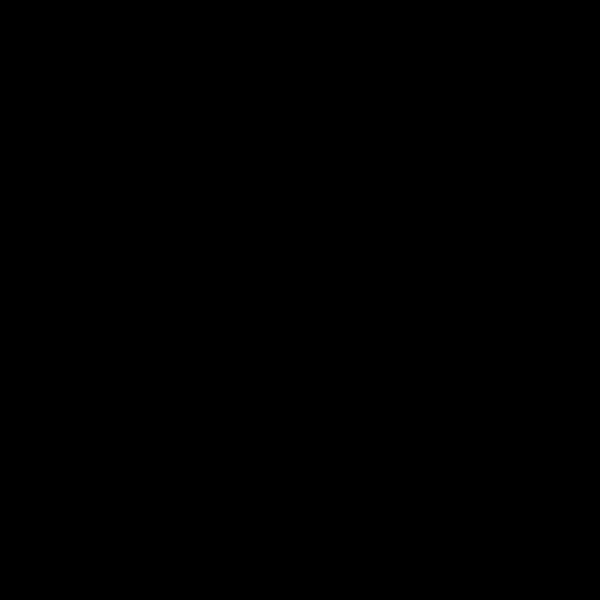}} & 
\subfloat[\scriptsize $\BGMldefault$ global (100\%)]{\includegraphics[scale=0.165]{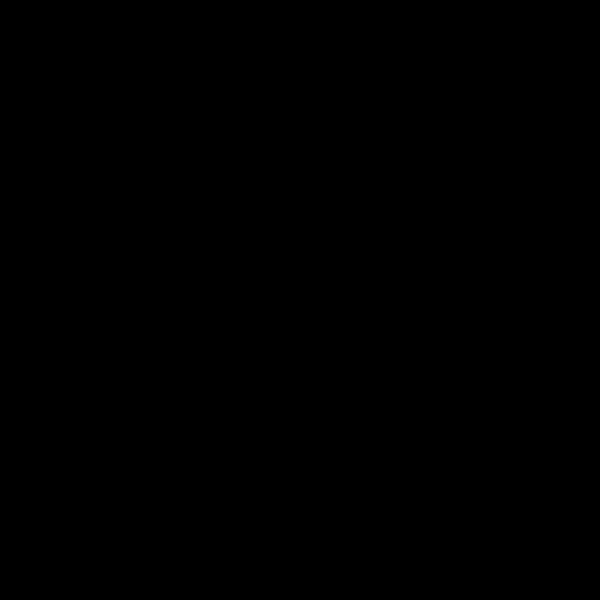}} & 
\subfloat[\scriptsize $\BGMlbest$ local (100\%)]{\includegraphics[scale=0.165]{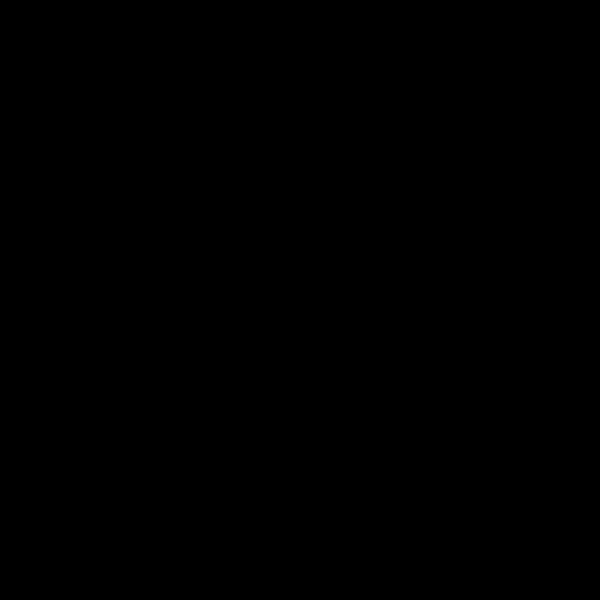}} &
\subfloat[\scriptsize $\BGMlbest$ global (100\%)]{\includegraphics[scale=0.165]{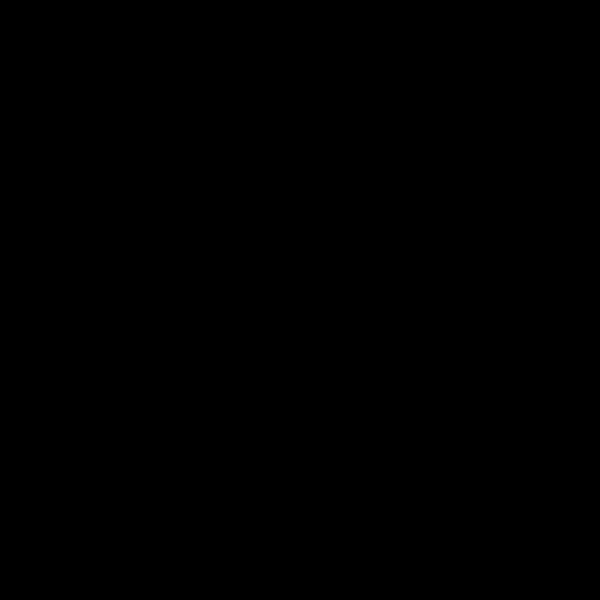}}
\\
 & & & \\ 
\large \textbf{ExParP} & & & \\
\subfloat[\scriptsize $\BGMldefault$ local (100\%)]{\includegraphics[scale=0.165]{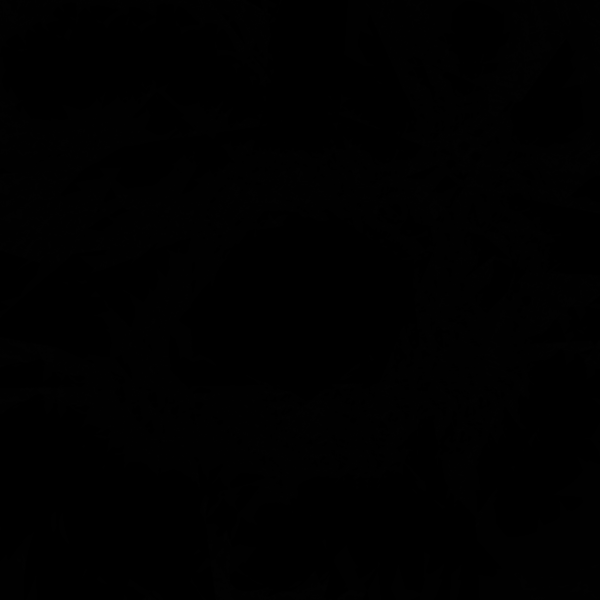}} & 
\subfloat[\scriptsize $\BGMldefault$ global (100\%)]{\includegraphics[scale=0.165]{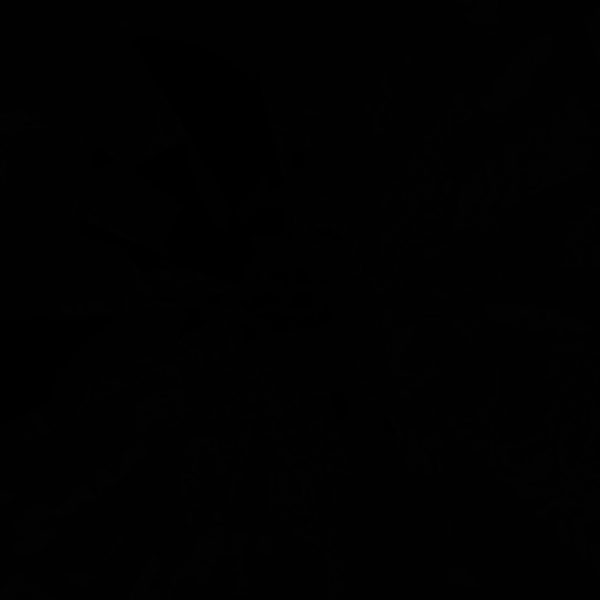}} & 
\subfloat[\scriptsize $\BGMlbest$ local (100\%)]{\includegraphics[scale=0.165]{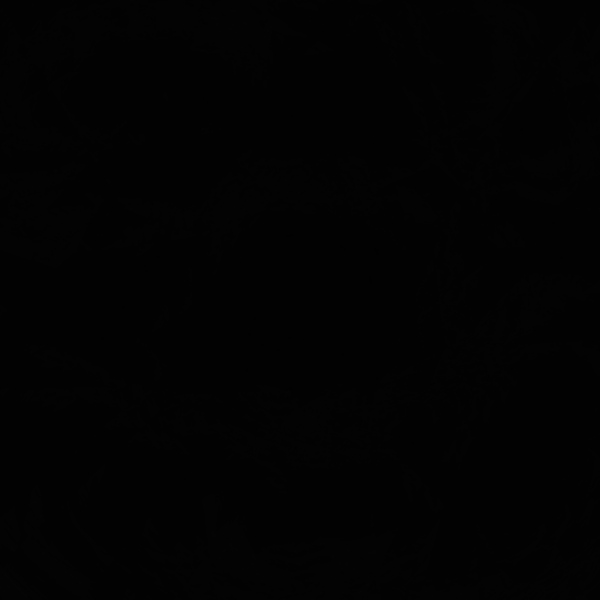}} &
\subfloat[\scriptsize $\BGMlbest$ global (100\%)]{\includegraphics[scale=0.165]{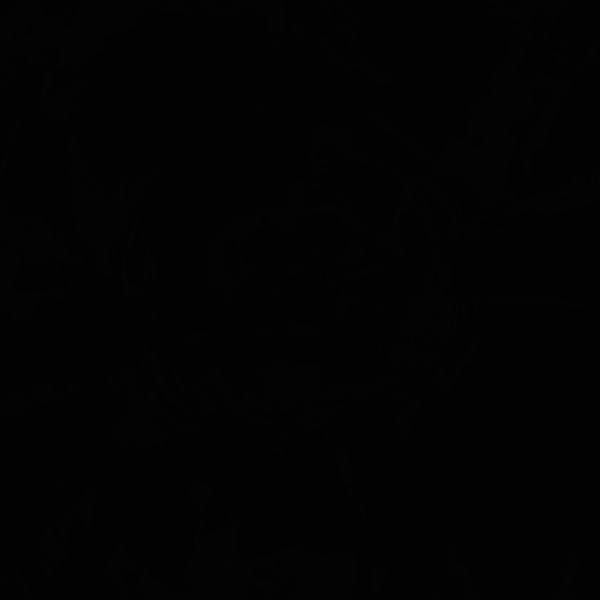}}
\\
 & & & \\
\large \textbf{DR} & & & \\
\subfloat[\scriptsize $\BGMldefault$ local (100\%)]{\includegraphics[scale=0.165]{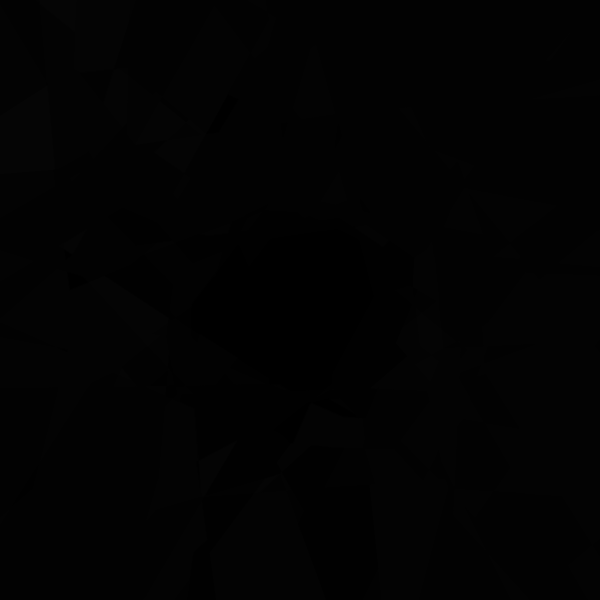}} & 
\subfloat[\scriptsize $\BGMldefault$ global (100\%)]{\includegraphics[scale=0.165]{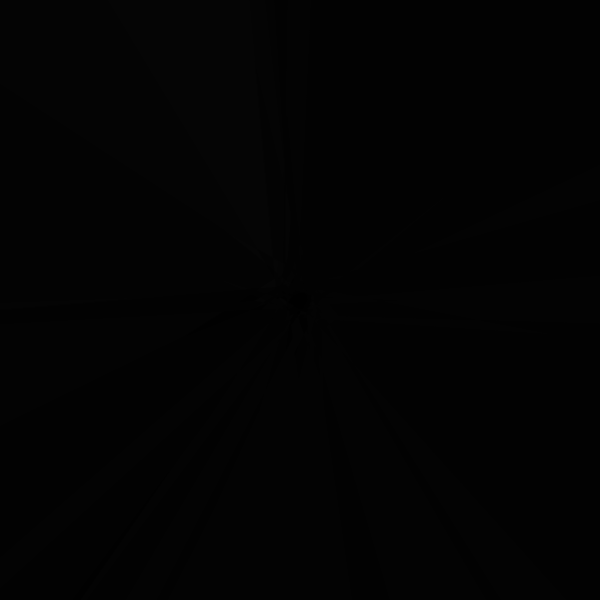}} & 
\subfloat[\scriptsize $\BGMlbest$ local (100\%)]{\includegraphics[scale=0.165]{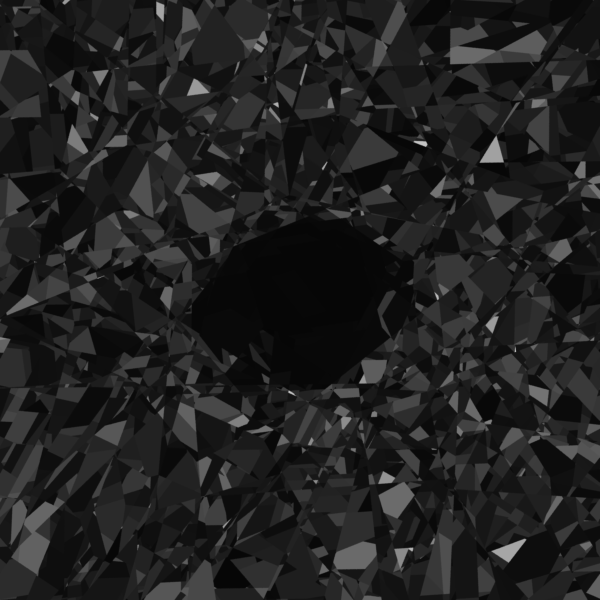}} &
\subfloat[\scriptsize $\BGMlbest$ global (100\%)]{\includegraphics[scale=0.165]{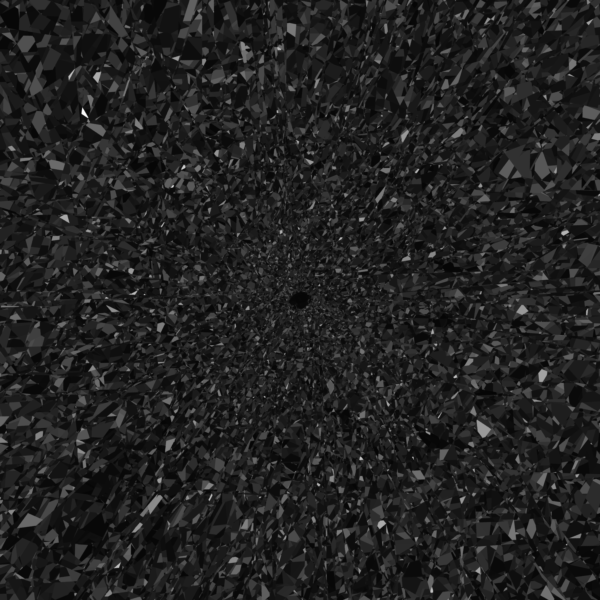}}
\\
 & & & \\
\large \textbf{CycDR} & & & \\
\subfloat[\scriptsize $\BGMldefault$ local (100\%)]{\includegraphics[scale=0.165]{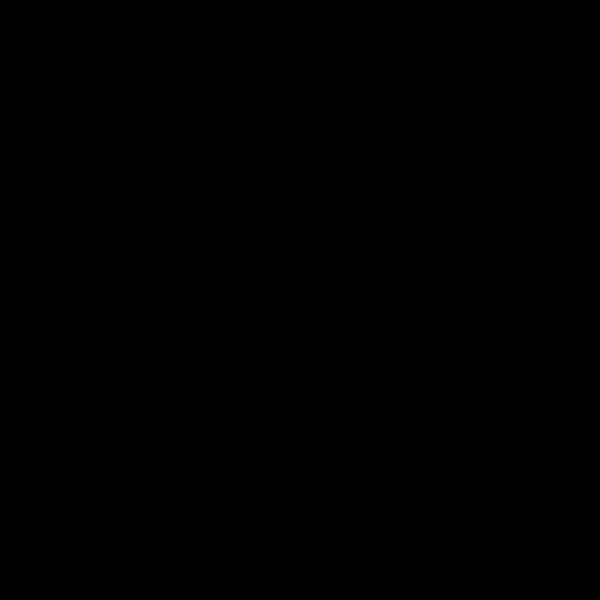}} & 
\subfloat[\scriptsize $\BGMldefault$ global (100\%)]{\includegraphics[scale=0.165]{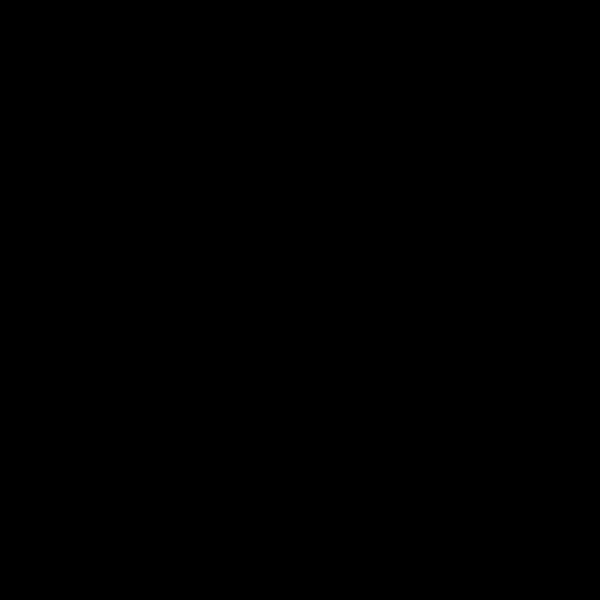}} & 
\subfloat[\scriptsize $\BGMlbest$ local (100\%)]{\includegraphics[scale=0.165]{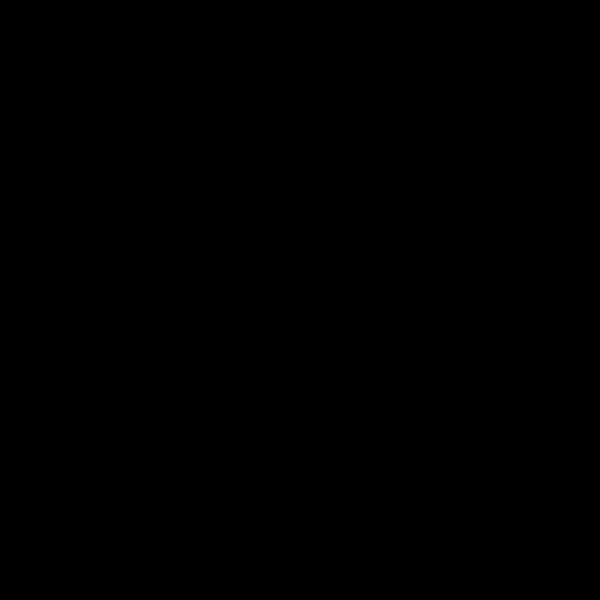}} &
\subfloat[\scriptsize $\BGMlbest$ global (100\%)]{\includegraphics[scale=0.165]{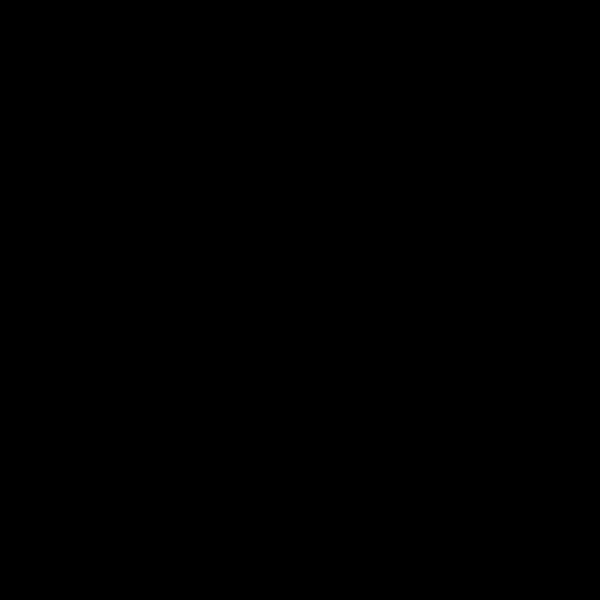}}
\end{tabular}
\caption{Behaviour of CycP, ExParP, DR, and CycDR for 
the many sets with few points constellation (success rates indicated in 
parentheses)}
\end{figure} 

\subsection{Many sets with many points}

\BGMvspace{0.5}

\begin{figure}[H]
\begin{tabular}{cccc}
\large \textbf{CycP} & & & \\
\subfloat[\scriptsize $\BGMldefault$ local (24\%)]{\includegraphics[scale=0.165]{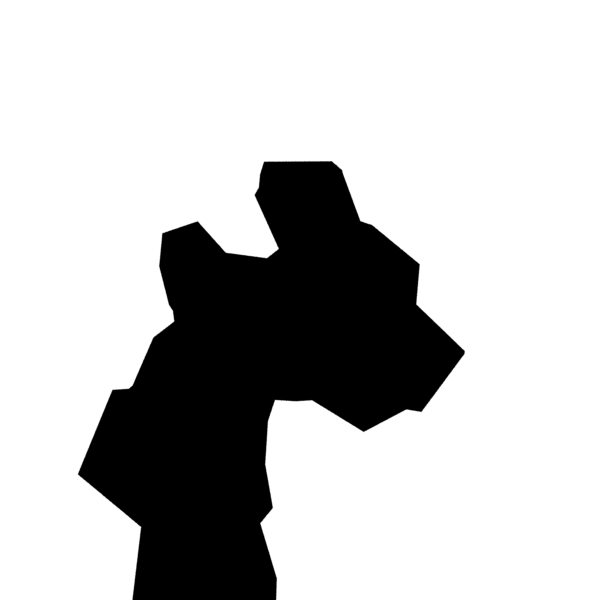}} & 
\subfloat[\scriptsize $\BGMldefault$ global (2.2\%)]{\includegraphics[scale=0.165]{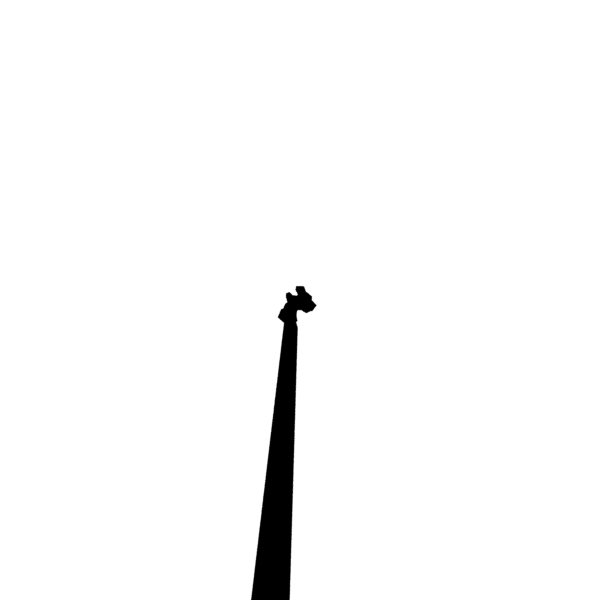}} & 
\subfloat[\scriptsize $\BGMlbest$ local (38\%)]{\includegraphics[scale=0.165]{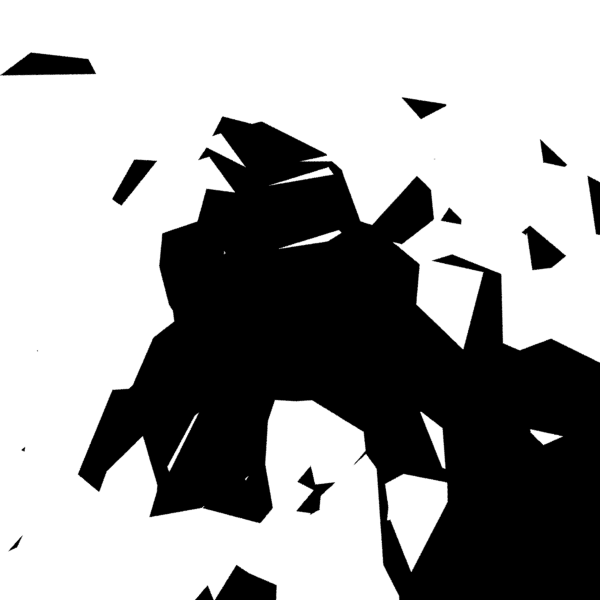}} &
\subfloat[\scriptsize $\BGMlbest$ global (47\%)]{\includegraphics[scale=0.165]{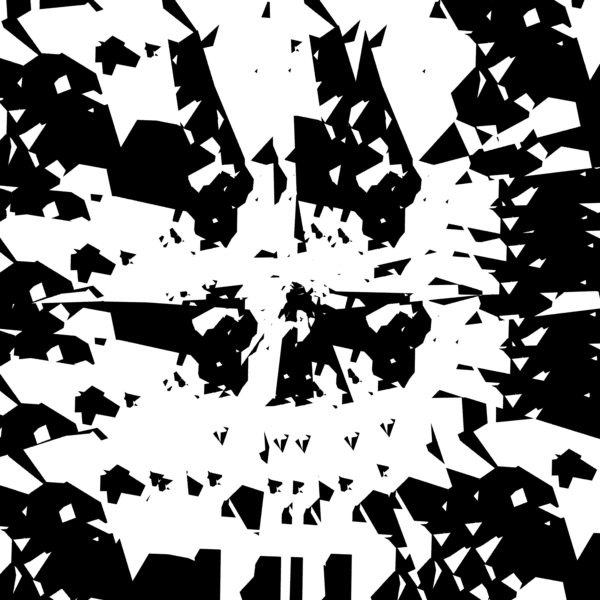}}
\\
 & & & \\ 
\large \textbf{ExParP} & & & \\
\subfloat[\scriptsize $\BGMldefault$ local (100\%)]{\includegraphics[scale=0.165]{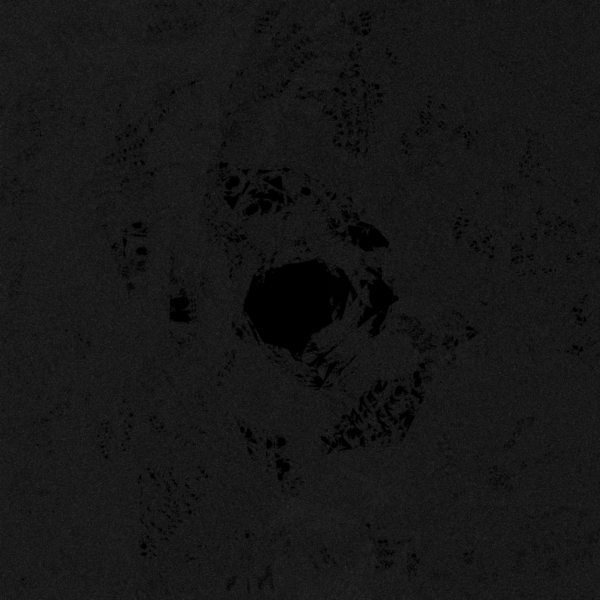}} & 
\subfloat[\scriptsize $\BGMldefault$ global (100\%)]{\includegraphics[scale=0.165]{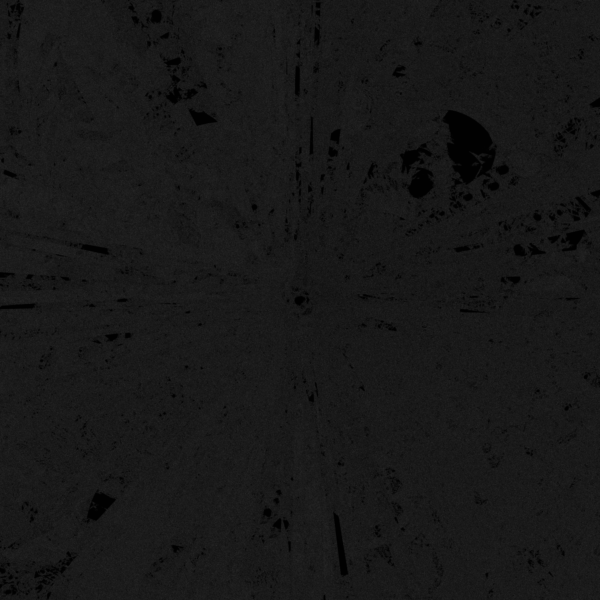}} & 
\subfloat[\scriptsize $\BGMlbest$ local (100\%)]{\includegraphics[scale=0.165]{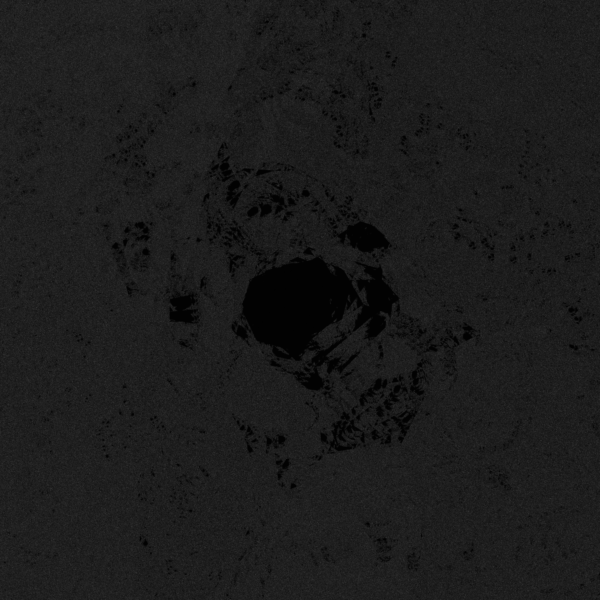}} &
\subfloat[\scriptsize $\BGMlbest$ global (100\%)]{\includegraphics[scale=0.165]{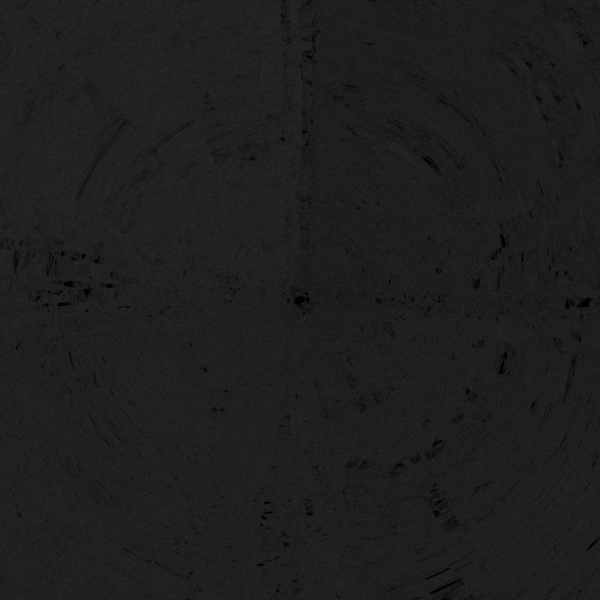}}
\\
 & & & \\
\large \textbf{DR} & & & \\
\subfloat[\scriptsize $\BGMldefault$ local (53\%)]{\includegraphics[scale=0.165]{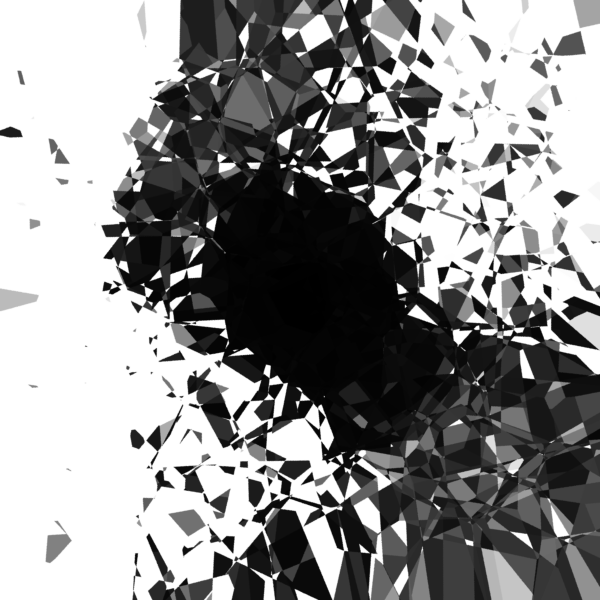}} & 
\subfloat[\scriptsize $\BGMldefault$ global (40\%)]{\includegraphics[scale=0.165]{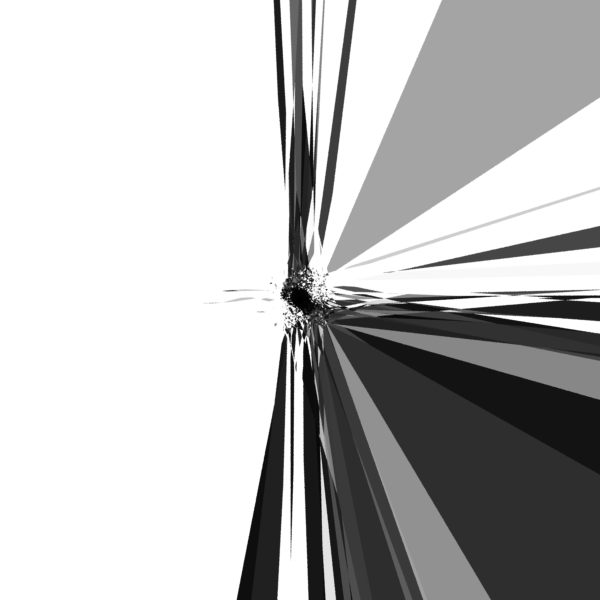}} & 
\subfloat[\scriptsize $\BGMlbest$ local (56\%)]{\includegraphics[scale=0.165]{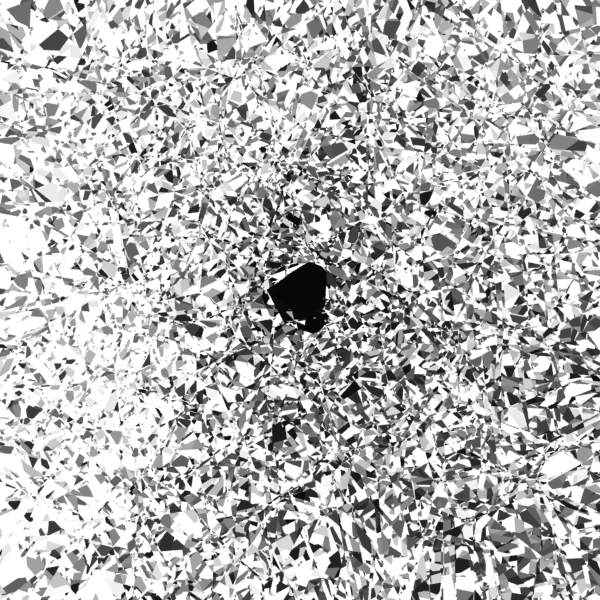}} &
\subfloat[\scriptsize $\BGMlbest$ global (57\%)]{\includegraphics[scale=0.165]{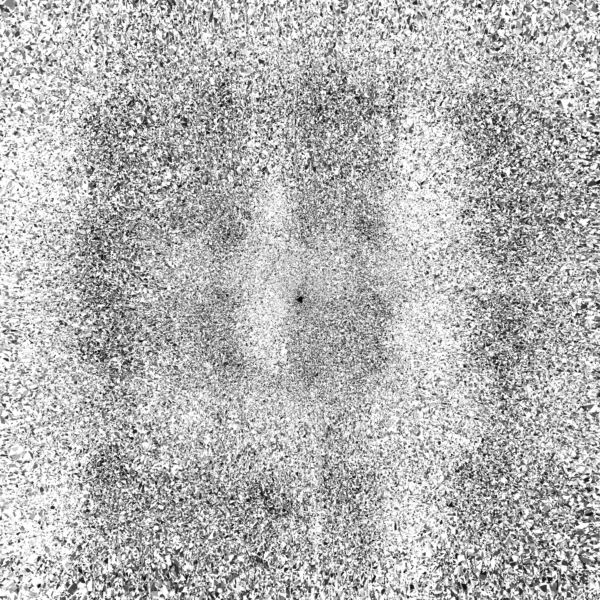}}
\\
 & & & \\
\large \textbf{CycDR} & & & \\
\subfloat[\scriptsize $\BGMldefault$ local (83\%)]{\includegraphics[scale=0.165]{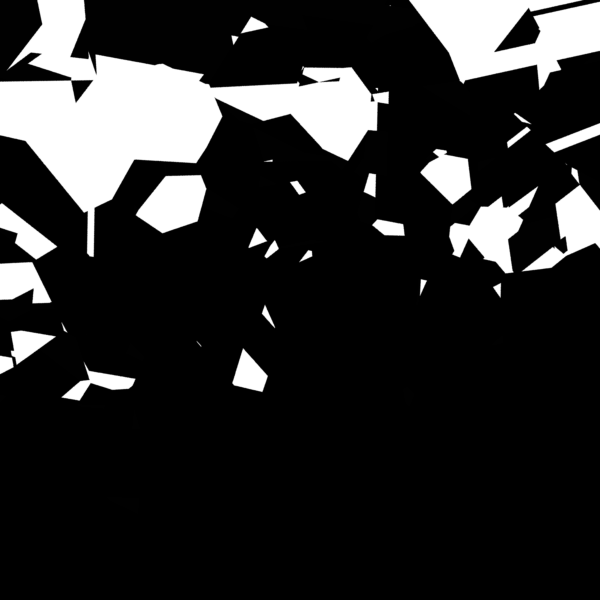}} & 
\subfloat[\scriptsize $\BGMldefault$ global (66\%)]{\includegraphics[scale=0.165]{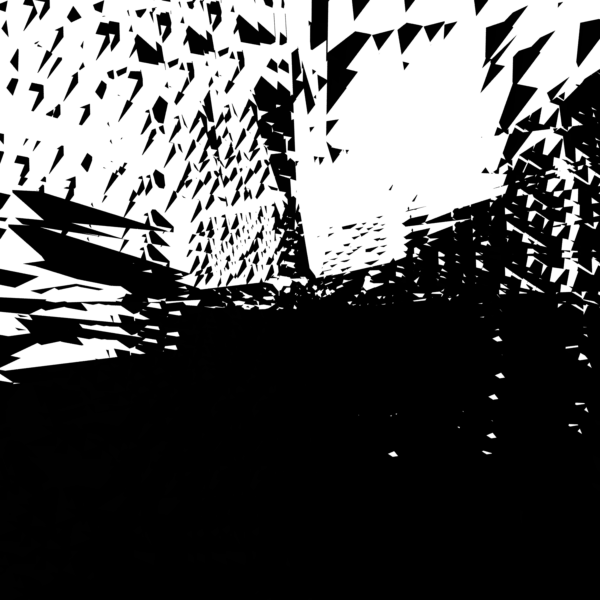}} & 
\subfloat[\scriptsize $\BGMlbest$ local (84\%)]{\includegraphics[scale=0.165]{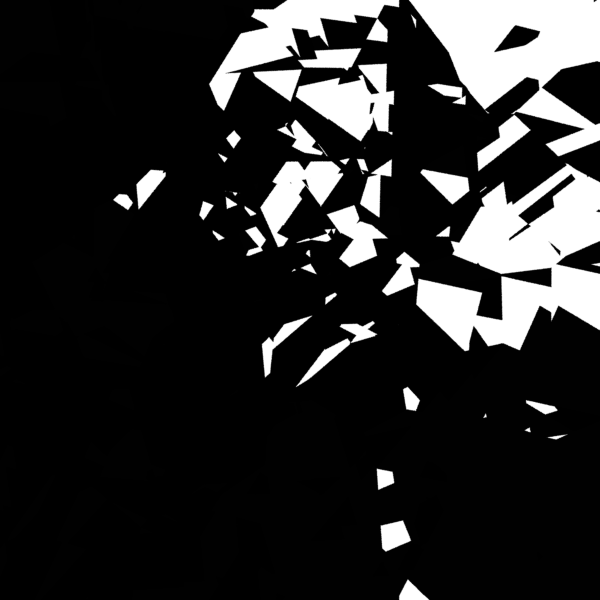}} &
\subfloat[\scriptsize $\BGMlbest$ global (82\%)]{\includegraphics[scale=0.165]{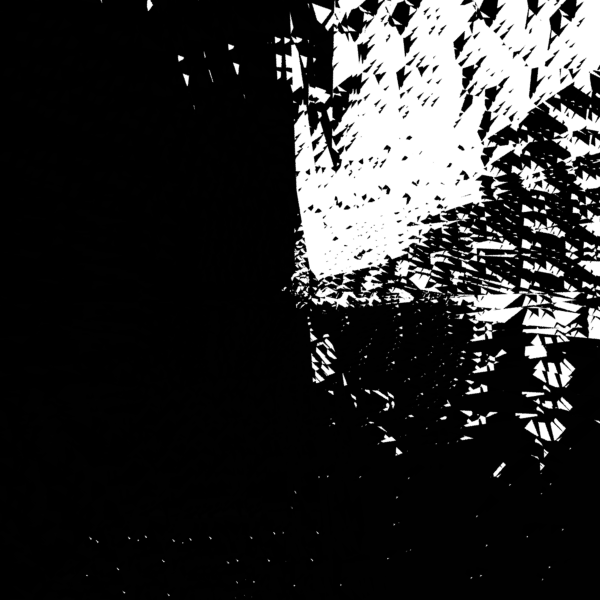}}
\end{tabular}
\caption{Behaviour of CycP, ExParP, DR, and CycDR 
for the many sets with many points constellation (success rates indicated in parentheses)}
\end{figure} 

\subsection{Discussion}
Comparing the success rates reported in the figures above,
it appears that ExParP, DR, and CycDR are good choices;
we recommend that CycP be not used. The effect of the tuning parameter
$\lambda$ is very striking for most algorithms when comparing performance
of $\BGMldefault$ with $\BGMlbest$.

\section{Divertissements}

\label{BGMsec:divert}

\begin{figure}
	\centering
		\begin{tabular}{c}
			\subfloat[]{\includegraphics[scale=0.33]{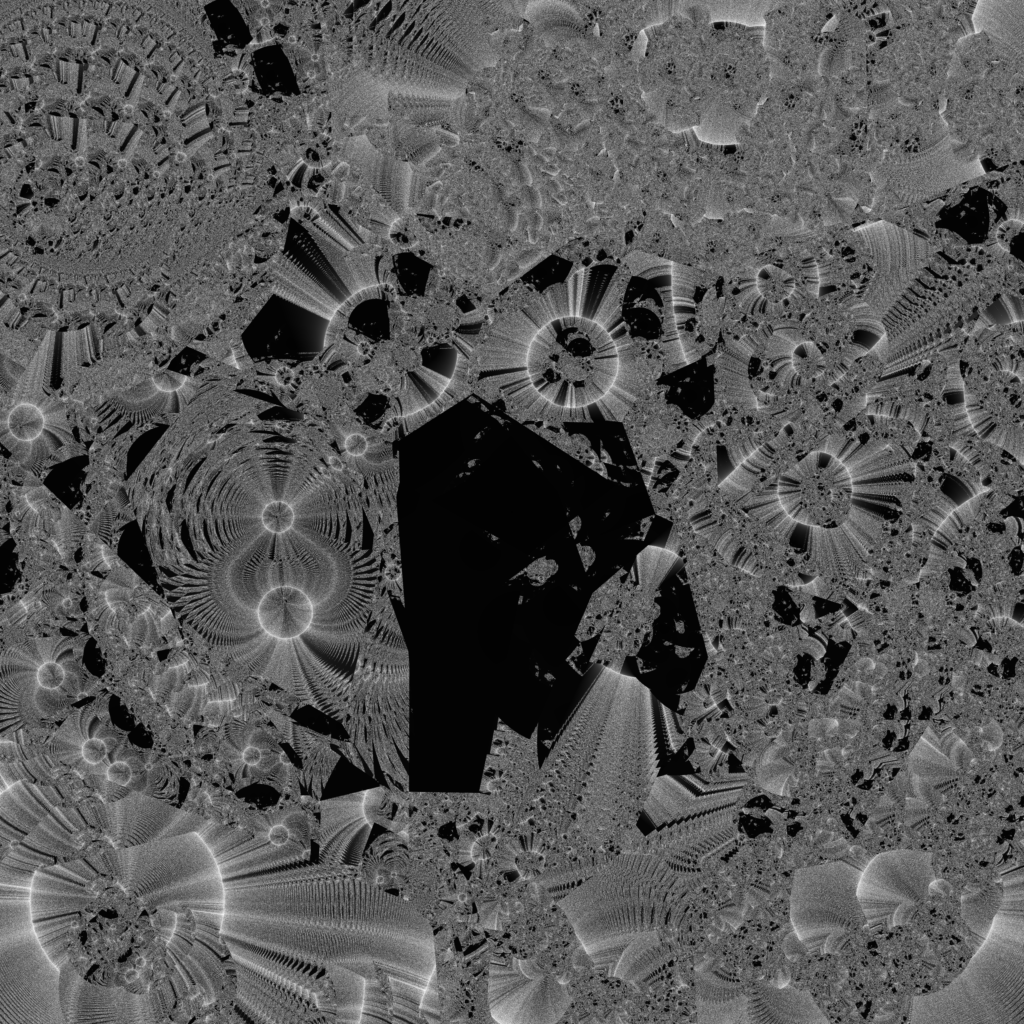}}\\
			\subfloat[]{\includegraphics[scale=0.33]{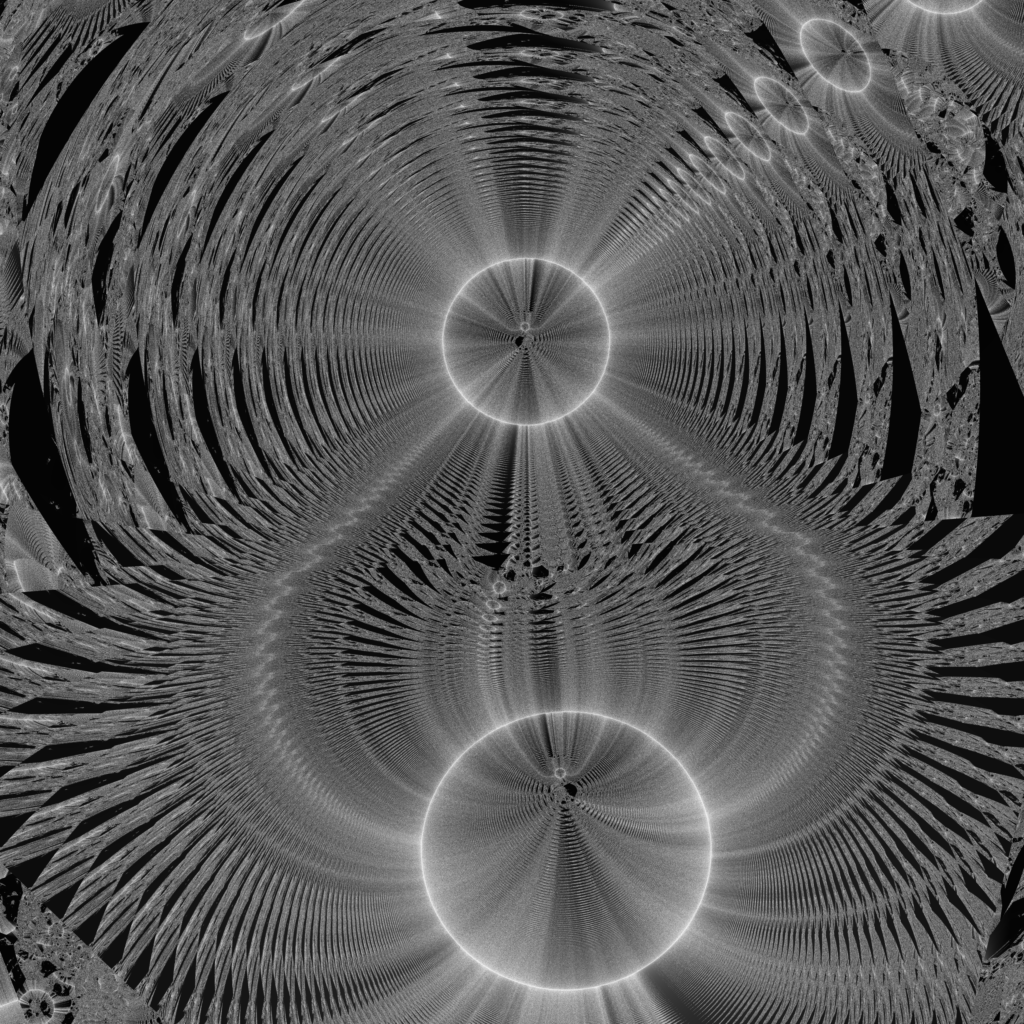}}
		\end{tabular}
		\caption{Shown in (a) is the performance of ExParP on a constellation consisting of 3 sets with 
		20 points each, with $\lambda = 0.998$, within the region $[-10,10]\times[-10,10]$. 
		A close-up of the centre-left region of (a) is presented in (b).}
		\label{fig:BGM_picwithcloseup}
	\end{figure}

\begin{figure}
	\centering
	\includegraphics [width=1.0\linewidth] {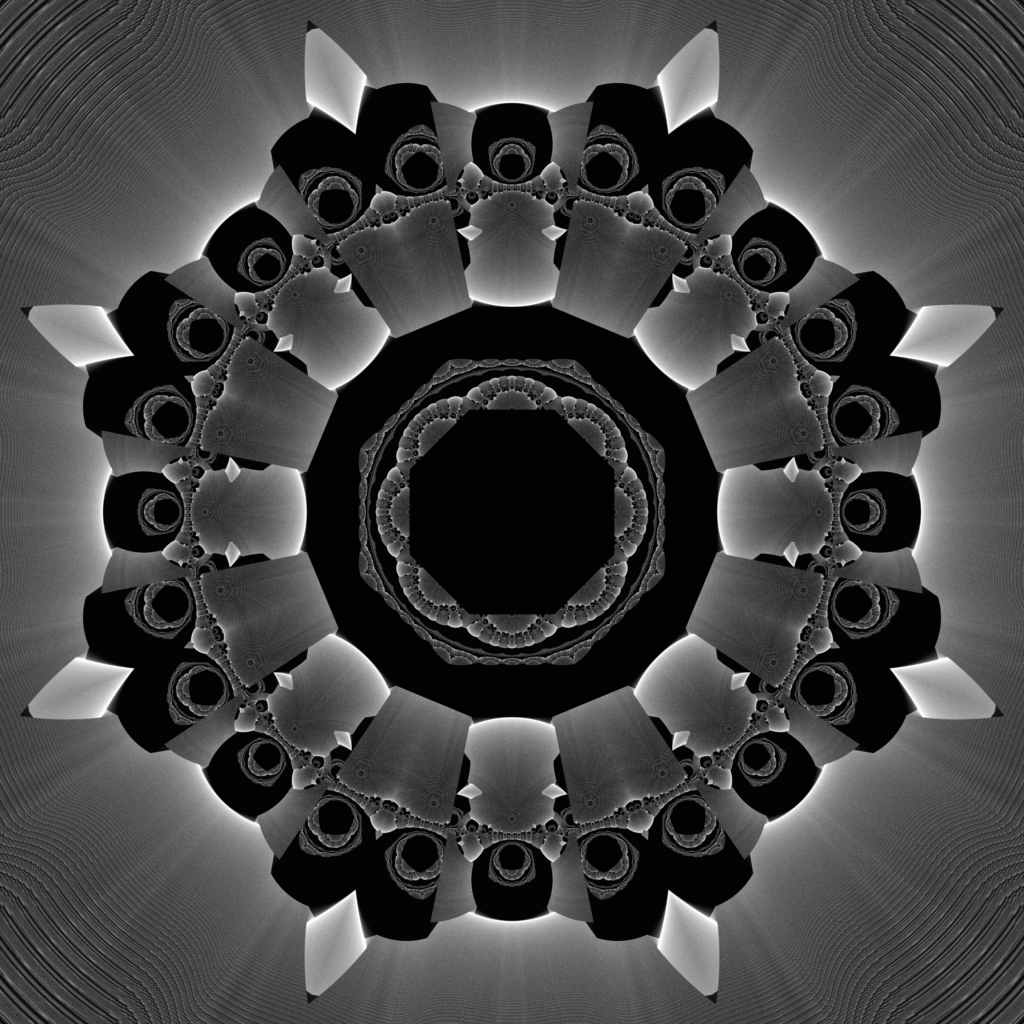}
	\caption {ExParP for a constellation with $\lambda = 0.995$, consisting of 2 subsets of concentric circles
	centred at the origin, with radii 4 and 8, containing 8 and 16 equispaced points, 
	respectively.}
		\label{fig:BGM_beauty}
\end{figure}

We experimented also with other constellations and encountered
some interesting behaviour of ExParP.
This algorithm seems to exhibit fractal-like behaviour for some constellations
--- whether they are created randomly or not. 
In the following, we present three images that we found particularly
delightful in Figure~\ref{fig:BGM_picwithcloseup} and Figure~\ref{fig:BGM_beauty}. 

\section{Concluding remarks}

\label{BGMsec:last}

We encountered a somewhat surprising complexity in the behaviour of four
algorithms for solving feasibility problems in a simple nonconvex case.
The importance of the tuning parameter $\lambda$ is apparent as
is the proximity to solutions (local vs global behaviour).
Further studies are needed to find effective guidelines for users in terms of 
choice of algorithms and the choice of the parameter $\lambda$. 
Finally, and similarly to \cite{BGMbib-Boretal}, 
we encountered \emph{beauty} in our numerical explorations.
It is our hope that others will join us and explore theoretically and
numerically this fascinating universe of constellations. 

\section*{Acknowledgements}
We thank the referee for constructive comments and suggestions.
The research of HHB was partially supported by NSERC.

\input{referenc}

\end{document}

%% file: referenc.tex
%
%
%
%